\documentclass[11pt]{article}

\usepackage{amsmath, amsthm}
\usepackage{amsmath, amsfonts}
\usepackage{amsmath, amssymb}
\usepackage{amsmath}
\usepackage{graphics}
\usepackage{graphicx}
\usepackage{color}
\usepackage{hyperref}
\usepackage[all]{xy}

\newtheorem{theorem}{Theorem}[subsection]
\newtheorem{definition}[theorem]{Definition}
\newtheorem{definition-lemma}[theorem]{Definition/Lemma}
\newtheorem{definition-explanation}[theorem]{Definition/Explanation}
\newtheorem{explanation-definition}[theorem]{Explanation/Definition}
\newtheorem{definition-fact}[theorem]{Definition/Fact}
\newtheorem{definition-notation}[theorem]{Definition/Notation}
\newtheorem{definition-conjecture}[theorem]{Definition/Conjecture}
\newtheorem{lemma}[theorem]{Lemma}
\newtheorem{lemma-definition}[theorem]{Lemma/Definition}

\newtheorem{remark}[theorem]{\it Remark}
\newtheorem{remark-notation}[theorem]{\it Remark/Notation}

\newtheorem{application-lemma}[theorem]{Application/Lemma}

\newtheorem{motif-question-project}[theorem]{Motif/Question/Project}

\newtheorem{example}[theorem]{Example}
\newtheorem{example-definition}[theorem]{Example/Definition}

\newtheorem{definition-prototype}[theorem]{Definition-Prototype}

\numberwithin{equation}{subsection}

\newtheorem{sdefinition-lemma}[stheorem]{Definition/Lemma}
\newtheorem{sdefinition-explanation}[stheorem]{Definition/Explanation}
\newtheorem{sexplanation-definition}[stheorem]{Explanation/Definition}
\newtheorem{sdefinition-fact}[stheorem]{Definition/Fact}
\newtheorem{sdefinition-notation}[stheorem]{Definition/Notation}
\newtheorem{sdefinition-conjecture}[stheorem]{Definition/Conjecture}

\newtheorem{slemma-definition}[stheorem]{Lemma/Definition}

\newtheorem{sremark-notation}[stheorem]{\it Remark/Notation}

\newtheorem{sapplication-lemma}[stheorem]{Application/Lemma}

\newtheorem{smotif-question-project}[stheorem]{Motif/Question/Project}

\newtheorem{sexample-definition}[stheorem]{Example/Definition}

\newtheorem{sdefinition-prototype}[stheorem]{Definition-Prototype}

\newtheorem{ssdefinition-lemma}[sstheorem]{Definition/Lemma}
\newtheorem{ssdefinition-explanation}[sstheorem]{Definition/Explanation}
\newtheorem{ssexplanation-definition}[sstheorem]{Explanation/Definition}
\newtheorem{ssdefinition-fact}[sstheorem]{Definition/Fact}
\newtheorem{ssdefinition-notation}[sstheorem]{Definition/Notation}
\newtheorem{ssdefinition-conjecture}[sstheorem]{Definition/Conjecture}

\newtheorem{sslemma-definition}[sstheorem]{Lemma/Definition}

\newtheorem{ssremark-notation}[sstheorem]{\it Remark/Notation}

\newtheorem{ssapplication-lemma}[sstheorem]{Application/Lemma}

\newtheorem{ssmotif-question-project}[sstheorem]{Motif/Question/Project}

\newtheorem{ssexample-definition}[sstheorem]{Example/Definition}

\newtheorem{ssdefinition-prototype}[sstheorem]{Definition-Prototype}

 \newcommand{\Azscriptsize}{{\mbox{\scriptsize\it A$\!$z}}}

\newcommand{\End}{\mbox{\it End}\,}

\newcommand{\Endsheaf}{\mbox{\it ${\cal E}\!$nd}\,}

\newcommand{\GL}{\mbox{\it GL}}

\newcommand{\Hom}{\mbox{\it Hom}\,}
\newcommand{\Homsheaf}{\mbox{\it ${\cal H}$om}\,}

\newcommand{\Imaginary}{\mbox{\it Im}\,}
\newcommand{\Image}{\mbox{\it Im}\,}

\newcommand{\Innsheaf}{\mbox{\it ${\cal I}\!$nn}\,}
\newcommand{\Int}{\mbox{\it Int}\,}

\newcommand{\Ker}{\mbox{\it Ker}\,}

\newcommand{\Mat}{\mbox{\it Mat}}

\newcommand{\Mod}{\mbox{\it Mod}\,}

\newcommand{\Qch}{\mbox{\it Qch}\,}

\newcommand{\Real}{\mbox{\it Re}\,}

\newcommand{\SL}{\mbox{\it SL}}

\newcommand{\SO}{\mbox{\it SO}\,}

\newcommand{\Spec}{\mbox{\it Spec}\,}
 \newcommand{\boldSpec}{\mbox{\it\bf Spec}\,}

\newcommand{\Supp}{\mbox{\it Supp}\,}

\newcommand{\Sym}{\mbox{\it Sym}}

\newcommand{\Tr}{\mbox{\it Tr}\,}

\newcommand{\gaugescriptsize}{\mbox{\scriptsize\it gauge}\,}

\newcommand{\nc}{{\mbox{\scriptsize\it nc}}}

\newcommand{\pr}{\mbox{\it pr}}

\newcommand{\rank}{\mbox{\it rank}\,}

\newcommand{\scriptsizesing}{{\mbox{\scriptsize\it sing}\,}}
\newcommand{\smoothscriptsize}{{\mbox{\scriptsize\it smooth}\,}}

\newcommand{\vol}{\mbox{\it vol}\,}

\newcommand{\longlongrightarrow}
 {\raisebox{.5ex}{\rule{1.2em}{.1ex}}\hspace{-1ex}\longrightarrow}    
\newcommand{\longlongleftarrow}
 {\longleftarrow\hspace{-1ex}\raisebox{.5ex}{\rule{1.2em}{.1ex}}}     
\newcommand{\longrightlongarrow}
 {\rule[.5ex]{.4em}{.1ex}\hspace{-1ex}\longrightarrow\hspace{-1ex}\rule[.5ex]{1.4em}{.1ex}}

\newcommand{\longrightaarrow}{\longrightarrow\hspace{-3ex}\longrightarrow}

\newcommand{\tinybullet}{{\raisebox{.2ex}{\tiny $\bullet$}}}							
						
\begin{document}

\enlargethispage{24cm}

\begin{titlepage}

$ $

\vspace{-1.5cm} 

\noindent\hspace{-1cm}
\parbox{6cm}{\small October 2024}\
   \hspace{7cm}\
   \parbox[t]{6cm}{\small
                arXiv:yymm.nnnnn [math.AG] \\
                D(16.1), NCS(3): \\
			    $\mbox{\hspace{1.2em}$C^\infty$ bundle settlement;}$\\
				$\mbox{\hspace{1.2em}$C^\infty$ noncommutative ringed space;}$\\
                $\mbox{\hspace{1.2em}dynamical D-brane, (kinetic) energy}$ 
				}

\vspace{2em}

\centerline{\large\bf
 A D-brane fantasy on noncommutative mirror symmetry, prelude:}
\vspace{1ex}
\centerline{\large\bf
 $C^\infty$ noncommutative ringed spaces from local noncommutative crepant}
\vspace{1ex}
\centerline{\large\bf 
 resolutions of a Calabi-Yau space with Gorenstein isolated singularities,} 
\vspace{1ex}
\centerline{\large\bf 
 dynamical D-branes thereupon, and questions beyond}
      
\vspace{3em}

\centerline{(Dedicated to the memory of {\it Joseph Polchinski},  1954-2018)}

\vspace{2.4em}

\centerline{\large
  Chien-Hao Liu   
            \hspace{1ex} and \hspace{1ex}
  Shing-Tung Yau
}

\vspace{2em}

\begin{quotation}
\centerline{\bf Abstract}
\vspace{0.3cm}

\baselineskip 12pt  
{\small
  In contrast to the world-sheet of a fundamental string, 
    the world-volume of stacked D-branes carries an Azumaya noncommutative structure ([L-Y1: Sec.\ 2] (D(1))),
	   allowing it to directly serve as a probe into noncommutative target-spaces.
  This feature leads to a D-brane fantasy:
    {\it Noncommutative Mirror Symmetry between noncommutative Calabi-Yau spaces may be realized
	        as different realizations of a supersymmetric D-brane world-volume quantum field theory,
			exactly like the string world-sheet aspect for Mirror Symmetry between (commutative) Calabi-Yau manifolds.
			}
  Driven by this fantasy, in the current notes a class of noncommutative ringed spaces
       shadowing over a $C^\infty$-manifold with corners
    are constructed from gluing local noncommutative crepant resolutions of Gorenstein isolated singularities.
  Dynamical D-branes on such noncommutative target-spaces are realized as maps/morphisms
       from an Azumaya manifold with a fundamental module with a connection $\nabla$ thereto.
  The notion of $\nabla$-adjusted kinetic energy for such a map is given via
    the basic noncommutative differential calculus developed earlier in [L-Y4] (D(11.1)).
  This provides an action functional for dynamical D-branes on such noncommutative spaces
      in parallel to the Polyakov action functional for fundamental bosonic strings on a commutative target-space.
  This sets up a basic stage to begin with for the realization of the D-brane fantasy on Noncommutative Mirror Symmetry.
  Questions beyond are sampled along the discussion.
 } 
\end{quotation}

\smallskip

\baselineskip 12pt
{\footnotesize
} 

 \bigskip

\noindent {\small MSC number 2020:  14A22, 81T30; 14E15, 57R22, 14J32
} 

\bigskip

\baselineskip 10pt
{\scriptsize
\noindent{\bf Acknowledgements.}
 We thank
 Andrew Strominger and Cumrun Vafa
   for influences to our understanding of strings, branes, and gravity.
 C.-H.L.\ thanks in addition
       Dori Bejleri,   Laurent C\^{o}t\'{e},    Thibault Decoppet, Yuriy Drozd,   Daniel~Freed,
	   Michael Hopkins, Vasily Krylov,    Stephen McKean,     Sunghyuk Park, Mihnea Popa,    Ming Hao Quek,
	   Matthew Reece, Subir~Sachdev,   Andrew Strominger,    Freid Tong,    Cumrun Vafa, Fan Ye
   for various topic courses,  fall 2022 - spring 2025, and
      (also to participants/speakers in Freed Group Meeting/Seminar  and Vafa Group Meeting/Seminar for)	 
        intellectual enrichment and discussions with some of them;
		Nancy Kanwisher,
            Gerard J.\ Tortora \& Bryan Derrickson,
			Philip~G.$\hspace{.7ex}$Zimbardo (1933-2024) \& Richard J.\ Gerrig ::		
		Jacob Collier,   Percy Goetschius (1853-1943), Daniel J.\ Levitin ::
		Jeong-hyun Seok ::
	    Deguchi Jin, Easy German Team, Jennifer Forrest,
	        Nguy$\tilde{\mbox{\^{e}}}$n Th\d{i} Thu H$\grave{\mbox{\u{a}}}$ng :: 	
	    Pema Ch\"{o}dr\"{o}n, Satya Narayana Goenka (1924-2013)
   for their books, masterclasses, open course, program, and/or lecture/lecture series that provide another kind of intellectual \&
     inner enrichment during the brewing years of the current work and new directions beyond;
	   Tale of Fantasy and Liyuan Yuan
   for their respective covers of Star Sky composed by Two Steps from Hell
         that accompany the typing/editing of the current notes;
  Pei-Jung Chen
      for the biweekly communications on work of J.S.\ Bach (1685-1750) and Ernesto K\"{o}hler (1849-1907)
      and
  Ling-Miao Chou
      for regular joint readings of  Tortora \& Derrickson,
   	       comments that improve the illustration, and the tremendous moral support.
} 
 
\end{titlepage}

\newpage

\enlargethispage{24cm}
\begin{titlepage}

$ $


\centerline{\small\it Dedicated to the memory of}
\centerline{\small\it Professor {\sl Joseph Polchinski, 1954-2018},}
\centerline{\small\it whose work motivated the D-project}

\vspace{3em}

\baselineskip 11pt

 \noindent
{\scriptsize 
 (From C.-H.L.)\hspace{1em}
 The time and place is  summer 1996 in a lecture hall of the Department of Physics, University of Colorado at Boulder,
   and the event is TASI Summer School,  well-known in the community of graduate students of physics.
 Through a friend I arranged a summer stay there  and audited all the lectures.
 Before that I had just learned about D-branes from my then-physics-advisor Prof.\ Orlando Alvarez at University of Miami and
     was very fascinated about them, though at the same time filled with puzzles and questions
	      that came naturally from a mathematical mind and needed to be resolved.
 It was in that summer school  I met Prof.\ Polchinski for the first time.
 He had a very special habit of closing his eyes whenever his lecture came to a high point as if he were completely immersed in
   an inner physics world he was just about to reveal  and explain to his audience.
 Not to lose such a wonderful opportunity,  I asked him several questions after his lectures.
 To my surprise and against original worries,
   he didn't kick me out of the lecture hall, thinking that this is a guy who came to embarrass him.
 Rather, to every puzzle on D-branes I brought to him, he acknowledged the issue but explained immediately
   how one can infer from a concrete, mathematically solid, physics example parallel to but simpler than the case
      he wanted to address and made conclusions and hence bypass all these mathematical obstacles.
 To him, physics looks almost like an art: an art of inference from the mathematical solid to whatever beyond.
 (Before him, I had seen Prof.\ Neil$\hspace{.7ex}$Turok
                                  doing the same in Cosmology course at Princeton University and
        Prof.\ Alvarez in the Quantum Field Theory course at University of California at Berkeley.
     And after him, Jacques Distler in String Theory Seminar at University of Texas at Austin  and
	   Shiraz Minwalla in the Quantum Field Theory course and String Theory course at Harvard University.)
 That's pretty amazing from a mathematician's eyes.
 (Based on cognitive neuroscience (cf.\ [K]),
      there could be fine differences in how each person's neural network in the brain is configured;
     thus such skill  is likely congenital,  unlearnable nor reachable from a fundamentally-mathematical-though-physics-oriented brain.)
     
 My second contact with him was in spring 2016.
 During the two decades in between,   a train of  luck had been given to me (by fate?, by destiny?, or by God/Higher Existence?)
  to promote my understanding of string theory both from the physics side and from the mathematics side
  and I had been in the middle of [L-Y'16] (D(13.1)) by then and once again had dwelled upon his work ([P'96], [P'98])
  as the source of inspiration, trying to resolve a sign puzzle,
  until I decided why not just consult him by an email.
 His clarification of the issue came soon after my email to him,
  ending with an apology (!!!) that he was in a medical treatment and the situation didn't allow him to elaborate more.
 I replied with thanks and a good wish to him to recover soon.
 Immediately afterwards, the demanding project continued to take my full focus.
    
 Then came the COVID-19 lockdown spring 2020 - summer 2022,
    which gave me an unexpected chance and period to explore -- and realized for the first time --
	how the world of education had been changed via the internet.
 Thanks to the joint efforts of many individuals and institutes, one can learn almost everything through the YouTube platform
     -- from how to fix a faucet in a house to high energy physics
	     or from a more popular language, e.g.\ German, to a rare ancient language, e.g.\ Pali -- all free.
 Among the courses I completed listening to is [Q].
 Curious about what the instructor had worked on, I made a search in the arXiv.
 It is then I came across [B-Q-W], from which I learned about the passing of Prof.\ Polchinski in 2018.
 A sense of loss sprang up unbidden despite only  limited contacts with him.
 From his own account [P'17], he should be in the first half year of his brain tumor treatment  when I emailed him.
 That must be a particularly painful period of time for him to accept the illness and adapt to the reality.
 What a generous person he must have been to still respond to my (must-be-very-trivial) question to him
  in such a life-threatening medical situation and for a stranger.
  
 The D-project was motivated by the intent to understand what dynamical D-branes really are when they stack
    as described in his TASI lecture [P'96] and by-now-classical textbooks [P'98].
 Over the years, the project has grown to addressing the question of
    what dynamical D-branes imply or can do based on this understanding.
 His work on D-branes has been in my constant rumination for nearly three decades
   -- including the first eleven years of completely quiet brewing before [L-Y'07] (D(1)) --
   and in this sense he has been my another physics advisor as well, though only through his work.
 %
 
 This prelude to a new direction of the D-project is thus dedicated to the memory of Prof.\ Polchinski.
    
\bigskip

\noindent
\parbox[t]{5em}{[B-Q-W]}
 \parbox[t]{48.2em}{R.\ Bousso,  F.\ Quevedo, and S.\ Weinberg,
   {\it Joseph Polchinski: a  biographical memoir},
    arXiv:2002.02371 [physics.hist-ph].
       }  
	
\medskip
	
\noindent
\parbox[t]{5em}{[K]} \parbox[t]{48.2em}{N.\ Kanwisher,
    {\sl The human brain},
   MIT 9.13, undergraduate course, spring 2019; MIT OpenCourseWare, accessible via YouTube.
   } 

\medskip

\noindent
\parbox[t]{5em}{[L-Y'07]} \parbox[t]{48.2em}{C.-H.\ Liu and S.-T.\ Yau,
 {\it Azumaya-type noncommutative space and morphisms therefrom:
    Polchinski's D-branes in string theory from Grothendieck's viewpoint},
 arXiv:0709.1515 [math.AG]. (D(1))	
   } 
  
\medskip
  
\noindent
\parbox[t]{5em}{[L-Y'16]} \parbox[t]{48.2em}{--------,
    {\it Dynamics of D-branes I,
	        The non-Abelian Dirac-Born-Infeld action, its first variation, and the equations of motion for D-branes
			--- with remarks on the non-Abelian Chern-Simons/Wess-Zumino term},
	  arXiv:1606.08529 [hep-th] (D(13.1)).
   } 
   
\medskip

\noindent
\parbox[t]{5em}{[P'96]} \parbox[t]{48.2em}{J.\ Polchinski,
 {\it Lectures on D-branes},
   in ``{\sl Fields, strings, and duality}", TASI 1996 Summer School,
   Boulder, Colorado, C.\ Efthimiou and B.\ Greene eds.,
   World Scientific, 1997.
 (arXiv:hep-th/9611050)
       } 

\medskip

\noindent
\parbox[t]{5em}{[P'98]} \parbox[t]{48.2em}{--------,
 {\sl String theory},
 vol.\ I$\,$: {\sl An introduction to the bosonic string};
 vol.\ II$\,$: {\sl Superstring theory and beyond},
 Cambridge Univ.\ Press, 1998.
      } 

 \medskip

 \noindent
 \parbox[t]{5em}{[P'17]} \parbox[t]{48.2em}{--------,
   {\it Memories of a theoretical physicist},
   arXiv:1708.09093 [physics.hist-ph].
   } 
 
\medskip

\noindent
\parbox[t]{5em}{[Q]} \parbox[t]{48.2em}{F.\ Quevedo,
   {\sl Supersymmetry and extra dimensions},
   graduate course, Department of Physics, University of Cambridge, spring 2006; accessible via YouTube.
  } 
   
} 

\end{titlepage}


\newpage
$ $

\vspace{-3em}

\centerline{\sc D-Brane Fantasy on Noncommutative Mirror Symmetry, Prelude}

\vspace{2em}


\begin{flushleft}
{\Large\bf 0. Introduction and outline}
\end{flushleft}

\begin{flushleft}
{\bf The D-Brane Fantasy\footnote{This
                                 meant-to-be-light  Introduction to why the D-brane Fantasy arises unavoidably involves vast amount of literatures
								  on the background that we are both unqualified and of no intention to summarize.								
							  Readers are referred to the sampled high-lighted words and the limited, incomplete sample list
						          in References, which influenced us greatly, for keyword search and further details.
                         } 
              }
\end{flushleft}
{\footnotesize
\begin{itemize}
 \item[$\tinybullet$] fantasy\:\; {\scriptsize noun}
   \vspace{-.6ex}
   \begin{itemize}
     \item[1:] the power or process of creating especially unrealistic or improbable mental images in\\ response to psychological need
     \item[2:] a creation of the imaginative faculty whether expressed or merely conceived: such as
      \vspace{-.6ex}
	  \begin{itemize}
	   \item[a:] a chimerical or fantastic notion
	   \item[b:] imaginative fiction featuring especially strange settings and grotesque characters
	   \item[c:] {\sc fantasia} sense 1 [a free usually instrumental composition not in strict form]
	   \item[d:] a fanciful design or invention	
	  \end{itemize}
	 \item[3:] $\cdots\cdots\cdots\cdots$
   \end{itemize}
   \hspace{-1em}[https:/\!/www.merriam-webster.com/dictionary/fantasy]
\end{itemize}	
} 

\bigskip

A {\it quantum field theory} has very rich contents,
  some of algebraic nature (e.g.\ {\it operator algebras/rings}, {\it state-operator correspondence}, {\it supersymmetry}),
  some of geometric nature
     (e.g.\ {\it moduli of vacua}, {\it nonperturbative sector},
	 {\it Wilson's theory-space} and the natural metric thereupon from the $2$-point function of fields involved,
	              {\it renormalization-group flow}, {\it supersymmetry} again),
  some of analytic  nature (e.g.\ {\it equations of motion}, {\it path-integrals},  {\it nonperturbative sector} again,
       {\it renormalization-group flow} again),
  some of combinatorial nature (e.g.\ {\it Feynman diagrams}),
  and some of categorical nature (e.g.\ {\it path-integral} again).
This makes it a fusion oven for different branches of mathematics.
A quantum field theory can have different presentations,
   which leads to the notion of {\it duality} both in physics and also on mathematical objects involved.
For example, seemingly different supersymmetric gauge theories can have isomorphic moduli of vacua
   but with their Coulomb branch and Higgs branch exchanged, (cf.\ {\it symplectic duality}).

Among all these dualities is {\it Mirror Symmetry}, a version of which states that
  \begin{itemize}
   \item[$\cdot$] {\bf [string world-sheet aspect of Mirror Symmetry]}\hspace{1em}
    {\it A supersymmetric string theory on different Calabi-Yau manifolds can be isomorphic
	        as a $2$-dimensional field theory, with different sectors of the operator rings switched.}
  \end{itemize}	
This leads to the notion of  {\it mirror pairs} $(M, W)$ of Calabi-Yau manifolds for which
    the complex (resp.\ complexified K\"{a}hler) deformation moduli of $M$ corresponds to
    the complexified K\"{a}hler (resp.\ complex) deformation moduli of $W$; and vice versa.
These deformations of the string target-space are encoded in different sectors of the operator ring of
       the $2$-dimensional supersymmetric string world-sheet theory
   and the theory with target $M$ and that with target $W$ differ only by an isomorphism that exchanges these sectors.
Many such mirror pairs were constructed, for example, in the {\sl Philip$\hspace{.7ex}$Candelas}'s group via toric geometry.
(See [Gr] for a review of (commutative) Mirror Symmetry and [MS] for sample work, both up to 1996,    and,
      e.g.\  [As2],   [A-B-C-D-G-K-M-S-S-W], [H-I-V], [H-K-K-P-T-V-V-Z]
	            for a glimpse of how the contents of Mirror Symmetry are further enriched by D-branes.)

However, the existence of rigid Calabi-Yau manifolds
   (i.e\ a K\"{a}hler manifold with vanishing first Chern class that has no complex deformations)
   renders such a mirror correspondence imperfect since the mirror to such a Calabi-Yau manifold won't exist in the above sense.
There are potentially a few cures to keep the mirror symmetry perfect, including
   (1) introduction of {\it non-K\"{a}hler Calabi-Yau manifolds},
   (2) taking other class of $2$-dimensional quantum field theories (e.g.\ {\it Landau-Ginzburg models})
                as the mirror candidate to the superstring theory on a rigid Calabi-Yau manifold.  and
   (3) introduction of {\it noncommutative Calabi-Yau spaces} and hope they provide the missing partners.
Method (3) is our focus here and there has been some mathematical progress made in this line,
    e.g.\  [Gi2], [Kel2], [Ko], [L-L-R], [Or]. [VdB2],
 though much remains to be done/understood for such {\it Noncommutative Mirror Symmetry} to reach the same level  of
   the current understanding of (Commutative) Mirror Symmetry.
  
Having ``reviewed" the world-sheet aspect of the stringy origin of  Mirror Symmetry
   -- another aspect comes from the lower-dimensional effective field theory from compactifications of
       Type IIA vs.\ Type IIB superstring theories , cf.\ Motif/Question/Project~2.2.7 --
a natural question arises:
  \begin{itemize}
   \item[{\bf Q.}]
 {\it What is the analogue statement to {\bf [string world-sheet aspect of Mirror Symmetry]} for Noncommutative Mirror Symmetry?}
  \end{itemize}
To reason the most probable answer to this question,
  assume that $Y^\nc$ is a noncommutative Calabi-Yau space in a setting
  that at least it makes sense to talk about local function-rings on $Y^\nc$, i.e.\ $Y^\nc$ is noncommutative ringed space
  supported on a topological space $Y$ with a noncommutative structure sheaf ${\cal O}_Y^\nc$.
To have an analogue to {\bf [string world-sheet aspect of Mirror Symmetry]},
 one needs {\it maps}
  $$
     \varphi\; :\;  ???\;  \longrightarrow\;  Y^\nc\,:=\, (Y, {\cal O}_Y^\nc)
  $$
  as dynamical fields on $???$ for which a suitable (most likely supersymmetric) quantum field theory on $???$ would be constructed.
Since local functions on $Y$ are pulled back under $\varphi$ to  local functions on $???$,
 one expects $???$ be a ringed space $(X, {\cal O}_X)$ as well.
Furthermore, since ${\cal O}_Y^\nc$  is noncommutative, one expects the structure sheaf on $X$ to be noncommutative as well,
 else local functions on $Y^\nc$ can be pulled back to $X$ only after passing through the commutative quotient sheaf
     ${\cal O}_Y^\nc/[{\cal O}_Y^\nc, {\cal O}_Y^\nc]$, which could be the zero-sheaf.
 Here,
    $[{\cal O}_Y^\nc, {\cal O}_Y^\nc]$ is the two-sided ideal sheaf of ${\cal O}_Y^\nc$
	   generated the commutator of elements in ${\cal O}_Y^\nc$.
Re-denote ${\cal O}_X$ by ${\cal O}_X^\nc$ to manifest the noncommutative feature on $X$     and
  let $X^\nc :=(X, {\cal O}_X^\nc)$;
 one now concludes that
 $$
     ???\; =\; \mbox{a {\it noncommutative} ringed space $X^\nc\,:=\, (X, {\cal O}_X^\nc)$}\,.
 $$
Thus,  to have an analogue to {\bf [string world-sheet aspect of Mirror Symmetry} for Noncommutative Mirror Symmetry,
 one has to seek a candidate of the noncommutative $X^\nc$.
Luckily, (Super) String Theory itself provides a candidate:
In contrast to the world-sheet of a fundamental string,
     the world-volume of stacked D-branes carries an Azumaya noncommutative structure
	     ([H-W], [Wi2] ; [L-Y1: Sec.\ 2] (D(1)),  [L-Y4] (D(11.1))),
	   allowing it to directly serve as a probe into noncommutative target-spaces.
Hence, one furher concludes that:
   $$
    \begin{array}{rcl}
        X^\nc     &   =    &   \mbox{\it world-volume of stacked dynamical D-branes} \\[1.2ex]
	                   &    =   &   X^{\!A\!z} \;
	                        :=\; (X,  {\cal O}_X^{A\!z}:=\Endsheaf_{{\cal O}_X^{\,\Bbb C}}({\cal E}))\,,
   \end{array}							
  $$
  where $X$ is the underlying $C^\infty$-manifold, possibly with boundary, for the world-volume of some stacked D-branes,
			  ${\cal E}$ is the Chan-Paton sheaf on $X$,   and
			  ${\cal O}_X^{A\!z}$ is the endomorphism sheaf of ${\cal E}$.
This train of reasoning leads to a D-brane fantasy:
	 \begin{itemize}
	  \item[]  {\bf [D-brane fantasy: D-brane world-volume aspect of Noncommutative Mirror Symmetry]}\hspace{1em}
      {\it Mirror Symmetry between noncommutative Calabi-Yau spaces (i.e.\ noncommutative Mirror Symmetry)
             corresponds to different realizations of a supersymmetric D-brane world-volume quantum field theory,
	         similar to/generalizing the string world-sheet aspect for Mirror Symmetry between (commutative) Calabi-Yau manifolds.}	
     \end{itemize}	
Once the world-volume theory candidates, i.e.\ dynamical stacked D-branes,  are chosen,  			
then, driven by this fantasy\footnote{This
                        D-brane fantasy should lie behind the non-abelian generalization of {\it gauged linear sigma model} of [Wi1]
                          studied by several string-theory groups.
                       The details of such a link deserves a separate work. 	
                                                         },  
 a long list of questions/themes follow, guided by the well-developed Commutative Mirror Symmetry, covered in three broad directions:
 {\it What are these noncommutative Calabi-Yau target space $Y^\nc$ (and structures thereupon)?}
 {\it What is the world-volume quantum field theory in this setting?}
 {\it Is there a justification of this D-brane fantasy?}

In the current notes,
  a class of noncommutative ringed spaces $Y^\nc$ shadowing over a $C^\infty$-scheme with corners
    are constructed from gluing local noncommutative crepant resolutions of Gorenstein isolated singularities.
Dynamical D-branes on such noncommutative target-spaces are defined as maps/morphisms
       from an Azumaya manifold with a fundamental module with a connection $\nabla$ thereto.
  The notion of $\nabla$-adjusted kinetic energy for such a map is given via
    the basic noncommutative differential calculus developed earlier in [L-Y4] (D(11.1)).
  This provides an action functional for dynamical D-branes on such noncommutative spaces
      in parallel to the Polyakov action functional for fundamental bosonic strings on a commutative target-space.
  This sets up a basic stage to begin with for the realization of the D-brane fantasy on Noncommutative Mirror Symmetry.
Questions beyond are sampled along the discussion.

\bigskip
\bigskip

\noindent
{\bf Convention.}
 References for standard notations, terminology, operations and facts are\\
  (1) algebraic geometry: [Hart], [Ei], [E-H];\;
              $C^{\infty}$ algebraic geometry: [Joy2], [FS-J], [Ka], [Bo-K];\;
              noncommutative algebraic geometry: [Gi1];\;
			  log algebraic geometry: [Og];\hspace{.6em}
  (2) noncommutative differential geometry: [GB-V-F], [B-M2];
              microlocal geometry:  [Na1], [Na2], [K-S]; \hspace{.6em}
  (3) string theory and D-branes: [G-S-W], [Pol4], [B-B-S]; [Pol3], [Ba], [Joh], [Sz].
 \begin{itemize}
  \item[$\cdot$]
   For clarity, the {\it real line} as a real $1$-dimensional manifold is denoted by ${\Bbb R}^1$,
    while the {\it field of real numbers} is denoted by ${\Bbb R}$.
   Similarly, the {\it complex line} as a complex $1$-dimensional manifold is denoted by ${\Bbb C}^1$,
    while the {\it field of complex numbers} is denoted by ${\Bbb C}$.
	
  \item[$\cdot$]	
  The inclusion `${\Bbb R}\subset{\Bbb C}$' is referred to the {\it field extension
   of ${\Bbb R}$ to ${\Bbb C}$} by adding $\sqrt{-1}$, unless otherwise noted.

   
  \item[$\cdot$]
   All manifolds are paracompact, Hausdorff, and admitting a (locally finite) partition of unity.
   We adopt the {\it index convention for tensors} from differential geometry.
    In particular, the tuple coordinate functions on an $n$-manifold is denoted by, for example,
    $(y^1,\,\cdots\,y^n)$.
   However, no up-low index summation convention is used.
   
  \item[$\cdot$]
   ${\Bbb C}^n$ as a ${\Bbb C}$-vector space vs as a complex manifold vs both.
   Similarly, for ${\Bbb R}^n$

  
  \item[$\cdot$]
  For this note, `{\it differentiable}', `{\it smooth}', and $C^{\infty}$ are taken as synonyms.
  
  %
  
  \item[$\cdot$]
   {\it section} $s$ of a sheaf or vector bundle vs.\ dummy labelling index $s$

  %
  %
  %
  
  \item[$\cdot$]
   $\Spec R $ ($:=\{\mbox{prime ideals of $R$}\}$)
         of a commutative Noetherian ring $R$  in algebraic geometry\\
   vs.\ $\Spec R$ of a $C^\infty$-ring $R$
  ($:=\Spec^{\Bbb R}R :=\{\mbox{$C^\infty$-ring homomorphisms $R\rightarrow {\Bbb R}$}\}$)

  \item[$\cdot$]
  {\it morphism} between schemes in algebraic geometry
    vs.\ {\it $C^{\infty}$-map} between $C^{\infty}$-manifolds or $C^{\infty}$-schemes
         	in differential topology and geometry or $C^{\infty}$-algebraic geometry
			
  \item[$\cdot$]
   group {\it action} vs.\  {\it action} functional for D-branes

  \item[$\cdot$]
   {\it metric tensor} $g$ vs.\ element $g^{\prime}$   in a {\it group} $G$ vs.\
      gauge coupling constant $g_{\gaugescriptsize}$

  %
  
  \item[$\cdot$]
   The `{\it support}' $\Supp({\cal F})$
    of a quasi-coherent sheaf ${\cal F}$ on a scheme $Y$ in algebraic geometry
     	or on a $C^\infty$-scheme in $C^\infty$-algebraic geometry
    means the {\it scheme-theoretic support} of ${\cal F}$
   unless otherwise noted;
    ${\cal I}_Z$ denotes the {\it ideal sheaf} of
    a (resp.\ $C^\infty$-)subscheme of $Z$ of a (resp.\ $C^\infty$-)scheme $Y$;
    $l({\cal F})$ denotes the {\it length} of a coherent sheaf ${\cal F}$ of dimension $0$.

  \item[$\cdot$]
   For a sheaf ${\cal F}$ on a topological space $X$,
   the notation `$s\in{\cal F}$' means a local section $s\in {\cal F}(U)$
      for some open set $U\subset X$.

  \item[$\cdot$]	
   For an ${\cal O}_X$-module ${\cal F}$,
    the {\it fiber} of ${\cal F}$ at $x\in X$	 is denoted ${\cal F}|_x$
	while the {\it stalk} of ${\cal F}$ at $x\in X$ is denoted ${\cal F}_x$.
  
  \item[$\cdot$] {\it Quantum topology} in the sense of quantization of General Relativity --
     in which even the {\it topology is fluctuating and needs to be summed over} --
     vs.\ that in the sense of {\it quantum invariants} in Low-Dimensional Topology

 \end{itemize}

 \bigskip
 \bigskip
   
\begin{flushleft}
{\bf Outline}
\end{flushleft}             		
\nopagebreak
{\small
\baselineskip 12pt  
\begin{itemize}
 \item[0]
  Introduction

 \item[1]
  Bundle settlements and $C^\infty$ Azumaya-type noncommutative ringed spaces
   \vspace{-.6ex}
   \begin{itemize}
     \item[1.1]
      $C^\infty$-vector-bundle settlements over a $C^\infty$-manifold with corners
      
	 \item[1.2]
	  Sheaves associated to vector-bundle settlements
	
	 \item[1.3]
	  Example: The local tensorial system from a seed system
	
	 \item[1.4]
	 $C^\infty$ Azumaya-type noncommutative ringed space from a mildly singular $C^\infty$-bundle settlement
	 on a mildly singular $C^\infty$-manifold with corners
   \end{itemize}
   
 \item[2]
  Seed systems over a Gorenstein singular Calabi-Yau space
   \vspace{-.6ex}
   \begin{itemize}
    \item[2.1]
	 Noncommutative crepant resolution of a Gorenstein isolated singularity and the associated apical algebra
	
	\item[2.2]
	 Seed systems over a Gorenstein singular Calabi-Yau space
  \end{itemize}
  
 \item[3]
  Dynamical D-branes on a $C^\infty$ Azumaya-type noncommutative ringed space
   \vspace{-.6ex}
   \begin{itemize}
     \item[3.1]
     ($C^\infty$) maps from a D-brane world-volume to a $C^\infty$ Azumayata-type noncommutative ringed space
 
     \item[3.2]
	 Basic ingredients for noncommutative Riemannian geometry on $Y^{nc}$
	
     \item[3.3]
	 An action functional for maps $\varphi$, given a connection $\nabla$ and metrics
   \end{itemize}
%
%
%
%
    %
\end{itemize}
} 

\newpage

\section{Bundle settlements and $C^\infty$ Azumaya-type noncommutative ringed spaces}

In Sec.\ 1.1, we introduce a new geometric object:
 {\it $C^\infty$ ${\Bbb C}$-vector-bundle settlements} over a $C^\infty$ manifold $Y$ with corners.
They serve as a precursor to a class of {\it $C^\infty$ noncommutative ringed space} that will be introduced in Sec.\ 1.2.
After giving a class of examples of such objects in Sec.$\hspace{.7ex}$1.3,
  we generalize the setting in Sec.\ 1.4 to allow mild isolated singularities on $Y$.
    %

\bigskip

\subsection{$C^\infty$-vector-bundle settlements over a $C^\infty$-manifold with corners}					

A few basic definitions are given,  which lead to the notion of $C^\infty$-vector-bundle settlement.
Direct sum and tensor product of such objects are also discussed.

\bigskip

\begin{flushleft}
{\bf $C^\infty$-vector-bundle settlements over a $C^\infty$-manifold with corners}
\end{flushleft}
\begin{definition} {\bf [$C^\infty$-scheme with corners]}\; {\rm
  Readers are referred
       to [Joy2: Sec.\ 2]     for the definition of {\it $C^\infty$-rings} and their homomorphisms;
       to  [Joy2: Sec.\ 4]    for the definition of  {\it $C^\infty$-schemes} and their morphisms;
	   to [FS-J: Sec.\ 4]   for the definition of  {\it $C^\infty$-rings with corners} and their homonorphisms;   and
       to [FS-J: Sec.\ 5] and [Ka: Sec.\ 4.8] for the definition of  {\it $C^\infty$-schemes with corners} and their morphisms.
 Some {\it logarithmic algebraic geometry} taste for the treatment of {\it corners} are referred to [G-M] and [Og].
 Assuming some background of Algebraic Geometry as in [E-H: Chapter I] or [Hart: Chapter II],
   then the following is enough for the current notes:
   \begin{itemize}
     \item[$\cdot$]
	  A $C^\infty$-scheme with corners is a ringed space $(Y, {\cal O}_Y)$
	     with the structure sheaf ${\cal O}_Y$ a sheaf of $C^\infty$-rings on the underlying topological space $Y$.
	  ${\cal O}_Y$ may contain {\it nilpotent elements}.

	 \item[$\cdot$]
	  The notion of an {\it ideal sheaf} ${\cal I}$  of ${\cal O}_Y$ is defined in the same way as in Algebraic Geometry.
	  The quotient ${\cal O}_Y$-algebra ${\cal O}_Y/{\cal I}$ defines a $C^\infty$-subscheme
	    $(Z, {\cal O}_Z:={\cal O}_Y/{\cal I})$ with corners in $Y$.
	
	 \item[$\cdot$]
	  The notion of {\it scheme-theoretic image} of a morphism between $C^\infty$-schemes with corners is defined in the same way
	    as in Algebraic Geometry.
     
	\item[$\cdot$]
	 The notion of {\it scheme-theoretic support} of an ${\cal O}_Y$-module ${\cal M}$ on $Y$
	   is defined in the same way as in Algebraic Geometry.
   \end{itemize}
}\end{definition}

\medskip

\begin{example} {\bf [$C^\infty$-manifold with corners]}\; {\rm
(Adapted from [Joy1: Sec.\ 2] and [Ka: Sec.\ 4.1] as guided by [Hi: Chapter I] to fix the terminology;
    readers are referred to ibidem for complete details.)
 A {\it pre-$n$-manifold with corners}, $n\ge 1$,  is a paracompact Hausdorff topological space $Y$ such that
  for all $p\in Y$, exactly one of the following two is satisfied:
   \begin{itemize}
    \item[(1)]
	 There are an open neighborhood $U\ni p$ and a continuous map $\rho: U \rightarrow {\Bbb R}^n$
	    such that $U$ and $\rho(U)$ are homeomorphic. In this case, we say that $p$ is a {\it manifold point}.
	  The set of manifold points of $Y$ is called the {\it interior} of $Y$, denoted by $\Int Y$.
   
   \item[(2)]
     There are an open neighborhood $U\ni p$ and
	    a continuous map $\rho: U\rightarrow {\Bbb R}^n_k:= [0, \infty)^k\times {\Bbb R}^{n-k}$, $1\le k\le n$,
		such that $\rho(p)\in \{\mathbf{0}_k\}\times {\Bbb R}^{n-k}$,
		                        where  $\mathbf{0}_k$ is the origin of $[0,\infty)^k$,   and
		         that $U$ and $\rho(U)$ are homeomorphic.
	  In this case, if $k=1$ (resp.\ $k\ge 2$) we say that $p$ is a {\it boundary point} (resp.\ {\it corner point}).
	  The set of all such points are denoted by $\partial Y$, named the {\it boundary} of $Y$.
	  We use also the phrase that {$p$ lies in the codimension-$k$ stratum of $\partial Y$}.	
   \end{itemize}
  The pair $(U, \rho)$ is called a {\it (local) chart} of $Y$.
  Two charts $(U_i, \rho_i)$ and $(U_j,\rho_j)$   are said to have {\it $C^\infty$/smooth overlap}
    if the {\it coordinate change}
	  $$
	    \rho_{ij}:= \rho_j \circ \rho_i^{-1}\;:\;
		   \rho_i(U_i\cap U_j)\;\longrightarrow\;  \rho_j(U_i\cap U_j)
	  $$
   is a $C^\infty$-diffeomorphism.
 Here:
    \begin{itemize}
	 \item[] {\bf Definition 1.2.1.1.  [$C^\infty$/smooth map between sets]}\;\;
	  Let $C_1\subset {\Bbb R}^{m_1}$ and $C_2\subset {\Bbb R}^{m_2}$ be arbitrary subsets.
      We say that a {\it map} $f:A\rightarrow B$	 is {\it $C^\infty$} (or synonymously,  {\it smooth})
	    if there exist an open $V_1\supset C_1$ in ${\Bbb R}^{m_1}$ and
		                      an open set $V_2\supset C_2$ in ${\Bbb R}^{m_2}$
        such that $f$ is the restriction of a $C^\infty$/smooth map $\widehat{f}: V_1\rightarrow V_2$.
   \end{itemize}	
 A collection ${\frak A}$ of charts of $Y$ that form a covering of $Y$ and
            such that every pair of charts in ${\frak A}$ has $C^\infty$ overlap    is called a {\it $C^\infty$ atlas} on $Y$.
 Any such ${\frak A}$  has a unique enlargement to a maximal $C^\infty$ atlas ${\frak A}^\prime$ on $Y$.
 We say that ${\frak A}^\prime$ is a {\it $C^\infty$ differential structure} on $Y$.
 A pre-$n$-manifold $Y$ with corners together with a choice of a maximal $C^\infty$ atlas
    is called a {\it $C^\infty$ $n$-manifold with corners}.
   
 A function $f:Y\rightarrow {\Bbb R}$ on a $C^\infty$-manifold with corners is said to be {\it $C^\infty$/smooth}
     if it is smooth with respect to every local chart(s).
  This determines the structure sheaf ${\cal O}_Y$ on $Y$ as a $C^\infty$-scheme with corners.
  
 Similarly, a map $f:Y\rightarrow Y^\prime$ between $C^\infty$-manifolds with corners is said to be {\it $C^\infty$/smooth}
    if $f$ is smooth respect to the atlases; namely,  	
	  for every $y^\prime\in Y^\prime$ and $y\in f^{-1}(y)$,
    	  there exist charts $(V^\prime, \rho^\prime)$ around  $y^\prime$ and
	                                     $(V, \rho)$ around $y$ such that $f(V)\subset V^\prime$ and that
		$\rho^\prime\circ f	\circ  \rho^{-1}:   \rho(V)\rightarrow \rho^\prime(V^\prime)$
    	is $C^\infty$/smooth.
}\end{example}

\bigskip

Throughout this subsection, we'll let
 $Y$ be a $C^\infty$ $n$-manifold with corners, treated as a $C^\infty$-scheme with corners;
 ${\cal O}_Y$ be its structure sheaf as such;  and
 ${\cal O}_Y^{\,\Bbb C}:= {\cal O}_Y\otimes_{\Bbb R}{\Bbb C} $  be the complexification of ${\cal O}_Y$.
  
\bigskip

\begin{definition} {\bf [regularly partitioned $C^\infty$-manifold with corners]}\;  {\rm				
 (1)
    A {\it partition} ${\cal P}$ of $Y$  is a decomposition $Y=\bigcup_{i\in I}Y_i$ of $Y$
	  into a collection of connected $C^\infty$ $n$-manifolds $Y_i$ with corners embedded in $Y$.
	Each component $Y_i$ in the decomposition is called a {\it block}.
	It is required that the interior of the blocks $Y_i$, $i\in I$, are disjoint from each other.
	We'll denote $Y$ equipped with a partition ${\cal P}$ by $_{\cal P}Y$ when need to emphasize this additional structure on $Y$.
	Note that if $Y^\prime$ is another $C^\infty$ manifold with corners with a partition ${\cal P}^\prime$,
	 then the product $Y\times Y^\prime$ is naturally equipped with a {\it product partition} ${\cal P}\times {\cal P}^\prime$
		via taking the product $Y_i\times Y^\prime_j$ of blocks $Y_i$ in ${\cal P}$ and blocks $Y^\prime_j$ in ${\cal P}^\prime$.
		
 (2)
    A partition ${\cal P}$ of $Y$ is called {\it regular} if the following two conditions are satisfied:
	    \begin{itemize}
		 \item[(a)]
		  Every interior point $p\in Y$ has a neighborhood $U$ such that, with the induced partition on $U$,
	       $_{\cal P}U$ is diffeomorphic to
		       a neighborhood of the origin of $_{\cal P}{\Bbb R}^k\times {\Bbb R}^{n-k}$ for some $k\in\{0,1, \cdots, n\}$,
	        where $_{\cal P}{\Bbb R}^k$ is ${\Bbb R}^k$ equipped with the {\it standard partition}
		    by the $2^k$-many coordinate quadrants  (with $_{\cal P}{\Bbb R}^0$ set to be a point by convention).
			
		 \item[(b)]
		  Every boundary or corner point $p\in\partial Y$ with a neighborhood $U\simeq$ a neighborhood of the origin of
	       $[0,\infty)^l\times {\Bbb R}^{n-l}$ has a neighborhood $U^\prime\subset U$, with the induced partition,
		   diffeomorphic to a neighborhood of the origin in
		       $[0,\infty)^l\times_ {\,\cal P}${${\Bbb R}^k \times {\Bbb R}^{n-l-k}$}, for some $k\in\{0,1, \cdots, n-l\}$.
	
	  \end{itemize}
	 The set of points of $Y$ that share the same $k$ in a regular partition ${\cal P}$
	    is called the {\it codimension-$k$ stratum} of ${\cal P}$, denoted $_{\cal P}^{\,k}Y$.		
	  By construction, $Y$ is the disjoint union $\amalg_{k=0}^n{_{\,\cal P\!}^{\:\:k}Y}$ of smooth submanifolds.
	  They satisfy Whitney's Condition (B) and  give a Whitney stratification of $Y$ that involves only codimension-$1$ submanifolds
	    with simple normal crossing singularities.		
   
 (3)
    Let ${\cal P}$ and ${\cal P}^\prime$ be partitions of $Y$.
	We say that ${\cal P}^\prime$ is a {\it refinement} of ${\cal P}$
	  if ${\cal P}^\prime$ is obtained by a partition, in the sense of (1), of each block of ${\cal P}$.
		
 (4)
     A partition ${\cal P}$ of $Y$ is called {\it weakly regular} if each $p\in Y$
	   has a neighborhood $U$ to which the restriction ${\cal P}|_U$ can be refined, in the sense of (3), to a regular partition of $U$.
   
{\sc Figure}~1-1-1.
      %
}\end{definition}
		
\begin{figure}[htbp]
  \bigskip
   \centering
   \includegraphics[width=0.80\textwidth]{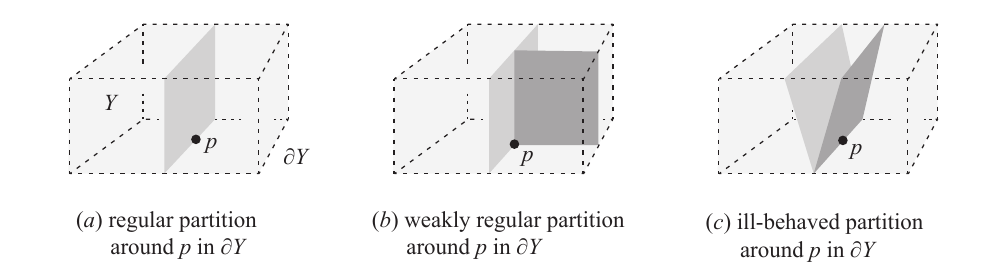}
 
  \bigskip
 \centerline{\parbox{13cm}{\small\baselineskip 12pt
  {\sc Figure}~1-1-1.
   Three partitions of a neighborhood of a $p\in \partial Y$ are illustrated:
     ({\it a}) is a regular partition;
	 ({\it b}) is a weakly regular partition which can be refined to a regular partition;
	({\it c}) is neither and there is no way it can be refined to a regular partition.
  }}
\end{figure}
      
\bigskip	

Loosely speaking, a `$C^\infty$-bundle settlement' is
  a collection of $C^\infty$-bundles $F_i$, $i\in I$,  over each block of a regularly partitioned $C^\infty$-manifold with corners   plus
  a gluing system of $C^\infty$ bundle embeddings among the restriction of $F_i$'s to the facets  of $\partial Y_i$, $i\in I$.
For the application in the current work, we take these bundles to be complex vector bundles. Then, in detail:
(Cf.\ Figure~1-1-2.)
		
\begin{figure}[htbp]
  \bigskip
   \centering
   \includegraphics[width=0.80\textwidth]{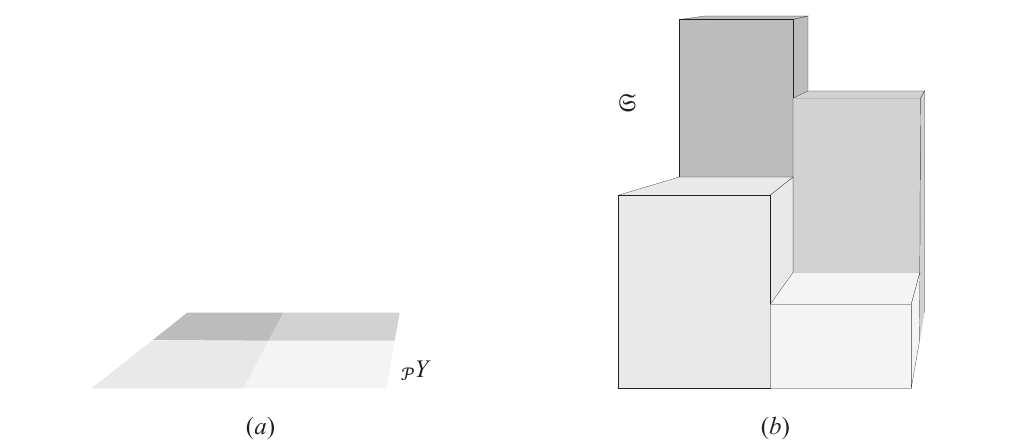}
 
  \bigskip
 \centerline{\parbox{13cm}{\small\baselineskip 12pt
  {\sc Figure}~1-1-2.
     ({\it a}) A regularly partitioned manifold $_{\cal P}Y$ with corners.
     ({\it b}) A $C^\infty$ ${\Bbb C}$-vector-bundle settlement ${\frak S}$	 over $_{\cal P}Y$ is indicated.
	                          Here the rank of a $C^\infty$ ${\Bbb C}$-vector-bundle block $F_i$ over $Y_i\in {\cal P}$
							     is indicated by the hight over the block.
							  A lower rank bundle block is glued to a higher rank bundle block in a consistent way
							    along the strata ${_{\cal P}}^{\:\ge 1}Y$ of $_{\cal P}Y$ where they meet.
  }}
\end{figure}
				
\bigskip

\begin{definition}{\bf [$C^\infty$-complex-vector-bundle settlement]}\; {\rm
A {\it $C^\infty$-complex-vector-bundle settlement}\footnote{\noindent{\it Cf.}$\hspace{.7ex}$From
                        Wikipedia on `{\it Human settlement}'
						   ({\tt https:/\!/en.wikipedia.org/wiki/Human\_settlement }):
						``In geography, statistics and archaeology, a {\it settlement}, locality or populated place is a community of people living
						     in a particular place.
						  The complexity of a settlement can range from a minuscule number of dwellings grouped together to the largest of cities
						    with surrounding urbanized areas.
						  Settlements include homesteads, hamlets, villages, towns and cities. $\cdots\cdots$."
                       In our adoption of the term, `{\it human}' is replaced by `{\it $C^\infty$-bundle}'.
                                                                                             }
	 \,${\frak S}$ over a $C^\infty$ regularly partitioned $n$-manifold $_{\cal P}Y$ with corners
 consists of the following data:
 \begin{itemize}	
  \item[(1)]	   [{\it $C^\infty$-bundle blocks}] \hspace{1em}
   A collection of $C^\infty$ ${\Bbb C}$-vector bundles $F_i$ over each block $Y_i$, $i\in I$,  in ${\cal P}$.	
   Here, by a {\it $C^\infty$ ${\Bbb C}$-vector bundle $F_i$ over each block $Y_i$},
     we mean that there is a $C^\infty$ ${\Bbb C}$-vector bundle $F_i^\prime$
	   on an open set $U_i^\prime\subset Y$  that contains $Y_i$ such that $F_i=F_i^\prime|_{Y_i}$.
   To serve as the reference for the notion of smoothness, we'll fix an $F_i^\prime$ for each $i$ and name it
     the {\it extended bundle block} for $F_i$.

 \item[(2)]   [{\it $C^\infty$-bundle gluing system along boundary facets}] \hspace{1em}
  A gluing system of $C^\infty$ ${\Bbb C}$-vector-bundle  embeddings
	   among the restriction of $F_i$'s to the facets of $\partial Y_i$'s, $i\in I$, defined as follows:
         %
    \vspace{-.6ex}
    \begin{itemize} 	
      \item[(2.1)]   [{\it $C^\infty$-gluing data}] \hspace{1em}
	   The gluing data is assigned inductively over the codimension-$k$ facets of the blocks in ${\cal P}$, for $k=1,\cdots, n$.
	
	  \item[]\hspace{1.2em}		
	   For $k=1$, let
	     $Z$ be a closed codimension-$1$ face in ${\cal P}$
		                               that lies in the intersection of two distinct blocks $Y_i$ and $Y_j$ in ${\cal P}$
	      with $\rank_{\Bbb C}F_i< \rank_{\Bbb C}F_j$.
	  Then assign a $C^\infty$ ${\Bbb C}$-vector-bundle embedding
	     $h_{Z, ij}: F_i|_Z \rightarrow F_j|_Z$ that is {\it smooth around $Z$}
	     in the sense that
			 \begin{itemize}
			  \item[$\cdot$]
			   there exists an open neighborhood $U\subset U^\prime_i\cap U^\prime_j$ of $Z$
			     such that there is a $C^\infty$ ${\Bbb C}$-vector-bundle embedding
			      $h_{U, ij}: F_i^\prime|_U \rightarrow F_j^\prime|_U$ with $h_{U, ij}|_Z=h_{Z, ij}$.
	        \end{itemize}	
		For the case as above but with $\rank_{\Bbb C}F_i=\rank_{\Bbb C}F_j$,
		 one fixes one direction, say $\overrightarrow{ij}$, and
		   require the $h_{Z, ij}$ above be a $C^\infty$ ${\Bbb C}$-vector-bundle isomorphism that is smooth around $Z$
		   (i.e.\ $h_{U, ij}$ is a smooth $C^\infty$ ${\Bbb C}$-vector-bundle isomorphism).
	  For convenience, an $h_{Z, ij}$ thus assigned will be called a {\it clutching map} along $Z$.
	  It may be encoded by the quiver $i\, \tinybullet\!\longrightlongarrow\!\tinybullet\, j$\,.
	
	  \medskip
	  \item[]
	  {\bf Definition 1.1.4.1. [clutching quiver]}\; {\rm
	       The quiver diagram that encode the $C^\infty$-gluing data over $Z$ is called
		     the {\it clutching quiver} of the bundle settlement along $Z$.
	     } 
	
	  \medskip
	  \item[]\hspace{1.2em}
      For $k=2, \cdots, n$,
         the $C^\infty$-gluing data are defined inductively
		   after introducing the notion of {\it normal cubes} of a point $p\in_{\,\cal P}${$Y$},
		   with the aid  of the clutching quivers that generalize the $k=1$ case for  bookkeeping.
		   		
	   \item[]\hspace{1.2em}
	    Let $p\in_{\,\cal P}${$Y$}, then $p\in _{\,\cal P\!}^{\:\:k}${$Y$} for a unique $k\in \{0,1,\cdots, n\}$
			    and there is a neighborhood of $p$ with the induced partition that is diffeomorphic to a neighborhood $U$ of the origin of
				$_{\cal P}{\Bbb R}^k\times {\Bbb R}^{n-k}$, if $p$ is an interior point of $Y$,   or of
				$[0,\infty)^l\times_{\,\cal P}${${\Bbb R}^k\times {\Bbb R}^{n-k-l}$}, for a unique $l$, if $p\in \partial Y$.
			 Under such a local diffeomorphism, define
			  a {\it cube} at $p\in_{\,\cal P}^{\:\:k}${$Y$} to be a submanifold with corners that is described by
			    $I^n_{p, \varepsilon}  := [-\varepsilon, \varepsilon]^n\subset U$, for $p$ interior,
				    $[0,\varepsilon]^l \times [-\varepsilon, \varepsilon]^{n-l}\subset U$, for $p\in\partial Y$, for some $\varepsilon>0$;
                  and
			  a {\it normal cube} at $p\in_{\,\cal P}^{\:\:k}${$Y\subset Y$} to be the submanifold with corners described by
			    $I^k_{p, \varepsilon}  := [-\varepsilon, \varepsilon]^k \times \mathbf{0}_{n-k}\subset U$, for $p$ interior,
				    $\mathbf{0}_l\times [-\varepsilon, \varepsilon]^k\times \mathbf{0}_{n-k-l}\subset U$, for $p\in\partial Y$,
					 for some $\varepsilon>0$.
			Here $\mathbf{0}_{n-k}$, $\mathbf{0}_l$, $\mathbf{0}_{n-k-l}$ are the origin of
			    ${\Bbb R}^{n-k}$, $[0,\infty)^l$, and ${\Bbb R}^{n-k-l}$ respectively.
			By construction,
  			  $_{\cal P}^{\:\:k}Y\cap I_{p,\varepsilon}^n= \mathbf{0}_k\times [-\varepsilon, \varepsilon]^{n-k}$
			       for $p$ an interior point of $Y$,
				 $[0,\varepsilon]^l\times\mathbf{0}_k\times [-\varepsilon, \varepsilon]^{n-k-l}$ for $p\in \partial Y$.
			 This justifies the name.
			
             \medskip			
			 \item[]
				{\bf  Lemma 1.1.4.2. [normal-cube bundle of facet is trivial]}\; {\it
			  Let $Z$ be a facet of the regular partition ${\cal P}$ of $Y$. Then $Z$ has a trivial normal-cube bundle.
			    } 
			
			 \medskip
			 \item[]  {\it Proof.}
			 Let $Z$ be a codimension-$k$ facet of ${\cal P}$.
			 Since all the blocks in ${\cal P}$ are embedded closed submanifolds with corners in $Y$,
			   every codimension-$1$ face of ${\cal P}$ has a trivial normal-interval  bundle.
			 Since ${\cal P}$ is regular, there are exactly $2^k$-many blocks that contain $Z$ as a codimension-$k$ facet in their boundary.
			 Thus, any of these blocks has exactly $k$-many boundary faces that contain $Z$.
			 The restriction of the normal-interval bundles of these faces to $Z$ constitute a globally trivialized $k$-frame bundle over $Z$.
             The lemma then follows.
	
			 \noindent\hspace{36em}$\square$
	
	  \item[]\hspace{1.2em}
	   Assume that the assignment of clutching maps is completed up to over codimension-$(k-1)$ facets of ${\cal P}$.
	   Let $Z$ be a codimension-$k$ facet of ${\cal P}$, $k\ge 2$.
	   Then
	     there are exactly $k$-many codimension-$1$ faces,  denoted $Z_1,\, \cdots,\, Z_k$, of ${\cal P}$
		      that contain $Z$ and intersect transversely along $Z$    and 	
	     the set of $2^k$-many vertices of any normal $k$-cube $I^k_{p, \varepsilon}$, $p\in Z\cap_{\,\cal P}^{\:\:k}${$Y$},
		  labels exactly the blocks of ${\cal P}$ that contain $Z$ in their boundary.
	   By induction,
	      the clutching maps over the $k$-many codimension-$(k-1)$ facets
		      $Z_1\cap \cdots Z_{l-1}\cap Z_{l+1}\cap \cdots Z_k$, $l=1,\,\cdots,\, k$, of ${\cal P}$
          are already defined. (Here, by convention, $Z_0:=Y$.)
      There restriction to over $Z$ gives part of the clutching maps over $Z$.
	  Let $Q_Z^\prime$ be the resulting preliminary clutching quiver along $Z$.
	
	   \medskip
	   \item[]
	   {\bf Definition 1.1.4.3. Definition [saturated/terminal clutching quiver]}\; {\rm
	      A quiver $Q_Z^\prime$ thus obtained is called {\it saturated} or synonymously {\it terminal}
		   if if there exists a distinguished vertex $v_\ast$ such that,
		     for any vertex $v$ of $Q_Z^\prime$,
			    there exists an edge-path in $Q_Z^\prime$ with initial vertex $v$ and terminal vertex $v_\ast$.
		 (Note that such a $v_\ast$ is not unique if $Q_Z^\prime$ contains an edge loop passing through $v_\ast$.)	
	    } 
		
      \medskip
      \item[]\hspace{1.2em}
      Continuing the discussion, 	
	  if $Q_Z^\prime$ is already saturated,
	  then  it means that there is a bundle block, say some $F_i$, of these $2^k$-many bundle blocks $F_{i^\prime}$ around $Z$,
	    such that, for all $i^\prime$, $F_i|_Z$ contains a copy of $F_{i^\prime}|_Z$ via a composition of existing clutching maps.
      In this case, we set the $C^\infty$-gluing data along $Z$
	        to be the set of the restriction of these existing $C^\infty$-gluing maps from induction.	
	  The quiver $Q_Z$ along $Z$ is thus $Q_Z^\prime$.
	
	  \medskip
      \item[]\hspace{1.2em}
	  Continuing the notation, if, on the other hand, $Q_Z^\prime$ is not yet saturated,
	  then one has to choose pairs of bundle blocks, say $F_{i_{j}}$ and $F_{i^\prime_j}$, from these  $2^k$-many bundle blocks
	   and add clutching map $h_{Z, \tinybullet\tinybullet}$ for the restriction of each pair to over $Z$
	       in such a way that $h_{Z, \tinybullet\tinybullet}$ is smooth around $Z$, as in the $k=1$ case,
	   until the corresponding clutching quiver $Q_Z$ is saturated. 	
      In this way, one obtains inductively a collection of clutching maps $\{h_{Z, ij}\}_{Z, ij}$ with $h_{Z,ij}$ smooth around $Z$,
	    where $Z$ runs over all the closed facets of ${\cal P}$.
		
     \item[(2.2)]	
	 [{\it consistency of gluings over $_{\cal P\!}^{\:k}Y$:
	                               cocycle conditions encoded in the clutching quivers $Q_Z$ along facets $Z$ of ${\cal P}$}]\hspace{1em}
       Note that after exhausting the clutching maps over all (codimension-$1$) faces/walls in ${\cal P}$,
	     the collection $\{\mbox{$F_i$ on $Y_i$}\}_{Y_i\in {\cal P}}$ automatically defines
	     a $C^\infty$ ${\Bbb C}$-vector-bundle settlement over $Y-\amalg_{k=2}^n{_{\,\cal P\!}^{\:\:k}Y}$.		 			
	   But to extend this bundle settlement further to over $\amalg_{k=2}^n{_{\,\cal P\!}^{\:\:k}Y}$,
	     a system of consistency conditions need to be imposed inductively.
   Suppose the consistency of gluings are imposed over all facets of ${\cal P}$ of codimension $\le k-1$.
   Let $Z$ be a codimension-$k$ facet of ${\cal P}$.
   Note that
     a composition of clutching maps over $Z$ corresponds to a unique edge-path in the quiver $Q_Z$; and vice versa.
   Then, the new consistency conditions on $\{h_{Z, ij}\}_{ij}$
       (i.e.$\hspace{.7ex}$those not from the restriction of clutching maps over codimension-$\le(k-1)$ facets containing $Z$)
	       are encoded by the set of equalities of distinct simple edge paths of the form
		   $$
		     \mbox{\it Path}_{\,1}\; =\; \mbox{\it Path}_{\,2}\,,
		   $$
		  where
			 \begin{itemize}
			  \item[$\cdot$]
			   {\it Path}$_{\,1}$ and {\it Path}$_{\,2}$ have the same initial vertex and the same final vertex,
			  \item[$\cdot$]
		      {\it Path}$_{\,1}$ contains a directed edge in $Q_Z-Q_Z^\prime$, and
			  \item[$\cdot$]
			{\it Path}$_{\,2}$ corresponds to the identity map when {\it Path}$_{\,1}$ is an edge-loop.
			\end{itemize}
    \end{itemize}  	
 \end{itemize}
{\sc Figure} 1-1-3.
      
 For convenience, we'll write `$F_i\in {\frak S}$' to mean that {\sl $F_i$ is a bundle block of ${\frak S}$}.
 Also, the following notations and sentences mean the same:
   {\sl ${\frak S}$ is a bundle settlement over $_{\cal P}Y$},
   {\sl $_{\cal P}{\frak S}$ is a bundle settlement over $Y$},
   {\sl $_{\cal P}{\frak S}$ is a bundle settlement over $_{\cal P}Y$}.
 When the regular partition ${\cal P}$ of $Y$ is meant to be kept implicit, we'll simply say that
       {\sl ${\frak S}$ is a bundle settlement over $Y$}.
}\end{definition}
				
\begin{figure}[htbp]
 \bigskip
  \centering
  \includegraphics[width=0.80\textwidth]{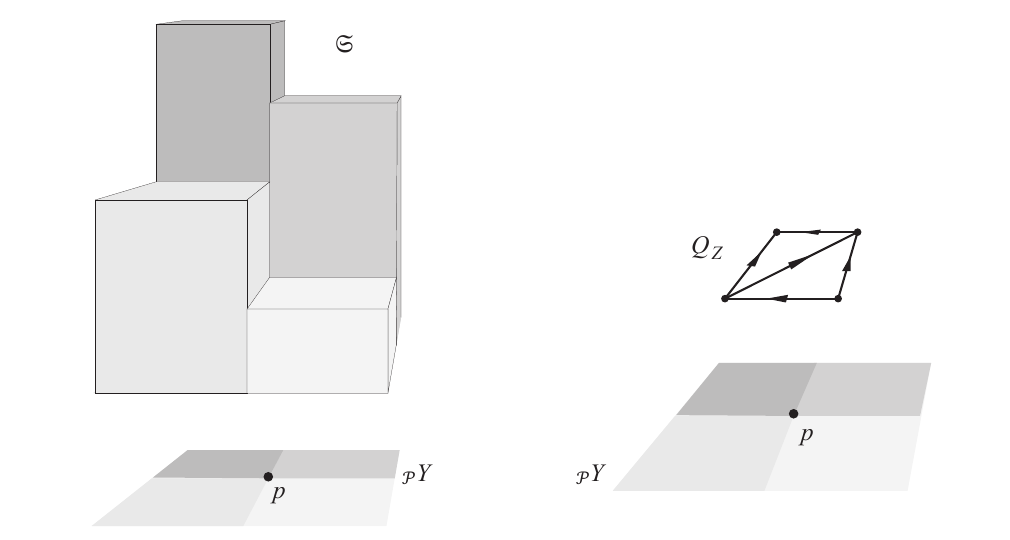}
 
  \bigskip
 \centerline{\parbox{13cm}{\small\baselineskip 12pt
  {\sc Figure}~1-1-3.
   The clutching quiver $Q_Z$ of a bundle settlement ${\frak S}$ along  a stratum $Z$
       of $_{\cal P}Y$ that contains $p$ in its interior is indicated.
    Here $Z$ is supressed and only a transverse slice to $Z\subset Y$ at $p$ is drawn.
  }}
\end{figure}

\medskip

\begin{definition}  {\bf [distinguished block, fiber and decorated fiber of vector-bundle settlement at point]}\;  {\rm
  Continuing Definition~1.1.4.
  Let $p\in  {_{\cal P}Y}$ be  a point in the interior of a codimension-$k$ facet $Z$ of $_{\cal P}Y$.
  By construction, the gluing system of $2^k$-many vector spaces $F_{i_1}|_p, \cdots, F_ {i_{2^k}}|_p$,
       from the blocks $F_{i_j}\in {\frak S}$ around $Z$,
    is encoded by a quiver with a distinguished vertex $v_\ast$.
  The bundle block $F_p:= F_{i_\ast}\in \{F_{i_1}, \cdots, F_{i_{2^k}}\}$ corresponding to $v_\ast$
    is called the {\it distinguished bundle block} of ${\frak S}$ at $p$.
  This is a bundle block of maximal rank along $\Int Z\ni p$.
  In terms of this,
    the {\it fiber} (resp.\ {\it decorated fiber}) of ${\frak S}$ at $p$ is defined to be
       the fiber $F_p|_p$ of the vector bundle $F_p$ at $p$
	     (resp.\ together with the partially ordered system of subspaces from the image of  $F_{i_j}|_p$, $i_j\ne i_\ast$, in $F_p|_p$
		               under a composition of the $C^\infty$ clutching maps in ${\frak S}$ along $Z$).
 In particular, the decorated fiber of a vector-bundle settlement is
  a vector space at an interior point of a block of $_{\cal P}Y$,
  a simple partial flag $F_i|_\tinybullet \subset F_j|_\tinybullet$ at an interior point of the walls, and
  a  vector space with a partially ordered system of subspaces (cf.\ multi-partial-flags) at a point in a higher-codimensional facet.
}\end{definition}

\medskip

\begin{remark} $[{\Bbb C}$-vector-bundle settlement in general$\,]$\; {\rm
 The above notion of bundle settlements can be generalized to over a weakly regular partition of $Y$.
 The consistency conditions for the gluing data remain encoded in stratum-wise locally constant quivers.
}\end{remark}

\medskip

\begin{motif-question-project}{\bf  [fortified/gated bundle settlement]}\; {\rm
 The above setting (more precisely choice) of the notion of a {\it bundle settlement} is enough for the purpose of the current work
   as such objects will give rise to a class of $C^\infty$ soft noncommutative schemes that serve as the target space for dynamical D-branes.
 However,
    from the point of view of a {\it log-like structure associated to a regular partition ${\cal P}$} on $Y$ and the input from the physics side,
	  i.e.\ {\it quantum field theory with boundary and corners} - in which case bulk, boundary, corners can each have fields living upon -,
  the version adopted here is only the most primitive form of the notion of a `bundle settlement'.
 One can add independently a smooth bundle over each codimension-$\ge 1$ irreducible stratum from ${\cal P}$ and
		glue the bundle-block $Y_i$, $i\in I$, through these additional bundles over boundary and corners.
 These additional bundles also have to glue among themselves similarly.
 Such a bundle settlement can be thought as {\it fortified} or {\it gated}.	
 One has to develop then the notion of {\it consistency of gluing} here and
  address even, for example, the notion of {\it $C^\infty$ global sections} of a bundle settlement in this context
  so that a good notion of associated sheaves can still be introduced.
}\end{motif-question-project}

\medskip

\begin{motif-question-project}{\bf  [cellular gauge theory \& n-category]}\; {\rm
 Pressing further along the line of Motif/Question/Project 1.1.7,
  the setting of bundle settlements, fortified or not, suggests a notion of {\it cellular gauge theory} and
  a certain {\it $n$-category} associated to it. Details?
 (Cf.\ Motif/Question/Project 1.2.8.)
 Once one goes this far, {\it can one generalize the cellular gauge theory to over a stratified space $Y$ of reasonably good singularities,
   e.g.\ the product of a manifold with corners with a $1$-dimensional web, or with a $2$-dimensional foam, or with ...?}	
}\end{motif-question-project}

\bigskip

\begin{flushleft}
{\bf Direct sum and tensor product  of concordant vector-bundle settlements}
\end{flushleft}
Let ${\frak S}_1 = (\{F_{1; i}\}_{i\in I}, \{h_{1; Z,  i_1 j_1}\}_{Z, i_1 j_1})$ and
      ${\frak S}_2 = (\{F_{2; i}\}_{i\in I}, \{h_{1; Z,  i_2 j_2}\}_{Z, i_2 j_2})$
  be two $C^\infty$ ${\Bbb C}$-vector-bundle settlements on $_{\cal P}Y$.

\bigskip

\begin{definition} {\bf [concordant bundle-settlements]}\;  {\rm
 We say that ${\frak S}_1$ and ${\frak S}_2$ are {\it concordant}
  if for all facets $Z$ of $_{\cal P}Y$, their corresponding clutching quivers $Q_{1; Z}$ and $Q_{2; Z}$ are identical.
}\end{definition}
 
\bigskip

\noindent
Assume that ${\frak S}_1$ and ${\frak S}_2$ are concordant.
Then not only the sets of the bundle blocks but also the $C^\infty$-gluing systems of the two share an identical  index system:
  $$
    {\frak S}_1 = (\{F_{1; i}\}_{i\in I}, \{h_{1; Z,  i j}\}_{Z, ij})
  	       \hspace{2em}\mbox{and}\hspace{2em}
    {\frak S}_2 = (\{F_{2; i}\}_{i\in I}, \{h_{1; Z,  i j}\}_{Z, ij}).
  $$
And the following data
   \begin{eqnarray*}
   {\frak S}_1 \oplus {\frak S}_2       & :=
	    & (\{F_{1; i}\oplus F_{2; i}\}_{i\in I}, \{h_{1; Z,  i j}\oplus h_{2; Z, ij}\}_{Z, ij})
		     \\[1.2ex]
	{\frak S}_1 \otimes_{\Bbb C} {\frak S}_2     & :=
	    & (\{F_{1; i}\otimes_{\Bbb C} F_{2; i}\}_{i\in I}, \{h_{1; Z,  i j}\otimes h_{2; Z, ij}\}_{Z, ij})
   \end{eqnarray*}
 define two new $C^\infty$ ${\Bbb C}$-vector-bundle settlements from the concordant ${\frak S}_1$ and ${\frak S}_2$.
 
\bigskip

\begin{definition}  {\bf [direct sum and tensor product]}\; {\rm
  Let ${\frak S}_1$ and ${\frak S}_2$ be concordant $C^\infty$ ${\Bbb C}$-vector-bundle settlements.
  Then,
     the $C^\infty$ ${\Bbb C}$-vector-bundle settlements
        ${\frak S}_1\oplus {\frak S}_2$  and ${\frak S}_1\otimes_{\Bbb C} {\frak S}_2$  as defined above
	 are called respectively the {\it direct sum} and the {\it tensor product} of ${\frak S}_1$ and ${\frak S}_2$.
}\end{definition}

\bigskip

\subsection{Sheaves associated to vector-bundle settlements}

Every vector bundle can be associated with its sheaf of sections. The same applies to bundle settlements.

\bigskip

\begin{flushleft}
{\bf Peripheral structures for a $C^\infty$ ${\Bbb C}$-vector-bundle settlement}
\end{flushleft}
Let ${\frak S}$ be a $C^\infty$ ${\Bbb C}$-vector-bundle settlement over $_{\cal P}Y$.
Continuing the terminology and notations in Definition~1.1.4   and  Definition~1.1.5.

\bigskip

\begin{definition} {\bf [bundle atlas for ${\frak S}$ from local distinguished blocks]}\; {\rm
  ${\frak S}$ is said to admit a {\it bundle atlas from local distinguished blocks}
   if for every $p\in Y$ there exists an open neighborhood $U_p$ of $p$ such that
     the clutching maps of bundle blocks around $p$ gives rise to a bundle embedding
	   ${\frak S}|_{U_p} \hookrightarrow F^\prime_p|_{U_p}$,
	 where $F_p$ is the distinguished bundle block at $p$ and $F^\prime_p$ is the extended bundle-block for $F_p$.
  $F^\prime_p|_{U_p}$ is called a {\it bundle chart} of ${\frak S}$ over a neighborhood of $p$.
 Such an ${\frak S}$ is called a {\it tame} $C^\infty$ ${\Bbb C}$-vector-bundle settlement.
}\end{definition}

\medskip

\begin{example} {\bf [local tensorial system is tame]}\; {\rm
 The local tensorial system ${\frak S}_{(H_{\frak U}, \otimes)}$, as constructed in Lemma/Definition~1.3.3,
  is a tame $C^\infty$ ${\Bbb C}$-vector-bundle settlement.
}\end{example}

\medskip

\begin{definition} {\bf [total space $|{\frak S}|$ of ${\frak S}$]}\;  {\rm
  The topological space $|{\frak S}|$ from gluing the bundle blocks $F_i\in {\frak S}$ by the clutching maps
    $\{h_{\tinybullet,\tinybullet\tinybullet}\}_{\tinybullet, \tinybullet\tinybullet}$ of ${\frak S}$
    is called the {\it total space} of ${\frak S}$.
  By construction, there is a built-in map $\pi^{\frak S} : |{\frak S}|\rightarrow {_{\cal P}Y}$
     with $(\pi^{\frak S})^{-1}(p) = F_p|_p$.
 With slight abuse of notations, we write also $\pi^{\frak S}:{\frak S}\rightarrow {_{\cal P}Y}$.
 For ${\frak S}$ tame, $\pi^{\frak S}$ is $C^\infty$
   in the sense that $\pi^{\frak S}|_{U_p}$ is the restriction of the bundle projection map
     $F^\prime_p|_{U_p}\rightarrow U_p$ (which is $C^\infty$)
	  for all chart $F^\prime_p|_{U_p}$ in the bundle atlas of ${\frak S}$.
}\end{definition}

\bigskip

\begin{flushleft}
{\bf Natural sheaves associated to $C^\infty$ ${\Bbb C}$-vector-bundle settlements}
\end{flushleft}

\begin{definition} {\bf [$C^\infty$-homomorphism between vector-bundle settlements]}\; {\rm
 Let $_{{\cal P}_1\!}{\frak S}_1$ and $_{{\cal P}_2\!}{\frak S}_2$
     be two tame $C^\infty$ ${\Bbb C}$-vector-bundle settlements on $Y$
     with regular partitions ${\cal P}_1$ and ${\cal P}_2$ respectively.
 A {\it $C^\infty$-homomorphism}
   from  ${_{{\cal P}_1\!}{\frak S}_1}$ to ${_{{\cal P}_2\!}{\frak S}_2}$
  is a continuous map
    $f: |{_{{\cal P}_1\!}{\frak S}_1}| \rightarrow |{_{{\cal P}_2\!}{\frak S}_2}|$
	on topological spaces that satisfies
	 \begin{itemize}
	  \item[(1)] ({\it fiber to fiber linearly})\hspace{1.2em}
	     The following diagram commutes
		   $$
             \xymatrix{
			    |{_{{\cal P}_1\!}{\frak S}_1}| \ar[rr]^-{f}   \ar[rd]_-{\pi^{{\frak S}_1}}
				   & & |{_{{\cal P}_2\!}{\frak S}_2}| \ar[ld]^-{\pi^{{\frak S}_2}}  \\
               & Y & \hspace{3em},
              }			
		   $$
		 with the restriction of $f$ to each fiber of ${_{{\cal P}_1\!}{\frak S}_1}$
		         a ${\Bbb C}$-linear map to a fiber of  ${_{{\cal P}_2\!}{\frak S}_2}$.
	
	  \item[(2)] ($C^\infty$/smoothness)\hspace{1.2em}
      For each $p\in Y$, there exist  a neighborhood $U$,
		a bundle chart $F^\prime_{1, p}|_U$ of ${_{{\cal P}_1\!}{\frak S}_1}$ over $U$, and
	    a bundle chart $F^\prime_{2, p}|_U$ of ${_{{\cal P}_2\!}{\frak S}_2}$ over $U$
		such that
		  the restriction of $f$ to over $U$ extends to
		  a $C^\infty$ ${\Bbb C}$-vector-bundle map $\hat{f}_U$ between the bundle charts:
		 $$
          \xymatrix{
		    F^\prime_{1, p}|_U \ar[rr]^-{\hat{f}_U}   && F^\prime_{2, p}|_U \\			
			  \left.|{_{{\cal P}_1\!}{\frak S}_1}| \rule{0ex}{1em}\right|_U \ar[rr]^-{f|_U}
                     \ar@{^{(}->}[u]			
			       &&  \;\left.|{_{{\cal P}_2\!}{\frak S}_2}| \rule{0ex}{1em}\right|_U\,. \ar@{^{(}->}[u]                          				
		  }
		 $$
	 \end{itemize}
 With slight abuse of notations, we write also
   $f: {_{{\cal P}_1\!}{\frak S}_1}\rightarrow {_{{\cal P}_2\!}{\frak S}_2}$.
 
 Denote the set of $C^\infty$-homomorphisms
     from ${_{{\cal P}_1\!}{\frak S}_1}$ to ${_{{\cal P}_2\!}{\frak S}_2}$
     interchangeably by $C^\infty({_{{\cal P}_1\!}{\frak S}_1}, {_{{\cal P}_2\!}{\frak S}_2})$
	  or $\Hom({_{{\cal P}_1\!}{\frak S}_1}, {_{{\cal P}_2\!}{\frak S}_2})$.
 This is naturally a $C^\infty(Y)^{\Bbb C}$-module.	
 When
    ${_{{\cal P}_1\!}{\frak S}_1}
	    = {_{{\cal P}_2\!}{\frak S}_2}  = {_{\cal P}{\frak S}} $,
  $f$ is called a {\it $C^\infty$-endomorphism} of $_{\cal P}{\frak S}$ and
 the $C^\infty(Y)^{\Bbb C}$-module
     $\Hom({_{\cal P}{\frak S}}, {_{\cal P}{\frak S}})$ is denoted by $\End({_{\cal P}{\frak S}})$.
}\end{definition}

\bigskip

With the above fundamental definition in place, basic sheaves that appear in Algebraic Geometry
 can be generalized naturally and immediately to the case of  tame vector-bundle settlements.
 
\bigskip

\begin{example} {\bf [$\Homsheaf$-sheaf and $\Endsheaf$-sheaf]}\; {\rm
 Let $_{{\cal P}_1\!}{\frak S}_1$ and $_{{\cal P}_2\!}{\frak S}_2$
     be two tame $C^\infty$ ${\Bbb C}$-vector-bundle settlements on $Y$
     with regular partitions ${\cal P}_1$ and ${\cal P}_2$ respectively.
 The presheaf defined by the assignment 	
      $$
	     U\; \longmapsto \;  \Hom(
		   {_{{\cal P}_1\!}{\frak S}_1}|_U,  {_{{\cal P}_2\!}{\frak S}_2}|_U)\,,
	  $$
	 $U\subset Y$  open, is a sheaf of ${\cal O}_Y^{\,\Bbb C}$-modules over $Y$,
  denoted by  $\Homsheaf_{{\cal O}_Y^{\;\Bbb C}}(
  		   {_{{\cal P}_1\!}{\frak S}_1},  {_{{\cal P}_2\!}{\frak S}_2})$.	
 In particular,
   $\Endsheaf_{{\cal O}_Y^{\,\Bbb C}}({_{\cal P}{\frak S}})
        := \Homsheaf_{{\cal O}_Y^{\,\Bbb C}}
		      ({_{\cal P}{\frak S}},  {_{\cal P}{\frak S}})$.	
}\end{example}

\medskip

\begin{example}{\bf [sheaf of sections and dual sheaf]}\;  {\rm
 Let ${\frak S}$ be a tame $C^\infty$ ${\Bbb C}$-vector-bundle settlement on $_{\cal P}Y$.
 The {\it sheaf of $C^\infty$-sections} of ${\frak S}$ is defined by
    ${\cal S}_{\frak S}
	    := \Homsheaf_{{\cal O}_Y^{\,\Bbb C}}({\cal O}_Y^{\,\Bbb C},  {\frak S})$.
 The $C^\infty(Y)^{\Bbb C}$-module of $C^\infty$ global sections of ${\frak S}$ is by definition
   $C^\infty({\frak S}):=  \Hom({\cal O}_Y^{\,\Bbb C},  {\frak S})$.

 As the clutching maps in ${\frak S}$ are ${\Bbb C}$-vector-bundle embeddings,
     the set of dual ${\Bbb C}$-vector-bundle blocks from ${\frak S}$  does not give rise to
     a $C^\infty$ dual ${\Bbb C}$-vector-bundle settlement ${\frak S}^\vee$ of ${\frak S}$ functorially
	 unless ${\frak S}$ is a ${\Bbb C}$-vector bundle.
 However, one can still define the {\it dual sheaf} of ${\frak S}$
     to be the dual ${\cal O}_Y^{\,\Bbb C}$-module of ${\cal S}_{\frak S}$, namely,
   ${\cal S}_{\frak S}^\vee
   	    := \Homsheaf_{{\cal O}_Y^{\,\Bbb C}}({\frak S} , {\cal O}_Y^{\,\Bbb C})$.
}\end{example}

\medskip

\begin{definition} {\bf [$C^\infty$ Azumaya-type noncommutative ringed space associated to bundle settlement]}\; {\rm
  Let
    $Y$ be a $C^\infty$-manifold with corners, with the structure sheaf  denoted by ${\cal O}_Y$ and
    ${\frak S}= {_{\cal P}{\frak S}}$ be a tame $C^\infty$ ${\Bbb C}$-vector-bundle settlement on a regularly partitioned $Y$.
  Then, the ringed space $(Y, {\cal O}_Y^\nc:= \Endsheaf_{{\cal O}_X^{\,\Bbb C}}({\cal S})_{\frak S} )$	
   is called the {\it $C^\infty$ Azumaya-type noncommutative space associated to ${\frak S}$}.
  Here,
  $C^\infty$ refers to the fact that the center ${\cal O}_Y^{\Bbb C}$ of ${\cal O}_Y^\nc$
     is the complexification of a sheaf of $C^\infty$-rings;
  {\it Azymaya-type} refers to the fact that ${\cal O}_Y^\nc$ is a sheaf of Azumaya algebras when restricted to over the interior
     of a block of $_{\cal P}Y$.
}\end{definition}

\bigskip

We conclude this subsection with the following open end:

\bigskip

\begin{motif-question-project} {\bf [bundle settlement with additional structure]}\; {\rm	
  Similar to the case of vector bundles, one can
      endow the bundle blocks $F_i$ with additional structures
	   (e.g.\ {\it Hermitian vector-bundle blocks}, {\it bundle-with-connection blocks})    and
	  require the clutching maps $h_{\tinybullet, \tinybullet\tinybullet}$ to be compatible with them
	      (e.g.\ {\it Hermitian vector-bundle embeddings}, {\it distribution-preserving bundle-embeddings}).
  One may also consider {\it principal-bundle settlements}  or {\it fortified principal-bundle settlements},
       their {\it associated vector-bundle settlements},     and
	   the corresponding gauge theories with matters
	         as the gluing of a collection of gauge theory with matters over a collection of  manifolds with corners.
  (Cf.\ Motif/Question/Project~1.1.7 and Motif/Question/ Project~1.1.8.)			
}\end{motif-question-project}

\bigskip

\subsection{Example: The local tensorial system from a seed system}

%
Presented here
 is a simple way to construct a $C^\infty$ ${\Bbb C}$-vector-bundle settlement,
 which will be employed in Sec.~2.
Let $Y$ be a $C^\infty$ $n$-manifold with corners,
  with its interior denoted by $\Int Y$ and its boundary $\partial Y:= Y-\Int Y$.

\bigskip

\begin{definition} {\bf [regular covering]}\; {\rm
  A finite open covering $\mathfrak{U}:= \{U_j\}_{j\in J}$ of $Y$ is called {\it regular}
	 if
	  \begin{itemize}
	   \item[(1)]
   	   the closure $\overline{U_j}$ of each $U_j$ is an embedded submanifold with corners in $Y$
	            such that
				 \begin{itemize}
				   \item[(1.{\it a})]
				      each $\overline{U_j}\cap \Int Y$ is a smooth manifold with boundary
			            (i.e.\ the codimension $\ge 2$ corners of $\overline{U_j}$, if nonempty, come solely from
				                    $\partial \overline{U_j}\cap \partial Y$)      and
									
				  \item[(1.{\it b})]
     				  each codimension-$k$ corner point of $\overline{U_j}$, $k\ge 2$, is
			             either a codimension-$k$ corner point of $Y$			
						 or from the transverse intersection of $\overline{\partial U_j\cap\Int Y}$
						              with a codimension-$(k-1)$-facet in $\partial Y$,
			   (in particular,  $\partial U_j = (\overline{\partial U_j\cap \Int Y})\amalg (U_j\cap \partial Y)$
						                 with $\overline{\partial U_j\cap \Int Y}$ and $\overline{U_j\cap \partial Y}$
											        meeting transversely  and
												 all the facets of $\overline{\partial U_j \cap \Int Y} $ of dimension $\le (n-2)$ lies
												     in $\overline{\partial U_j\cap \Int Y} \cap \overline{U_j\cap \partial Y}$);
			    \end{itemize}
				and						
				
	   \item[(2)]
 	   the collection of the $(n-1)$-dimensional strata in $\overline{\partial U_j\cap \Int Y}$, $j\in J$, and $\partial Y$
		     intersect transversely at where they intersect..
      \end{itemize}			
 Cf.$\hspace{.7ex}${\sc Figure} 1-3-1.
                  %
}\end{definition}
	
\medskip
			
\begin{lemma-definition} {\bf [from regular covering to regular partition]}\;
 A regular finite open covering $\mathfrak{U}$ of $Y$ determines a regular partition ${\cal P}_\mathfrak{U}$ of $Y$.
 {\rm  ${\cal P}_{\mathfrak{U}}$  is called the {\it regular partition of $Y$ associated to $\mathfrak{U}$}.
               (Cf.\ {\sc Figure}~1-3-1.)}
\end{lemma-definition}
	
\medskip

\noindent{\it Proof}. An explicit  construction of ${\cal P}_\mathfrak{U}$ is given as follows.
   Continuing the notion in the current subsection, define a characteristic map $\chi: Y\rightarrow 2^J$ by the assignment
     $p \mapsto  J_p\subset J$, where $J_p\, := \{j\in J\, |\, p\in U_j \}$.
   Since $\mathfrak{U}$ is a regular finite open covering,
     $J_p\subset J_{p^\prime}$ for all $p^\prime$ in a small enough neighborhood of $p$    and
    for $J^\prime\subset J$ such that 	$\chi^{-1}(J^\prime)$	is nonempty,
    the closure $Y_{J^\prime}:= \overline{\chi^{-1}(J^\prime)}$
	    is a possibly disconnected, embedded submanifold with corners in $Y$ of the same dimension as $Y$.
   The interior of $Y_{J^\prime}$ consists of points in $Y$ whose small enough neighborhoods are all mapped to $J^\prime$
     while the boundary $\partial Y_{J^\prime}$ of $Y_{J^\prime}$ are those points in $Y$ that are mapped to $J^\prime$
	   but always have a point in their neighborhood mapped to some $J^{\prime\prime}\ne J^\prime$.
   Let ${\cal P}_\mathfrak{U}:=\{Y_i\}_{i\in I}$ be the collection of connected components of $Y_{J^\prime}$'s.
   Then,
     the transversality of $\overline{\partial U_j \cap \Int Y}$, $j\in J$, among themselves and with $\partial Y$
	 implies that ${\cal P}_\mathfrak{U}$ is a regular partition of $Y$.
	
  \noindent\hspace{40.8em}$\square$

  \begin{figure}[htbp]
   \bigskip
    \centering
      \includegraphics[width=0.80\textwidth]{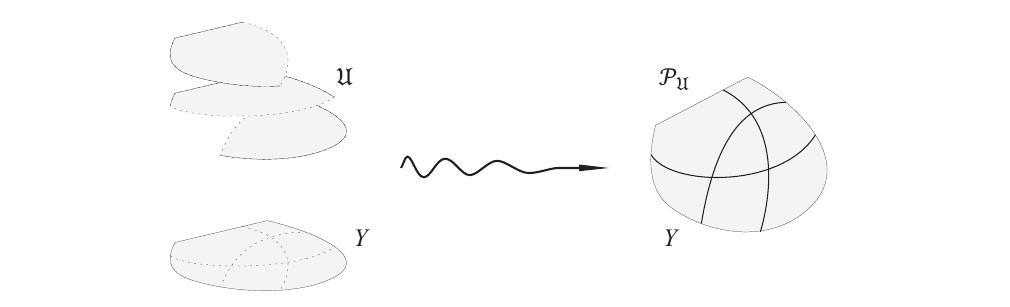}
 
    \bigskip
    \centerline{\parbox{13cm}{\small\baselineskip 12pt
    {\sc Figure}~1-3-1.
     From a regular covering ${\frak U}$ of $Y$ to a regular partition ${\cal P}_{\frak U}$ of $Y$.
    }}
  \end{figure}

\bigskip
	
Let
  $\mathfrak{U}=\{U_j\}_{j\in J}$ be a regular finite covering of $Y$   and
  $H_j$ be a $C^\infty$ ${\Bbb C}$-vector bundle on $\overline{U_j}$
    (i.e.$\hspace{.7ex}H_j$ is the restriction of
       a $C^\infty$ ${\Bbb C}$-vector bundle $H_j^\prime$ on an open set $U_j^\prime\supset \overline{U_j}$).
	
 \medskip

 \noindent
 {\bf Assumption [existence of nowhere-vanishing section]}\;\;
{\it Assume that, for all $j\in J$,
         $H_j$ has a $C^\infty$-section $s_j$ that is nowhere-vanishing over
		                        $\overline{U_j}\cap \bigcup_{j^\prime\ne j}\overline{U_{j^\prime}}$.}
 
 \medskip

\noindent
For convenience, put an order on $J$, i.e.\ an identification $J=\{1, 2, \cdots, |J|\}$ with the natural order $1<2<\cdots< |J|$.
Then
  the collection $H_\mathfrak{U}:= \{\mbox{$H_j$ over $\overline{U_j}$}\}_{j\in J}$ determines
  a $C^\infty$ ${\Bbb C}$-vector-bundle settlement on $_{{\cal P}_\mathfrak{U}\!}Y$,
     where ${\cal P}_{\mathfrak{U}}= \{Y_i\}_{i\in I}$ is the regular partition of $Y$ from Lemma/Definition~1.3.2,
  as follows:
   \begin{itemize}
    \item[$\tinybullet$]
	{\it $[C^\infty$ ${\Bbb C}$-vector bundle blocks$\,]$}\hspace{1em}
	  From the proof of  Lemma/Definition~1.3.2,
	    for each $Y_i$,
		 there is a unique $J(i)=\{j_1 < \cdots < j_l\}\subset J$ such that $Y_i$ is a component of $Y_{J(i)}$. 	
	  Let
	    $$
		    F_i\;  :=\;   H_{j_1}|_{Y_i}\otimes \cdots \otimes H_{j_l}|_{Y_i}\,.
		$$
  
	\item[$\tinybullet$]
	{\it $[C^\infty$-gluing data$\,]$}\hspace{1em}
	   For $Z\subset Y_i\cap Y_{i^\prime}$ a common (codimension-$1$) face of $Y_i$ and $Y_{i^\prime}$,
	     the regularity of $\mathfrak{U}$ implies that there is a unique $j\in J$ such that $Z\subset \partial \overline{U_j}$.
	   Then
        $$		
	     \mbox{either\:\:\:\: `\,$Y_i\subset Y-U_j$  and  $\Int Y_{i^\prime}\subset U_j$\,'
		  \:\:\:\:or\:\:\:\:           `\,$Y_{i^\prime}\subset Y-U_j$ and $\Int Y_i\subset U_j$\,'.}
        $$		
	  (Equivalently, either $j\in J(i^\prime)$ or $j\in J(i)$.)
       Up to relabelling of $i$ and $i^\prime$. we may assume the former.
	   Then assign the direction $\overrightarrow{ii^\prime}$   and
	     note that $J(i^\prime)= J(i)\cup \{j\}$ as sets and
		         that, as ordered sets, $J(i^\prime)$ is the insertion of $j$ to $J(i)$ at a unique site:
				 $ j_1 < \cdots < j_{l^{\prime\prime}}<j < j_{l^{\prime\prime}+1}<\cdots < j_l$.
       Define the clutching map $h_{Z, ii^\prime}: F_i|_Z\rightarrow F_{i^\prime}|_Z$
	      to be the built-in ${\Bbb C}$-vector-bundle embedding
		  \begin{eqnarray*}
		    \lefteqn{H_{j_1}|_Z\otimes \cdots \otimes H_{j_l}|_Z\;
			                   \stackrel{\sim}{\longrightarrow }
			                          H_{j_1}|_Z\otimes \cdots \otimes H_{j_{l^{\prime\prime}}}|_Z
			                             \otimes s_j|_Z \otimes
			                                H_{j_{l^{\prime\prime}+1}}|_Z \otimes \cdots \otimes H_{j_l}|_Z} \\			
               && \hspace{10em}
			        \subset
                        H_{j_1}|_Z\otimes \cdots \otimes H_{j_{l^{\prime\prime}}}|_Z
				           \otimes H_j|_Z \otimes
				              H_{j_{l^{\prime\prime}+1}}|_Z \otimes \cdots \otimes H_{j_l}|_Z\,.
		  \end{eqnarray*}
	 \end{itemize}
	
\begin{lemma-definition}	 {\bf [local tensorial system]}\;
  The above data defines a $C^\infty$ ${\Bbb C}$-vector-bundle settlement
     $\mathfrak{S}_{(H_\mathfrak{U}, \otimes)}$ over $_{{\cal P}_\mathfrak{U}\!}Y$,
  {\rm called the} {\it local tensorial system generated by $H_\mathfrak{U}$}.
\end{lemma-definition}

\medskip

\noindent{\it Proof.}
 Note that, by construction, all the clutching maps $h_{\tinybullet, \tinybullet\tinybullet}$ are $C^\infty$ around the faces.
 We give a detailed description of the resulting quivers for the strata of $_{\cal P}Y$,
    from which it is immediate that the consistency of gluing over each stratum is automatically satisfied.
 
 Let $Z$ be a codimension-$k$ facet of ${\cal P}$, $k\ge 2$.
 Then associated to $Z$ is the subset $J(Z)\subset J$ of cardinality $k$ such that
   $\partial\overline{U_j}$, $j\in J(Z)$, intersect transversely along $Z$.
 When crossing $\overline{\partial \overline{U_j}\cap \Int Y}$ from $Y-U_j$ to $U_j$,
  the fiber of the bundle block gets an extra tensor factor $E_j$.
 It follows that the clutching maps $h_{Z, \tinybullet\tinybullet}$ around $Z$ follows exactly the inward-normal direction of
    the boundary $\overline{\partial \overline{U_j}\cap \Int Y}$ of $U_j$, for $j\in J(Z)$.
 The associated quiver is the $1$-skeleton of the normal cube $I^k_{p.\varepsilon}$, for any $p\in \Int Z$,
   with the coordinate axes following the inward-normal direction of the corresponding
                       $\overline{\partial \overline{U_j}\cap \Int Y}$ of $U_j$, for $j\in J(Z)$, and
   the direction of an edge following the axis to which it is parallel to.
 This quiver has $2^k$-many vertices, $k\cdot 2^{k-1}$-many edges,
     a unique source and a unique sink,  and any vertex lies in an edge-path that goes from the source to the sink;
	 cf.$\hspace{.7ex}${\sc Figure}~1-3-2.
	         %
			 %
  Thus, in particular,
    it is already saturated and there is no need to add any more directed edges to it.
  \begin{figure}[htbp]
   \bigskip
    \centering
      \includegraphics[width=0.80\textwidth]{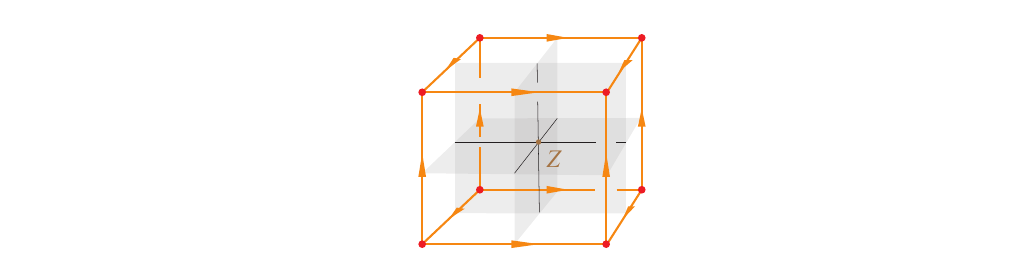}
 
    \bigskip
    \centerline{\parbox{13cm}{\small\baselineskip 12pt
    {\sc Figure}~1-3-2.
     The clutching quivers of a local tensorial system always have a sink and a source.
	  Illustrated is the clutching quiver (with red vertices and orange directed edges) of a codimension-$3$ stratum $Z$ (in brown),
	    with all the tangential directions to $Z$ suppressed.
	  The adjacent codimension-$1$ and codimension-$2$ strata to $Z$ are indicated (in gray and black respectively).	
    }}
  \end{figure}
  
  Given an edge-path of the quiver, let
    $v_i$ and $v_f$ be its initial and final vertices,
    $F_{v_i}$	 and $F_{v_f}$ be the bundle blocks associated to these vertices,
	$J(v_i)$ and $J(v_f)$ be the subsets of $J$ associated to these blocks with the induced order.
 Then $J(v_i)\hookrightarrow J(v_f)$	as ordered sets.
 This inclusion together with the restrictions $s_j|_Z$, for $j\in J(v_f)-J(v_i)$,
    specifies a natural embedding
	 $$
	   h_{Z;\, v_iv_f}\, :\,  F_{v_i}|_Z\; \longrightarrow\; F_{v_f}|_Z
	 $$ of ${\Bbb C }$-vector bundles
	 by mapping each tensor factor of $F_{v_i}|_Z$ to its identical factor in $F_{v_f}|_Z$ via the identity map and
	      inserting $s_j|_Z$ to the corresponding missed factor for $j\in J(v_f)-J(v_i)$.
 By construction, the composition of clutching maps along an edge-path with the initial vertex $v_i$ and the final vertex $v_f$
  corresponds exactly to $h_{Z;\, v_iv_f}$.
 This depends only on the initial and the final blocks and, hence, the consistency of gluing is always satisfied.			
   
 This completes the proof.
  
\noindent\hspace{40.8em}$\square$   

\bigskip
 
\begin{example}{\bf [natural sheaves associated to local tensorial system]}\; {\rm
 Continuing Example~1.2.2.
 Since a local tensorial system  ${\frak S}:= \mathfrak{S}_{(H_\mathfrak{U}, \otimes)}$ is tame,
   associated to  $\mathfrak{S}_{(H_\mathfrak{U}, \otimes)}$ are, for examples,
   its {\it sheaf of $C^\infty$-sections} ${\cal S}_{\frak S}$,
   $C^\infty(Y)^{\Bbb C}$-{\it module of global sections} $C^\infty({\frak S})$,
   {\it endomorphism sheaf} $\Endsheaf_{{\cal O}_Y^{\,\Bbb C}}({\frak S})$.
}\end{example}

\medskip

\begin{definition} {\bf [seed system]}\; {\rm
  The collection $H_{\frak U}$ of bundles with a distinguished section over a regular finite covering ${\frak U}$
    is called a {\it seed system} on $Y$
   as it generates (i.e.$\hspace{.7ex}$``grows to") a $C^\infty$ Azumaya-type noncommutative ringed space,
   following Example~1.3.4. ({\sc Figure}~1-3-3.)
}\end{definition}

\begin{figure}[htbp]
   \bigskip
    \centering
      \includegraphics[width=0.80\textwidth]{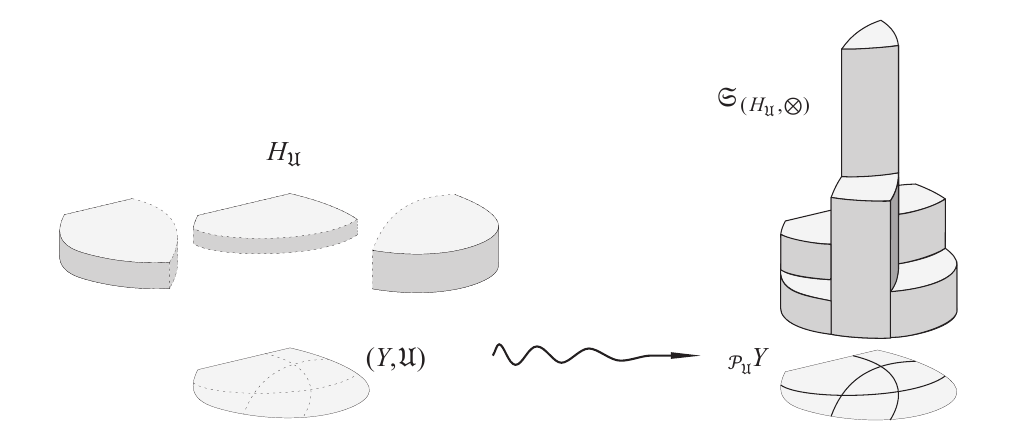}
 
    \bigskip
    \centerline{\parbox{13cm}{\small\baselineskip 12pt
    {\sc Figure}~1-3-3.
     A seed system $H_{\frak U}$ over a regular covering ${\frak U}$ of $Y$ gives rise to a bundle settlement over
	     ${_{{\cal P}_{\frak U}}}\!Y$, where ${\cal P}_{\frak U}$ is the regular partition of $Y$ determined by ${\frak U}$.
    }}
\end{figure}

\bigskip

\subsection{$C^\infty$ Azumaya-type noncommutative ringed space from a mildly singular $C^\infty$-bundle settlement
        on a mildly singular $C^\infty$-manifold with corners}
     %

With a generalization of the notion of a {\it local chart} of a manifold with corners (resp.\ a {\it local trivialization} of a vector bundle)
     so that the notion of {\it $C^\infty$ maps} between manifolds with corners
	     (resp.\ {\it $C^\infty$-homomorphisms} between vector bundles) can be defined in a sensible way,
 the notion of a $C^\infty$ ${\Bbb C}$-vector-bundle settlement over a $C^\infty$-manifold with corners can be generalized to
  a singular $C^\infty$ ${\Bbb C}$-vector-bundle settlement over a singular $C^\infty$-manifold with corners.
In this subsection, we discuss a class of such singular objects that will be used in later sections.

\bigskip

\begin{remark} {\it $[$chart, regularity, and locally-finitely-generatedness$]$}\; {\rm
 Though the notion of smoothness (equivalently, $C^\infty$) in Differential Geometry
     are originally for maps between open sets in ${\Bbb R}^l$'s,
   one learns from Algebraic Geometry  and Moduli Problems that one may still address such a notion for singular spaces
    by introducing appropriate local charts;
   cf.$\hspace{.7ex}$Kuranishi structures, tangent-obstruction complexes, finite generation of rings and modules.
 The setting below follows the same line of thought.
}\end{remark}

\bigskip

\begin{flushleft}
{\bf Mildly singular $C^\infty$ ${\Bbb C}$-vector-bundle settlements and related sheaves}		
\end{flushleft}
\begin{definition} {\bf [mildly singular $C^\infty$-manifold with corners]}\; {\rm
 A topological space $Y$ is called a {\it mildly singular $C^\infty$ $n$-manifold with corners}
  if every non-manifold-nor-corner point of $Y$ is isolated and the following conditions are satisfied:
  (Such a point is called a {\it singularity} of $Y$.
     Denote the set of such points by $Y_\scriptsizesing$,
                 the set of corner points of $Y$ by $\partial Y$,  and
                $Y_\smoothscriptsize:= Y-Y_\scriptsizesing$,
                $\Int Y:= Y-\partial Y$.)
	 \begin{itemize}
	  \item[$\cdot$]
	  $Y_\smoothscriptsize$ is a $C^\infty$ $n$-manifold with corners  and
	  $Y_\scriptsizesing$ is discrete and closed.
		
      \item[$\cdot$]
	  ${\partial Y}$ is closed. (Equivalently, $\Int Y$ is open.)		
	  Thus, $Y_\scriptsizesing$ is contained in $\Int Y$.
		
	  \item[$\cdot$]  [{\it local chart at singularity}\,]\hspace{2em}
	  Every $p\in Y_\scriptsizesing$ has a fixed open neighborhood $U_p\subset \Int Y$
      such that
	     $U_p$ is topologically a cone over a connected closed $(n-1)$-manifold with the apex $p$
	      (and hence contains no other singularities of $Y$)   and that
		 there is an embedding $\alpha_p: U_p\hookrightarrow {\Bbb R}^{n_p}$, for some $n_p> n$,
		    that is $C^\infty$ on $U_p-\{p\}$. {\sc Figure~1-4-1}
			                                               %
	 \end{itemize}
 A function $f$ defined on an open set $U\subset Y$ is said to be $C^\infty$ (synonymously, {\it smooth})
   if $f$ is smooth on $U-Y_\scriptsizesing$ and
      the restriction $f|_{V\cap U_p}$ for each $p\in Y \cap Y_\scriptsizesing$ is the restriction of a $C^\infty$ function
	    on a neighborhood of $\alpha_p(V\cap U_p)$ in ${\Bbb R}^{n_p}$ under $\alpha_p$.
 Denote the sheaf of $C^\infty$-functions on $Y$ still by ${\cal O}_Y$  and name it the {\it structure sheaf} of $Y$.
}\end{definition}
	
\begin{figure}[htbp]
  \bigskip
   \centering
   \includegraphics[width=0.80\textwidth]{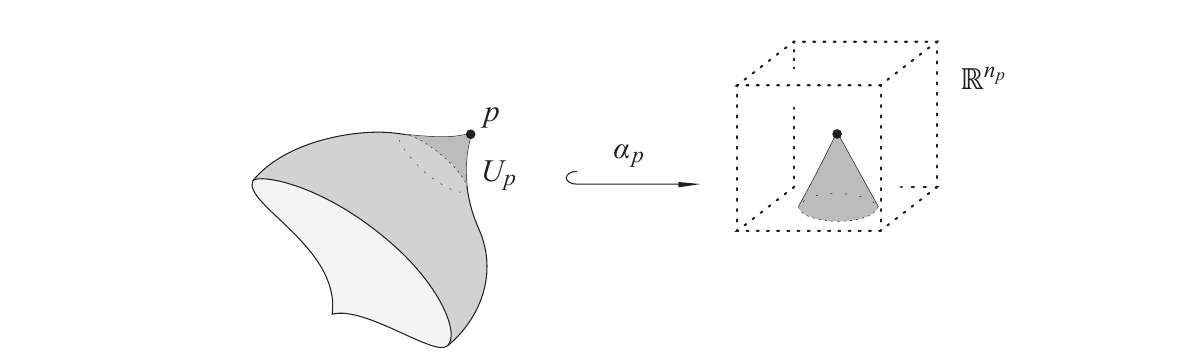}
 
  \bigskip
 \centerline{\parbox{13cm}{\small\baselineskip 12pt
  {\sc Figure}~1-4-1.
    A local chart  at an isolated singularity $p$ in the interior of a manifold $Y$ with corners.
	Similar to the role it plays for a regular point, a chart at $p$ determines the $C^\infty$-ring of local smooth functions around $p$
	  and hence the notion of a map from or to $Y$ smooth at $p$.
	 Equivalently, all the local charts together determines the structure  sheaf ${\cal O}_Y$ of $Y$ as a $C^\infty$-scheme
	   and hence the notion of a morphism to or from $Y$ from or to other $C^\infty$-schemes.
  }}
\end{figure}

\medskip

\begin{definition} {\bf [regular partition]} \; {\rm
  (Cf.$\hspace{.7ex}$Definition~1.1.3.)
 Let $Y$ be a mildly singular $C^\infty$ $n$-manifold with corners.
 A decomposition $Y=\bigcup_{i\in I}Y_i$ of $Y$
     into a collection of connected mildly singular $C^\infty$ $n$-manifolds $Y_i$ with corners
  is called a {\it regular partition $\cal P$} of $Y$ if ${\cal P}$ induces a regular partition on $Y_\smoothscriptsize$.
 In particular,\;\;
   $Y_\scriptsizesing\;
       \subset\; Y - \amalg_{k=1}^n{_{\cal P}^{\;k}Y}\;
	   =\; \amalg_{i\in I}\Int Y_i$.
}\end{definition}

\medskip

\begin{definition} {\bf [mildly singular $C^\infty$ ${\Bbb C}$-vector bundle]}\; {\rm
 Let $Y$ be a mildly singular $C^\infty$-manifold with corners.
 A {\it mildly singular $C^\infty$ ${\Bbb C}$-vector bundle} on $Y$
    is a sheaf ${\cal F}$ of ${\cal O}_Y^{\,\Bbb C}$-modules
  such that
    ${\cal F}|_{Y_{smooth}}$ is locally free of finite rank
      (equivalently, ${\cal F}|_{Y_{smooth}}$ is the sheaf of sections of a $C^\infty$ ${\Bbb C}$-vector bundle
          $F$ on $Y-Y_\scriptsizesing$)
	and that, for $p\in Y_\scriptsizesing$, there is an open $U_p^\prime\subset U_p$ such that
	    ${\cal F}_{U_p^\prime}$ is finitely generated.
 Will denote the whole singular bundle also by $F$ and set
    the global-Hom $\Hom(F_1, F_2):= \Hom_{{\cal O}_Y^{\,\Bbb C}}({\cal F}_1, {\cal F}_2)$  and
	the sheaf-Hom $\Homsheaf(F_1, F_2):=\Homsheaf_{{\cal O}_Y^{\,\Bbb C}}({\cal F}_1, {\cal F}_2)$
	   for mildly singular $C^\infty$ ${\Bbb C}$-vector bundles $F_1$ and $F_2$ on $Y$.		
}\end{definition}

\bigskip

Recall Definition~1.1.4.
Since, by design,   the set $Y_\scriptsizesing$ of singularities of a mildly singular $C^\infty$-manifold with corners lies in the
   interior of the blocks of a regular partition ${\cal P}$ of $Y$,
 the following definition follows immediately:
 
\bigskip

\begin{definition} {\bf [mildly singular $C^\infty$-bundle settlement]}\; {\rm
 Let $Y$ be a mildly singular $C^\infty$-manifold with corners  and ${\cal P}=\{Y_i\}_{i\in I}$ be a regular partition of $Y$.
 A gluing system
     $_{\cal P}{\frak S}= (\{F_i\}_{i\in I},
	                                                   \{h_{\tinybullet, \tinybullet\tinybullet}\}_{\tinybullet, \tinybullet\tinybullet}) $
	 over $Y$,
     defined in the same way as in Definition~1.1.4
	    but with each $C^\infty$ ${\Bbb C}$-bundle block $F_i$
		  generalized to a mildly singular $C^\infty$ ${\Bbb C}$-vector bundle,
	  is called a {\it mildly singular $C^\infty$ ${\Bbb C}$-vector-bundle settlement} over $Y$.
}\end{definition}

\medskip

\begin{example} {\bf [from Algebraic Geometry$/{\Bbb C}$ and Complex Geometry]}\; {\rm
  Algebraic varieties over ${\Bbb C}$ and analytic spaces with only Gorenstein isolated singularities
      are examples of mildly singular $C^\infty$-manifolds
  while torsion-free coherent sheaves over the above spaces that are locally free over the regular/smooth part of the spaces
      are examples of mildly singular $C^\infty$ ${\Bbb C}$-vector bundles. 	
}\end{example}

\medskip

\begin{definition} {\bf [regular covering \& seed system]}\; {\rm
 (Cf.$\hspace{.7ex}$Definition$\hspace{.7ex}$1.3.1, Definition$\hspace{.7ex}$1.3.5.)
 Let $Y$ be a mildly singular $C^\infty$-manifold with corners.
 A finite open covering  ${\frak U}:= \{U_j\}_{j\in J}$ of $Y$ is called {\it regular}
   if (1) each $p\in Y_\scriptsizesing$ is contained in a unique $U_{j_p}$ and in such a way that
                $p\in U_{j_p}-\bigcup_{j\ne j_p}\overline{U_j}$; and
      (2) ${\frak U}$ restricts to a regular covering of $Y-Y_\scriptsizesing$.
 For such a ${\frak U}$,
   a collection $H_{\frak {U}}:=\{H_j\}_{j\in J}$,
      where $H_j$ is a mildly singular $C^\infty$ ${\Bbb C}$-vector bundle on $\overline{U_j}$
	  (i.e.$\hspace{.7ex}H_j$ is the restriction of a mildly singular $C^\infty$ ${\Bbb C}$-vector bundle $H_j^\prime$
	           on an open set $U_j^\prime\supset \overline{U_j}$)
	   such that, for all $j\in J$,
         $H_j$ has a $C^\infty$-section $s_j$ that is nowhere-vanishing over
		                        $\overline{U_j}\cap \bigcup_{j^\prime\ne j}\overline{U_{j^\prime}}$
  is called a {\it seed system} on $Y$.
}\end{definition}

\bigskip

Given a seed system $H_{\frak U}$ on $Y$, the same construction in Sec.\ 1.3
 applies to give
   a {\it local tensorial system} ${\frak S}:= {\frak S}_{(H_{\frak U}, \otimes)}$,
        cf.\ Lemma/Definition~1.3.3, and then
   a $C^\infty$ Azumaya-type noncommutative ringed space
      $Y^\nc:=(Y, {\cal O}_Y^{\nc}:=\Endsheaf_{{\cal O}_{Y}^{\,\Bbb C}}({\cal S}_{\frak S}))$,
	  cf.\ Example~1.3.4.

\bigskip

\begin{flushleft}
{\bf Morita equivalence and ${\cal O}_Y^\nc$-modules }
\end{flushleft}
The fundamental (left) action of $\End_{\Bbb C}({\Bbb C}^l)$ on ${\Bbb C}^l$ induces
   a right action of $\End_{\Bbb C }({\Bbb C})$ on the dual $({\Bbb C}^l)^\vee$ and
  a natural  ${\Bbb C}$-vector-space isomorphism
  $({\Bbb C}^l)^\vee \otimes_{End_{\,\Bbb C}({\Bbb C}^l)} {\Bbb C}^l \rightarrow {\Bbb C}$,
    $f\otimes v \mapsto f(v)$.
On the other hand, one has a natural ${\Bbb C}$-vector-space  isomorphism
  ${\Bbb C}^l\otimes_{\Bbb C}({\Bbb C}^l)^\vee \rightarrow \End_{\Bbb C}({\Bbb C}^l)$,
   $v\otimes f\mapsto v\cdot f(\tinybullet)$.
These induce an adjoint-pair of equivalences
 $$
   \xymatrix{
   {\Bbb C}\mbox{-}\Mod \ar@<+.6ex>[rrrr]^-{{\Bbb C}^l\otimes_{\Bbb C}(\tinybullet)}
        &&&& \End_{\Bbb C}({\Bbb C}^l)\mbox{-}\Mod \ar@<+.3ex>
		         [llll]^-{({\Bbb C}^l)^\vee\otimes_{End_{\,\Bbb C}({\Bbb C}^l)} (\tinybullet)}\,.
   }
 $$
This is the most elementary form of Morita equivalence between two categories of modules of different rings/algebras
(here ${\Bbb C}$ versus $\End_{\Bbb C}({\Bbb C}^l)$),
which generalizes immediately to ${\Bbb C}$-vector bundles over a manifold with corners:
  \begin{itemize}
   \item[$\cdot$]
   {\it Natural $C^\infty$ ${\Bbb C}$-vector-bundle isomorphisms}
  $$
    F^\vee \otimes_{End_{\,\Bbb C}({\Bbb C}^l)}F \;\longrightarrow\;  \underline{\Bbb C}_Y
  $$
   and
  $$
   F\otimes_{\Bbb C}F^\vee\; \longrightarrow\;  \End_{\Bbb C}(V)\,
  $$
   where
     $F$ is a $C^\infty$ ${\Bbb C}$-vector bundle of rank $l$ over a $C^\infty$-manifold $Y$ with corners,
     $F^\vee$ the dual ${\Bbb C}$-vector bundle of $F$, and
	 $\underline{\Bbb C}_Y$ is the trivial ${\Bbb C}$-vector bundle of rank $1$ over $Y$
	 (whose sheaf of $C^\infty$-sections is ${\cal O}_Y^{\,\Bbb C}$);
   
  \item[$\cdot$]
  which give rise to an {\it adjoint-pair of equivalences} of categories of sheaves of modules
  $$
   \xymatrix{
   {\cal O}_Y^{\,\Bbb C}\mbox{-}\Mod \ar@<+.6ex>
         [rrrr]^-{{\cal F}\otimes_{{\cal O}_Y^{\,\Bbb C}}(\:\tinybullet\:)}
        &&&& \Endsheaf_{{\cal O}_Y^{\;\Bbb C}}({\cal F})\mbox{-}\Mod \ar@<+.3ex>
		         [llll]^-{{\cal F}^\vee
				                 \otimes_{{\cal E}nd_{\,{\cal O}_Y^{\,\Bbb C}}({\cal F})} (\:\tinybullet\:)}\,,
    }
  $$
  where ${\cal F}$ is the sheaf of $C^\infty$-sections of $F$ and ${\cal F}^\vee$ is the dual sheaf.
 \end{itemize}
     %

Such an equivalence between categories of modules is called {\it Morita equivalence}.
Readers are referred to [Le: Sec.\ 2] and [Mo: Sec.\ 3] for more thorough discussions.

Given a mildly singular $C^\infty$ ${\Bbb C}$-vector-bundle settlement ${\frak S}:=\{F_i\}_{i\in I}$
  over a regularly partitioned mildly singular $C^\infty$-manifold with corners,
  let ${\cal S}_{\frak S}$ be its sheaf of $C^\infty$-sections  and
  ${\cal O}_Y^\nc:= \Endsheaf_{{\cal O}_Y^{\,\Bbb C}}({\cal S}_{\frak S}) $
     is the structure sheaf for the $C^\infty$ Azumaya-type noncommutative ringed space $Y^\nc$.
Then each $F_i$ is a $C^\infty$ ${\Bbb C}$-vector bundle on $Y_i-Y_\scriptsizesing$,
       where $Y_i\in {\cal P}$ is the defining block of $F_i$,  and
   overall ${\cal S}_{\frak S}$ is locally free on $Y_0 := Y-(Y_\scriptsizesing\cup  {\,_{\cal P}^{\,\ge 1}Y})$.
Thus, Morita equivalence relates
  $\Endsheaf_{{\cal O}_{Y_i}^{\,\Bbb C}}({\cal F}_i) |_{Y_i-Y_{sing}}$-$\Mod$
     to ${\cal O}_{Y_i-Y_{sing}}$-$\Mod$;
  and ${\cal O}_{Y_0}^\nc$-$\Mod$ to ${\cal O}_{Y_0}$-$\Mod$.	
(Here, ${\cal F}_i$ be the sheaf of $C^\infty$-sections of $F_i$ over $Y_i$.)

\bigskip

\begin{motif-question-project} {\bf [full Morita equivalence]}\; {\rm
 Let $Y$ be a mildly singular $C^\infty$-manifold with corners.
 Work out criterions on Gorenstein singularities of $Y$ and (presumably combinatorial) conditions on ${\frak S}_{\cal P}$
  so that a full Morita equivalence ${\cal O}_Y^\nc$-$\Mod\rightleftarrows{\cal O}_Y$-$\Mod$ holds.
}\end{motif-question-project}

\medskip

\begin{remark} {$[$locally finite regular covering$]$}\; {\rm									
 Before proceeding to the next section, we remark that
  the notion of `regular covering' and related constructions in this section generalize to {\it locally finite regular coverings}.
}\end{remark}

\medskip

\begin{remark} {$[$bundle settlement vs.\ soft noncommutative scheme\,$]$}\; {\rm
 The notion of bundle settlement is based on a notion of generalized gluing over $Y$
     that gives rise to an Azumaya-type noncommutative ringed space $Y^\nc$ shadowing over $Y$
  while the notion of soft noncommutative schemes, as experimented in [L-Y7] (D(15.1), NCS(1))	  and [L-Y8] (NCS(2)),
    is based on another notion of generalized gluing and leads a noncommutative space $Y^{\nc \prime}$ into which $Y$ is embedded.
 For the local affine charts, the former uses extensions of commutative algebras to noncommutative algebras while the latter realizes
   commutative algebras as quotients of noncommutative algebras.
 This is the origin of the difference of {\it shadowing over} vs.\ {\it embedded into} feature
   of $Y^\nc$ vs.\ $Y^{\nc \prime}$ with respect to $Y$.
 The string-theoretic implication of this difference via D-branes remains to be understood.
 Cf.\ {\sc Figure} 1-4-2.
           %
 %
 \begin{figure}[htbp]
  \bigskip
   \centering
   \includegraphics[width=0.80\textwidth]{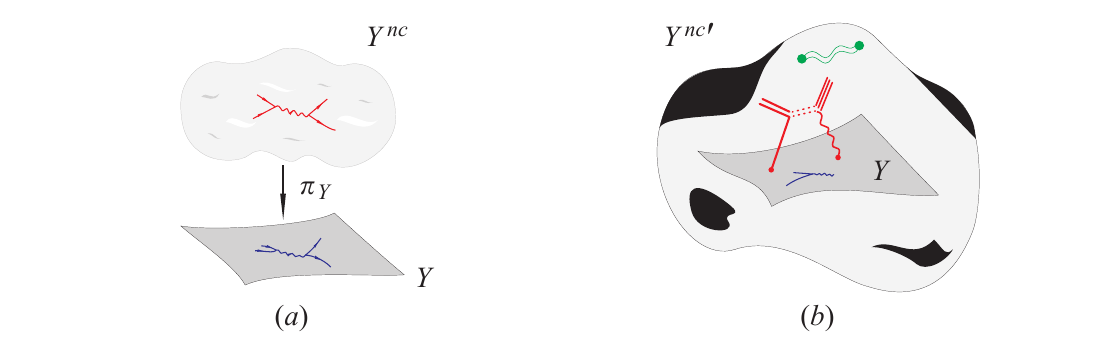}
 
  \bigskip
 \centerline{\parbox{13cm}{\small\baselineskip 12pt
  {\sc Figure}~1-4-2.
  Two classes of noncommutative schemes/ringed-spaces as generalizations of commutative schemes:
  ({\it a}) {\it noncommutative clouds $Y^{nc}$ over a commutative scheme} $Y$, as in the current notes; and
  ({\it b}) {\it soft noncommutative schemes} $Y^{nc\, \prime}$
                   from  {\it extension of local commutative affine charts of a commutative scheme $Y$ to noncommutative charts}  and
                   gluing of charts via inclusions to gluing of charts via morphisms (i.e.\ {\it soft gluing}),
				   as in [L-Y7]  (D(15.1), NCS(1)) and [L-Y8] (NCS(2)).
 The notion of morphisms from an Azumaya/endomorphism ringed-space (i.e.\ D-brane probe)
                   to a noncommutative space in either class can be naturally defined.
  For a physical world supported on a space in class ({\it a}),
    everything observed in the commutative world $Y$, both at the classical and the quantum levels,
     is a shadow of its parallel (classical or quantized)  in the noncommutative cloud $Y^{nc}$.
  In contrast, for a physical world supported on a space in class ({\it b}),
    as a commutative probe sees only the commutative subscheme $Y$(and its quantization)
     	of the noncommutative soft scheme $Y^{nc\,\prime}$,
	 D-branes serve as mediators/messengers/sonar for an observer in the commutative world $Y$
	     to probe the noncommutative world $Y^{nc\,\prime}$ around.
	 In general one expects that there could be ``dark regions" in $Y^{nc\,\prime}$ that remain dark to $Y$
	 (i.e.\ cannot can be probed via scattering D-branes into them).
	Some scattering processes are illustrated in each class of targets.
 (There are other classes of noncommutative spaces, e.g.\ via graded algebras or categories, in the literature.
      If they admit the notion of morphisms from a matrix point into them, then they are all potential target spaces for dynamical D-branes.
	  Such details/possibilities have not yet been explored.)\\
  $\mbox{\hspace{1.2em}}$It
  is noteworthy that for a noncommutative space $Y^{nc}$ in Class ({\it a}), $Y^{nc}$ remains to have the topology $Y$
    while for $Y^{nc \prime}$ in Class ({\it b}) the topology is gone in general.
  The latter reminds one of quantum topology in quantization of gravity theory
	   though the connection of the two phenomena is not clear to us.
	 }}
 \end{figure}					
}\end{remark}

\bigskip

\section{Seed systems over a Gorenstein singular Calabi-Yau space}

We review in Sec.\ 2.1 the notion and examples of (local) noncommutative crepant resolutions    and
 use them in Sec.\ 2.2 to create seed systems on a Calabi-Yau space with Gorenstein isolated singularities.

\bigskip

\subsection{Noncommutative crepant resolution of a Gorenstein isolated singularity and the associated apical algebra}

In this subsection, we'll stay in the realm of ordinary Algebraic Geometry as in [Hart].

Let
 $Y=\Spec R$ be an affine variety associated to an ${\Bbb C}$-algebra $R$ such that
       $Y$ is smooth except  at an isolated Gorenstein singularity $y_\ast$;
   $I_{y_\ast}\subset R$ (resp.$\hspace{.7ex}{\cal I}_{y_\ast}\subset {\cal O}_Y$)
      be the ideal (resp.$\hspace{.7ex}$ideal sheaf) associated to $y_\ast$.
  
\bigskip

\begin{definition} {\bf [noncommutative crepant resolution]}\; {\rm ([VdB2: Definition~4.1].)
  A {\it noncommutative crepant resolution} of $R$ is a homologically homogeneous $R$-algebra of the form
  $A =\End_R(M)$ where $M$ is a reflexive $R$-module.
 (Cf.\ Details in  [VdB2: Sec.~3 \& Sec.~4] of {\sl Michel Van den Bergh}; see also [Le: Sec.\ 12], [B-O: Sec.\ 5], [We: Sec.\ 2].)
}\end{definition}

\bigskip
 
\noindent
Denote the coherent ${\cal O}_Y$-module associated  to the $R$-module $A$ by ${\cal A}$    and
let $Y^\nc:=(Y, {\cal A})$ be the associated noncommutative ringed space.
  
\bigskip

\begin{definition}   {\bf [apical algebra of resolution of isolated singularity]}\; {\rm
 The (finite-dimensional) noncommutative  ${\Bbb C}$-algebra $A/(I_{y_\ast}A)$
      is called the {\it apical algebra} of the given noncommutative crepant resolution.
 It is the fiber ${\cal A}/({\cal I}_{y_\ast}{\cal A})$ of the ${\cal O}_Y$-algebra ${\cal A}$
   over the isolated singularity $y_\ast\in Y$ (i.e.$\hspace{.7ex}$ the fiber of $Y^\nc$ over $Y$ at $y_\ast$).
}\end{definition}

\medskip

\begin{remark} $[$existence of ${\Bbb C}$-point over $y_\ast\,]$\;  {\rm
 As we'll see through the conifold example below,
   the fiber of $Y^\nc$ over $Y$ at $y_\ast$ may contain a ${\Bbb C}$-point
     even when the general fiber doesn't.
}\end{remark}
  
\medskip

\begin{example} {\bf [ordinary double point/conifold]}\; {\rm
In this local model,
   $Y=\Spec R$, where
     $$
	    R={\Bbb C}[y_1, y_2, y_3, y_4]/(y_1y_3-y_2y_4)\,,
	 $$
	 and
   $I\subset R$ be the ideal $(y_1, y_2)$.
Then, as a hypersurface in ${\Bbb A}^4:= \Spec{\Bbb C}[y_1, y_2, y_3, y_4]$,
  $Y$ has an ordinary double point $q$ at the origin and
  $I$ is the ideal of $R$ associated to a Weil divisor $\simeq \Spec {\Bbb C}[y_3, y_4] \simeq {\Bbb A}^2$ on $Y$
    that contains the double point.
The blowing-up of $R$ along $I$ gives rise to a crepant resolution of $Y$ with the exceptional locus ${\Bbb P}^1$.
(Cf.\ {\sc Figure}~2-1-3-1.)
   %
 \begin{figure}[htbp]
  \bigskip
   \centering
   \includegraphics[width=0.80\textwidth]{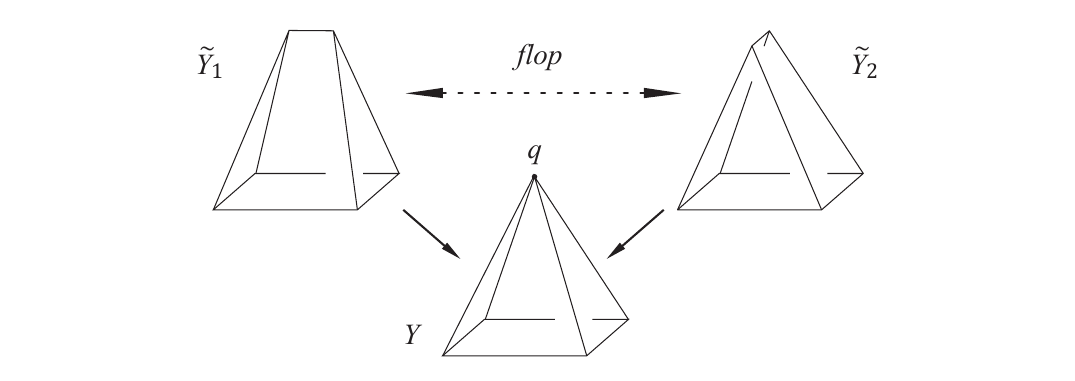}
 
  \bigskip
 \centerline{\parbox{13cm}{\small\baselineskip 12pt
  {\sc Figure}~2-1-4-1.
   (Commutative) crepant resolutions of a conifold $Y$, illustrated in terms of toric geometry.
   The two resolutions $\widetilde{Y}_i\rightarrow Y$, $i=1, 2$,  are related by a flop.
  }}
 \end{figure}
Set $\{1\}= \{1\}\amalg \emptyset$.
 Then
  in terms of the contravariantly equivalent language of rings and modules
                                      to affine schemes and quasi-coherent sheaves over affine schemes,
   we define
    the ${\cal O}_Y$-module
      ${\cal I}_{\{1\}}\oplus {\cal I}_{\emptyset}$ by $I \oplus R$;     and
    the ${\cal O}_Y$-algebra
      ${\cal A}_{(\{1\}, \emptyset)}$ by the $R$-algebra
      $$
          \End_R(I \oplus R)\;
	         \simeq \left(
	                  \begin{array}{lc}
					    R & I \\ I^{-1} & R
					  \end{array}
	                     \right)\,,
      $$
	  which acts on  $I\oplus R$ (with elements treated as column vectors) from the left.
     Here, letting $K(R):= R[(R-\{0\})^{-1}]$ be the total quotient ring of $R$,    then
      $I^{-1}:= \{s\in K(R)| sI\subset R\}$ a fractional ideal of $R$.
     Explicitly, $I^{-1}=    y_1^{-1}(y_1, y_4)= y_2^{-1}(y_2, y_3)$.
 The $R$-algebra $\End_R(I \oplus R)$
  gives a {\it noncommutative crepant resolution} of $R$.
 Equivalently in geometrical terms,
  it is the function ring of an (algebraic) {\it noncommutative smooth Calabi-Yau $3$-space} over $Y$.
 Cf.$\hspace{.7ex}$[VdB2: proof of Proposition 7.3].
  
 Some immediate observations follow:
  \begin{itemize}
   \item[(1)]	[{\it as quotient of an algebra supported on a free module}$\,$]\hspace{1em}
  The surjective $R$-module-homomor- phisms
       $R\otimes y_1 \oplus R\otimes y_2 \rightarrow I$,
	    with $r_1\otimes y_1 + r_2\otimes y_2 \mapsto r_1y_1+r_2y_2$ in $R$, and\\
       $R\otimes 1 \oplus R \otimes (y_1^{-1}y_4)\rightarrow I^{-1}$,
        with $r_3\otimes 1 + r_4\otimes (y_1^{-1}y_4)\mapsto r_3+r_4\cdot (y_1^{-1}y_4)$ in $K(R)$,	
               from the generators of $I$ and $I^{-1}$
     and the multiplication of $I$ and $I^{-1}$ in $K(R)$
   give rise to a surjective $R$-algebra-homomorphism
   $$
    \left(
	  \begin{array}{cc}
	    R  & R\otimes y_1 \oplus R \otimes y_2 \\
		R \otimes 1 \oplus R\otimes (y_1^{-1}y_4)	& R
	  \end{array}
	\right)\;
	 \longrightarrow \;
	  \left(
	     \begin{array}{lc}
		     R & I \\ I^{-1} & R
		 \end{array}
	 \right)\,.	
   $$	
  This realizes $ \End_R(I\oplus R)$ as a quotient algebra of an $R$-algebra with an underlying rank-$6$ free $R$-module.

  \item[(2)]  [{\it fiber over ordinary double point}$\,$]\hspace{1em}
 (2.1)
  ${\cal I}_{\{1\}}\oplus {\cal I}_\emptyset$ is locally free of rank $2$ on $Y_\ast := Y-\{q\}$,
   with the fiber over $q$ isomorphic to ${\Bbb C}^3$.
 
 (2.2)
  ${\cal A}_{(\{1\},\emptyset)}$ is a sheaf of rank-$2$ Azumaya/matrix algebras on $Y_\ast$
    with the fiber ${\Bbb C}$-algebra over $q$
	         isomorphic to the $6$-dimensional ${\Bbb C}$-algebra
																		 in $2\times 2$ matrix form
	$$
	  \left(
	   \begin{array}{cc}
	    {\Bbb C}   & {\Bbb C}\cdot\underline{y}_1\oplus {\Bbb C}\cdot\underline{y}_2  \\[1.2ex]
	    {\Bbb C}\cdot \underline{1}\oplus {\Bbb C}\cdot (\underline{y}_1^{-1}\underline{y}_4)
		                    & {\Bbb C}
	   \end{array}
	  \right)\,,
	$$
	where 	
	 $$
	  \begin{array}{lcl}
	   \underline{1}\cdot \underline{y}_1  =  \underline{y}_1\cdot \underline{1}=0\,,	
	       && (\underline{y}_1^{-1}\underline{y}_4)\cdot \underline{y}_1
	                   = \underline{y}_1\cdot (\underline{y}_1^{-1}\underline{y}_4)=0\,, \\[1.2ex]		      
	   \underline{1}\cdot \underline{y}_2  =  \underline{y}_2\cdot \underline{1}=0\,,			
		   &&  (\underline{y}_1^{-1}\underline{y}_4)\cdot \underline{y}_2
	                        = \underline{y}_2\cdot (\underline{y}_1^{-1}\underline{y}_4)=0\,.	
	  \end{array}
	 $$	
	
  \item[(3)]  [{\it projective module/irreducible representation}$\,$]\hspace{1em} 		
   ${\cal A}_{(\{1\}, \emptyset)}$ has two projective modules\\   $R\oplus I^{-1}$ and $I\oplus R$
    while its fiber at $q$ has
	 two $3$-dimensional irreducible representations
	 $$
	  \left[
	   \begin{array}{c}
	    {\Bbb C}   \\[1.2ex]
	    {\Bbb C}\cdot \underline{1}\oplus {\Bbb C}\cdot (\underline{y}_1^{-1}\underline{y}_4)
	   \end{array}
	  \right]\,,\;\;
	  \left[
	   \begin{array}{c}
	    {\Bbb C}\cdot\underline{y}_1\oplus {\Bbb C}\cdot\underline{y}_2  \\[1.2ex]
	    {\Bbb C}
	   \end{array}
	  \right]\,,
	$$
    and two ${\Bbb P}^1$-families of equivalent irreducible $1$-dimensional representations
 	$$
	  \left[
	   \begin{array}{c}
	    0   \\[1.2ex]
	    {\Bbb C}\cdot
		  (a\cdot \underline{1}+b \cdot (\underline{y}_1^{-1}\underline{y}_4))
	   \end{array}
	  \right]\,,\;\;
	  \left[
	   \begin{array}{c}
	    {\Bbb C}\cdot
		    ( a\cdot \underline{y}_1 +  b\cdot\underline{y}_2 )  \\[1.2ex]
	     0
	   \end{array}
	  \right]\,,
	$$
	$(a, b)\ne (0,0)$,
	which are associated respectively to the ${\Bbb C}$-algebra-epimorphisms
    $$
	  \left(
	   \begin{array}{cc}
	    {\Bbb C}   & {\Bbb C}\cdot\underline{y}_1\oplus {\Bbb C}\cdot\underline{y}_2  \\[1.2ex]
	    {\Bbb C}\cdot \underline{1}\oplus {\Bbb C}\cdot (\underline{y}_1^{-1}\underline{y}_4)
		                    & {\Bbb C}
	   \end{array}
	  \right)\;
	  \longrightarrow \;
	  \left(\begin{array}{cc} 0 & 0 \\ 0 & {\Bbb C}\end{array}\right)
	$$
	and
	$$
	  \left(
	   \begin{array}{cc}
	    {\Bbb C}   & {\Bbb C}\cdot\underline{y}_1\oplus {\Bbb C}\cdot\underline{y}_2  \\[1.2ex]
	    {\Bbb C}\cdot \underline{1}\oplus {\Bbb C}\cdot (\underline{y}_1^{-1}\underline{y}_4)
		                    & {\Bbb C}
	   \end{array}
	  \right)\;
	  \longrightarrow \;
	  \left(\begin{array}{cc} {\Bbb C} & 0 \\ 0 & 0 \end{array}\right)\,.	  	
	$$	
    %
  \end{itemize}
\noindent\hspace{40.8em}$\square$
}\end{example}

\medskip

\begin{example}  {\bf [quotient singularity]}\; {\rm
 (Cf.$\hspace{.7ex}$[VdB2: Example 1.1] and [Le: Sec.$\hspace{.7ex}$11].)
 In this local model,
  there is a finite subgroup $G \subset \SL_n({\Bbb C})$ that acts linearly and effectively
    on the ${\Bbb C}$-vector space ${\Bbb C}^n$.
 Fix  a basis $(y^1, \,\cdots\,, y^n)$ of ${\Bbb C}^n$,
   there is then an induced action on the polynomial ring
     $S={\Bbb C}[y^1,\,\cdots\,, y^n] = \Sym^\tinybullet ({\Bbb C}^n)$ over ${\Bbb C}$
     and hence on ${\Bbb A}^n_{\Bbb C}:= \Spec ({\Bbb C}[y^1,\,\cdots\,, y^n])$.
 Under the setting, $R=S^G$ the subring of $G$-fixed elements of  $S$ and the inclusion $R\subset S$ defines
     the quotient map $\rho: {\Bbb A}^n_{\Bbb C}\longrightaarrow Y=\Spec R$ of the $G$-action
	 on ${\Bbb A}^n_{\Bbb C}$.
 The (maximal) ideal associated to the singularity $y_\ast\in Y$ is given by
   $I_{y_\ast}={\frak m}_{y_\ast}=R-({\Bbb C}-\{0\})$.
 Denote the ${\cal O}_Y$-module associated to the $R$-module $S$ by ${\cal S}$,
   which is $\rho_\ast {\cal O}_{{\Bbb A}^n_{\Bbb C}}$.
 Since $G$ acts on ${\Bbb A}^n_{\Bbb C}-\{\mathbf{0}\}$  freely,
  ${\cal S}$ is locally free over $Y-\{y_\ast\}$ of rank the order $|G|$ of the group $G$.
 (Cf.\ {\sc Figure}~2-1-5-1.)
     %
     %
 \begin{figure}[htbp]
   \bigskip
    \centering
    \includegraphics[width=0.80\textwidth]{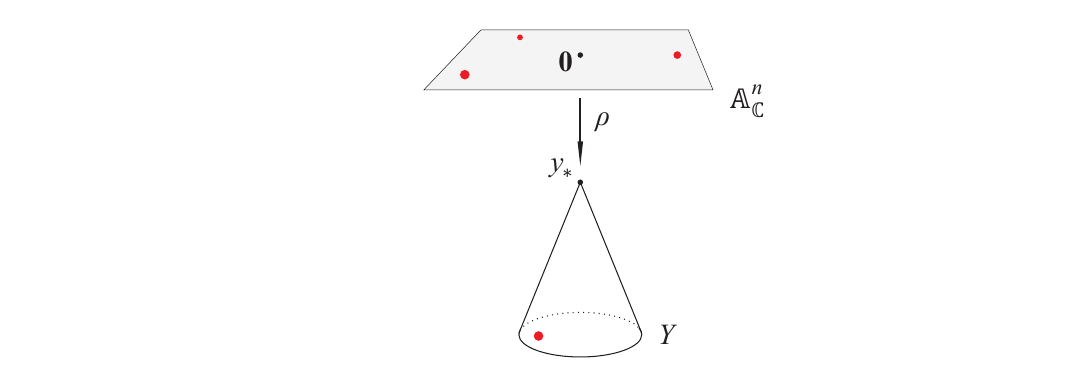}
 
   \bigskip
  \centerline{\parbox{13cm}{\small\baselineskip 12pt
   {\sc Figure}~2-1-5-1.
    Isolated quotient singularities.
   }}
 \end{figure}
      
 A noncommutative crepant resolution of $Y$ is given by the {\it endomorphism algebra}
  $$
      A\; =\; \End_R(S)\,, \hspace{1em}\mbox{where $S$ is regarded as a (left) $R$-module}\,.
  $$
 Denote the associated noncommutative ringed space by
    $Y^\nc := (Y, {\cal O}_Y^\nc:= {\cal A}:= \Endsheaf_{{\cal O}_Y}({\cal S}))$.
 There is an isomorphism $\gamma: S\rtimes G \stackrel{}{\longrightarrow} A$,
    defined by $\gamma (s,g)(s^\prime)= s\,g(s^\prime)$,
    where $S\rtimes G$  is the {\it twisted group algebra} $S\rtimes G$,
   whose ring-multiplication is defined by $(s_1, g_1)\cdot (s_2, g_2) = (s_1\,g_1(s_2), g_1g_2)$.
 The $R$-module-homomorphism
    $$
	  \pi\; :\;  S\; \longrightarrow\;  R\,, \hspace{2em}
                     s\;\longmapsto\;  \frac{1}{|G|}\sum_{g\in G} g(s)\,
    $$	
    leads explicitly to the following lemma
   \begin{itemize}
    \item[] \parbox[t]{38em}{{\bf Lemma~2.1.5.1.  [$R$ is a direct summand of $S$]}\; {\it
    There is an $R$-module projection map $\pi: S\rightarrow R$ that restricts to the identity map on $R$ .
	 It follows that, as $R$-modules, $S= R\oplus \pi^{-1}(0)=: R\oplus S^\prime$.
      }} 
   \end{itemize}
  Consequently, $A$ has a presentation in a block $2\times 2$ matrix form:
   $$
      A\; =\; \End_R(S)\;  =\;
	    \left[
		  \begin{array}{cc}
		   \End_R(R, R))           & \Hom_R(S^\prime, R) \\
		   Hom_R(R, S^\prime)  & \End_R(S^\prime)
		  \end{array}
		 \right]\;                        =\;
	   \left[
		  \begin{array}{cc}
		    R                     & (S^\prime)^\vee \\
		    S^\prime       &   \End_R(S^\prime)
		  \end{array}
		 \right]\,,	
   $$
   with the product of the $2\times 2$ matrix-entries defined by the composition of $R$-module-homomor-phisms:
      $a_{ij}a_{jk}:= a_{ij}\circ a_{jk}$.
   \begin{itemize}
    \item[] \parbox[t]{38em}{{\bf Lemma~2.1.5.2.
	            [sufficient condition for existence of ${\Bbb C}$-point over $y_\ast$]}\; {\it
	 If the $R$-module $S^\prime$ does not contain $R$ as a direct summand, then
	  the apical algebra $A/(I_{y_\ast}A)$ has a ${\Bbb C}$-algebra quotient ${\Bbb C}$.
      }} 
	
	\medskip
	\item[] \parbox[t]{38em}{{\it Proof}.\;
	Since $R$ is trivially projective as an $R$-module,
	 any epimorphism $S^\prime\rightarrow P$ would give rise to a split exact sequence
	   $0 \rightarrow S^{\prime\prime} \rightarrow S^\prime \rightarrow R \rightarrow 0$  and, hence,
	   realizes $R$ as a direct summand of $S^\prime$.
    It follows that, if $R$ is not a direct summand of $S^\prime$,
	  then every $R$-module-homomorphism $S^\prime\rightarrow R$ takes values in $R-({\Bbb C}-\{0\})$,
	  which is exactly the maximal ideal ${\frak m}_{y_\ast}=I_{y_\ast}$ of $y_\ast$.
    In this case, all the compositions in
       $\Hom_R(S^\prime, R)\times \Hom_R(R, S^\prime)\rightarrow \End_R(R, R)=R$,
            $(a_{12}, a_{21})\mapsto a_{12}\circ a_{21}$,
	   take values in ${\frak m}_{y_\ast}$.
	 Consequently, the $R$-module-epimorphism
	    $$
		   A\; \longrightarrow\;  R\,, \hspace{2em}
		    \left[ \begin{array}{ll}a_{11} & a_{12} \\ a_{21}  &  a_{22}\end{array}\right]\;
		    \longmapsto\;
			\left[ \begin{array}{cc} a_{11}   & 0 \\   0  &  0  \end{array} \right]\,,		
		$$
      	descends to a ${\Bbb C}$-algebra quotient after modding out by ${\frak m}_{y_\ast}$.
	  This proves the lemma.\\
	     $\mbox{\hspace{36em}}\square$
	   }
   \end{itemize}
   Since ${\cal A}$ is  locally free over $Y-\{y_\ast\}$, with fiber the simple algebra $\GL_{|G|}({\Bbb C})$,
     a ${\Bbb C}$-point of $Y^\nc$, if exists, can only occur over $y_\ast$.
   \begin{itemize}
     \item[]\parbox[]{38em}{{\bf Example 2.1.5.3. [$A_1$ surface singularity]}\;
	  Let $S={\Bbb C}[y^1, y^2]$ and $G={\Bbb Z}/2$ act on $S$ by $(y^1, y^2)\mapsto (-y^1, -y^2)$.
	  Then $R=S^G$ consists of polynomials, with ${\Bbb C}$-coefficients, in variables $y^1$ and $y^2$ with
	      all their  monomial terms of even total degree  and
	 $\Spec R$ gives a local model of the $A_1$ surface singularity.
    In this case
       $S=R\oplus S^\prime$, where $S^\prime$	consists of polynomials, with ${\Bbb C}$-coefficients,
	      in variables $y^1$ and $y^2$ with  all their  monomial terms of odd total degree;
      and the direct-sum decomposition is simply the grouping of even-total-degree terms and odd-total-degree terms.
	Let $Y=\Spec R$. Then, as an ${\cal O}_Y$-module, the coherent sheaf ${\cal S}$ on $Y$ associated
	    to the $R$-module $S$ has the generic fiber of dimension $2$ and the fiber
	    at the singularity $y_\ast\in Y$ of dimension $3$.
	It follows that $S^\prime$ cannot have $R$ as a direct summand, else $S^\prime$	would have torsion elements.
    Thus, $\End_R(S)$ has a ${\Bbb C}$-algebra quotient ${\Bbb C}$;
	 equivalent;ly, $Y^\nc=(Y, \Endsheaf_{{\cal O}_Y}({\cal S}))$
	   contains a ${\Bbb C}$-point over $y_\ast$.
           } 
   \end{itemize}	
\noindent\hspace{40.8em}$\square$
}\end{example}

\bigskip

\begin{example} {\bf [Gorenstein affine toric Calabi-Yau threefold]}\;{\rm
  (Cf.\ [Br], [Ken], [Gu], [St].)
  The field theory on the world-volume $X$ of stacked D$3$-branes
       sitting at the singularity $y_\ast$ of an affine toric Calabi-Yau threefold
	       $Y/{\Bbb C}=\Spec R$ with Gorenstein isolated singularity
	   in the Type IIB Superstring Theory compactification on $Y$
     is given by a superconformal quiver gauge theory with a superpotential $W$
	   (a formal noncommutative polynomial with matter superfields in the world-volume field theory as noncommuting variables).
  The $10$-dimensional target space-time in this case is ${\Bbb M}^4\times Y  $,
        where ${\Bbb M}^4$ is the $4$-dimensional Minkowski space-time,   and
	$X$ is embedded in ${\Bbb M}^4\times Y$ with image ${\Bbb M}^4\times\{y_\ast\}$.
  Given the supersymmetry type and kinds of superfields involved,
    the kinetic term for all superfields and the potential term for the gauge/vector super fields in the Lagrangian density
	  are completely determined and
	it is the superpotential $W$ that parameterizes most degrees of freedom of the related Wilson's theory-space
	(i.e.$\hspace{.7ex}$the
             moduli stack that parameterizes supersymmetric field theories of a given combinatorial type specified by the quiver).
 The data of the monomials that appear in the $W$ can be encoded in a bi-partite tiling of the (real $2$-)torus
 (named a {\it dimer model}, cf.\ {\sc Figure}~2-1-6-1).
      %
      %
 One can associate to the superpotential $W$ its {\it noncommutative Jacobian ring} $A$.
 \begin{figure}[htbp]
   \bigskip
    \centering
    \includegraphics[width=0.80\textwidth]{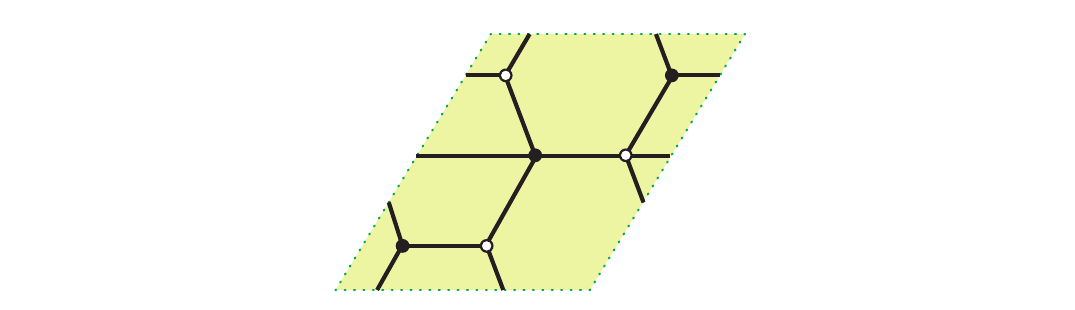}
 
   \bigskip
  \centerline{\parbox{13cm}{\small\baselineskip 12pt
   {\sc Figure}~2-1-6-1.
    A {\it dimer} (i.e.\ bi-partite unoriented graph) on the (topological, real $2$-)torus is indicated, following [Ken: Fig.\ 4].
	Here the torus is presented as an open fundamental domain of a ${\Bbb Z}^2$-action on ${\Bbb R}^2$.		
   }}
 \end{figure}
 
 Beginning with a dimer model,  {\sl Nathan Broomhead} shows in [Br: Chapters 7 \& 8] that
   \begin{itemize}
    \item[$\cdot$] \parbox[t]{38em}{\it When
	 the dimer model associated to the superpotential $W$ satisfies a certain algebraic consistency condition,
	$A$ is a Calabi-Yau $3$-algebra in the sense of {\rm [Gi2]}      and
	    is of the form $\End_{R^\prime}(M)$,
		where $R^\prime$ is the center of $A$ and $M$ is a reflexive $R^\prime$-module.}
   \end{itemize}			
 This renders $A$ a noncommutative crepant resolution of $Y^\prime := \Spec R^\prime$.
 (From stringy resolution of singularities, one would expect $R^\prime=R$ one begins with in the compactification
     of Type IIB string theory on $Y=Spec R$.)
 Together with the work [Gu] of {\sl Daniel Gulotta}	and [St] of {\sl Jan~Stienstra}
      that associate a geometrically consistent (which implies algebraically consistent) dimer model
               to any given Gorenstein affine toric threefold, Broomhead proves:
	 \begin{itemize}
	  \item[$\cdot$] \parbox[t]{38em}{[Br: Theorem 8.6]\, {\bf [existence of NCCR via dimer model]} \; {\it
	   Every Gorenstein affine toric threefold admits a noncommutative crepant resolution, which can be obtained
	   via a geometrically consistent dimer model.
      }} 
	 \end{itemize}
  
  In [Br], a toric-geometry description of $A$ from its associated dimer model is given as follows:
  (subject only to mild notation changes)
   \begin{itemize}
    \item[(1)]
	 Given a dimer model, i.e.$\hspace{.7ex}$a bipartite graph on the (real $2$-)torus $T$ with the vertices colored
	   by, say, either black or white,
	 assume that the dimer model is {\it balanced}
	 -- meaning that the number of black vertices is the same as the number of white vertices --
	 and let $Q$ be the dual quiver of the bipartite graph in $T$, following an orientation convention.
	Denote the set of vertices of $Q$ (resp.$\hspace{.7ex}$directed edges of $Q$, oriented facets in $T-Q$) by $Q_0$
	(resp.$\hspace{.7ex}Q_1$, $Q_2$).
	This gives a polygonal decomposition of $T$.
	
	\item[(2)]
	On one hand, one can read off the superpotential $W$ from the dimer model on $T$ (or equivalently from $Q\subset T$)
	 and associate to the dimer the noncommutative Jacobian ring $A$ of $W$.
	
	\item[(3)]
	On the other hand,
	the boundary map $\partial$ on chains from the decomposition with integer coefficients defines a complex of
	  finitely generated Abelian groups:
	  $$
	     {\Bbb Z}_{Q_0}\; \stackrel{\partial}{\longlongleftarrow}\;
		   {\Bbb Z}_{Q_1}\; \stackrel{\partial}{\longlongleftarrow}\;
		   {\Bbb Z}_{Q_2}
	  $$
	  and the boundary map $d$ on cochains from the decomposition with integer coefficients defines the dual complex
	   $$
	    {\Bbb Z}^{Q_0}\; \stackrel{d}{\longlongrightarrow}\;
		   {\Bbb Z}^{Q_1}\; \stackrel{d}{\longlongrightarrow}\;
		   {\Bbb Z}^{Q_2}\,.
	   $$
	  Let $\underline{1}_2\in {\Bbb Z}^{Q_2}$  be the function on $Q_2$ that sends all the facets to $1$.
	  Then the balancedness condition on the dimer model  implies that $\underline{1}_2$ is exact.
	  Let  $N:= d^{-1}({\Bbb Z}\cdot \underline{1}_2)$.
	  Since $d\circ d =0$, one has a map of lattices
	   $$
	     {\Bbb Z}^{Q_0}\; \stackrel{d}{\longlongrightarrow} \; N\,.
	   $$
	 The kernel of this map is ${\Bbb Z}\cdot \underline{1}_0$, where $\underline{1}_0$ is the map on $Q_0$
	   that sends all the vertices to $1$.
	 For a strongly convex rational polyhedral cone $\sigma$ in $N\otimes_{\Bbb Z}{\Bbb R}\supset N$,
	   let $N^+_\sigma:= \sigma \cap N$.
	
	 \item[]\hspace{1.2em}
	 Consider the corresponding dual map of the dual lattices
	    $$
		   {\Bbb Z}_{Q_0}\; \stackrel{\partial}{\longlongleftarrow}\;
		     M:= N^\vee := \Hom_{\Bbb Z}(N, {\Bbb Z})
		$$
	 and the dual cone of $N^+_\sigma$
	    $$
		  M^+_\sigma\;  :=\;   (N^+_\sigma)^\vee\;
		     =\; \{u\in M\,|\, \langle u, v \rangle \ge 0,\, \mbox{for all $v\in N^+_\sigma$}\}\,.
		$$
 	
	 \smallskip
	 \item[] \parbox[t]{38em}{{\bf Definition 2.1.6.1. [noncommutative toric data \& noncommutative toric algebra]}\;
	 ([Br: Definition 5.1].)
	    Let
		  $$
		    \begin{array}{l}
 		     M_{ij}\; :=\;  \partial^{-1}(j-i)\; \mbox{for $i,j\in Q_0$}\,, \hspace{2em}
		        M^+_{\sigma, ij}\; :=\;  M_{ij}\cap M^+_\sigma\,, \hspace{1em}\mbox{ and}  \\[1.2ex]
		     M_\sigma^{+, \flat}\; :=\;
			    \mbox{the submonoid of $M^+_\sigma$ generated by $\bigcup _{i,j\in Q_0}M^+_{\sigma, ij}$.}
		   \end{array}
		  $$		
	  The tuple $(N, Q_0, d, \sigma)$ is called a {\it noncommutative toric data} if it satisfies the following two conditions
           \begin{itemize}
		     \item[(i)]
			 $M^+_{\sigma, ij}$ is non-empty for all $i,j\in Q_0$;
			
			 \item[(ii)]
			 $N^+_\sigma = (M_\sigma^{+, \flat})^\vee:
			    = \{v\in N\,|\, \langle u, v\rangle \ge 0\,, \, \mbox{for all $u\in M_\sigma^{+, \flat}$}\}$.
           \end{itemize}
		 Given such a data, one can associate a {\it noncommutative toric algebra}
		    $$
			  B_\sigma\;:=\; \left[
			        \begin{array}{ccc}
					  {\Bbb C}[M^+_{\sigma, 11}]    & \cdots   & {\Bbb C}[M^+_{\sigma, 1n}] \\
					      \vdots  & \cdots  & \vdots \\
				      {\Bbb C}[M^+_{\sigma, n1}]    & \cdots   & {\Bbb C}[M^+_{\sigma, nn}] 					
					\end{array}
			         \right]
			$$
			Here
			  $Q_0$ is presented as $\{1,\, \cdots\,,\, n\}$  and
			  the multiplication
			    ${\Bbb C}[M^+_{\sigma, ij}] \times {\Bbb C}[M^+_{\sigma, jk}]
				     \rightarrow {\Bbb C}[M^+_{\sigma, ik}]$
				is from the restriction of the multiplication
                    ${\Bbb C}[M^+_\sigma] \times {\Bbb C}[M^+_\sigma] \rightarrow {\Bbb C}[M^+_\sigma]$
				in the group ring ${\Bbb C}[M_\sigma^+]$.
		   Note that $M_{\sigma, ii}=\partial^{-1}(0)=: M_o$ and thus
		     ${\Bbb C}[M^+_{\sigma, ii}]$ 	is independent of $i\in Q_0$.
		  Denote it by ${\Bbb C}[M_{\sigma, o}^+]$.			
	 }
	
	\item[]	
   {\sc Figure} 2-1-6-2.
      %
      %
    \begin{figure}[htbp]
      \bigskip
       \centering
       \includegraphics[width=0.80\textwidth]{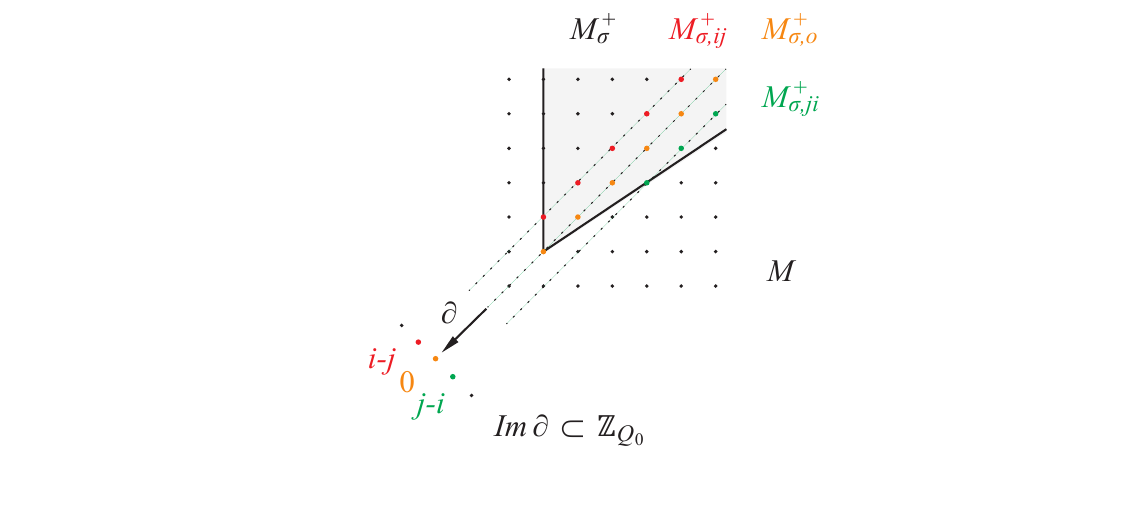}
 
      \bigskip
     \centerline{\parbox{13cm}{\small\baselineskip 12pt
      {\sc Figure}~2-1-6-2.
	    Part of a {\it noncommutative toric data} is indicated, following [Br: Example 5.4].
		In terms of the toric geometry language in [Fu],
		   a dual cone in the $M$ lattice,
		   a projection to a base lattice, and
		   an additively closed finite set  that contains the origin and lies in the image of the lattice map
		 determine a block-matrix ring.
      }}
    \end{figure}
		
	\smallskip
	\item[]\hspace{1.2em}
	 $B_\sigma$ has the following properties:
	 \begin{itemize}
	  \item[({\it a})]
	   {\it The center of $B_\sigma$ is ${\Bbb C}[M_{\sigma, o}^+]$.}
	
	  \item[({\it b})]
	   {\it $B_\sigma \otimes_{{\Bbb C}[M_{\sigma,o}^+]} {\Bbb C}[M_o]
	            = \Mat_n({\Bbb C}[M_o])\,$, the $n\times n$ matrix ring with entries in ${\Bbb C}[M_o]$.}
	
	  \item[({\it c})]
	   {\it There is a natural ${\Bbb C}$-algebra-homomorphism ${\frak h} : A\rightarrow B_\sigma$.}
	 \end{itemize}
	
	 \smallskip
	 \item[] \parbox[t]{38em}{{\bf  Definition 2.1.6.2. [algebraic consistency]}\; ([Br: Definition 5.12].)
	  A dimer model is called {\it algebraically consistent}
	    if the morphosm ${\frak h}$ in ({\it c}) above is an isomorphism.
	  } 
	
	\smallskip
	\item[(4)] \parbox[t]{38em}{[Br: Sec.$\hspace{.7ex}$8.1] {\bf [reflexivity]}\; {\it
	 Let $A\simeq B_\sigma$ be the algebra obtained from an algebraically consistent dimer model.
	 Then, for any $i\in Q_0$, there is a natural  ${\Bbb C}[M_{\sigma, o}^+]$-module-isomorphism
	   $$
	      {\Bbb C}[M_{\sigma, jk}^+]\;  \stackrel{\sim}{\longlongrightarrow}\;
		      \Hom_{{\Bbb C}[M_{\sigma, o}^+]}
			        ({\Bbb C}[M_{\sigma, ij}^+]\, ,\, {\Bbb C}[M_{\sigma, ik}^+])\,.
	   $$	
	  Consequently, ${\Bbb C}[M_{\sigma, jk}^+]$ is a reflexive ${\Bbb C}[M_{\sigma,o}^+]$-module
	    for all $j, k\in Q_0$.
	  }} 
	
	 \medskip
     \item[] \hspace{1.2em}	
	 Let $F_i := \oplus_{j\in Q_0}{\Bbb C}[M_{\sigma, ij}^+]$.
	 It follows that
	 {\it for an algebraically consistent dimer model,
	 $A\simeq B_\sigma \simeq \End_{{\Bbb C}[M_{\sigma, o}^+]}(F_i)$
	   is a noncommutative crepant resolution of $Y^\prime := \Spec {\Bbb C}[M_{\sigma, o}^+]$.}\\[1.2ex]
              $\mbox{\hspace{1.2em}}$
     The work [Gu] of Gulotta and [St] of Stienstra says that  {\it there exists a geometric consistent dimer model
        on the (real 2-)torus such that one can have $\sigma$ chosen (via taking the polygon of {\it perfect-matchings})
		so that ${\Bbb C}[M_{\sigma, o}]=R$.} 	
	 Thus, it can be arranged in the toric construction via a dimer model so that $Y^\prime = Y$.	
	 Together with [Br: Chapter 6\, {\it Geometric consistency implies algebraic consistency}],
	   this proves [Br: Theorem 8.6] quoted above.\\[1.2ex]
              	 $\mbox{\hspace{1.2em}}$
	 This concludes our highlight of [Br] (for the part we need).
  \end{itemize}
 
 \medskip
 \noindent
 It follows from the above highlight and the same reasoning as in Example~2.1.5 that
 %
   \begin{itemize}
	\item[]\parbox[t]{38em}{
	 {\bf  Corollary 2.1.6.3. [$R$ direct summand; existence of ${\Bbb C}$-point]}\;
	{\it The reflexive $R$-module $F_i$ admits a decomposition $F_i=R\oplus F_i^\prime$ as $R$-modules,
	            for all $i\in Q_0$.
	        Furthermore, if there exists an $i\in Q_0$ such that $F_i^\prime$ does not contain $R$ as a direct summand,
			  then the apical algebra of the noncommutative crepant resolution as constructed above via a dimer model
			  has a ${\Bbb C}$-point.}
			                                           } 
   \end{itemize}
\noindent\hspace{40.8em}$\square$ 		
}\end{example}

\smallskip

\begin{remark} $[$cone over del Pezzo surface\,$]$\; {\rm
 (Cf.$\hspace{.7ex}$[VdB2: Sec.$\hspace{.7ex}$7].)
  Let $Z$ be a del Pezzo surface (i.e.~a
      complex projective surface with ample anticanonical line bundle $\omega_Z^{-1}$).
  Then
   $Z$ is either ${\Bbb P}^2$,
   ${\Bbb P}^1\times {\Bbb P}^1$,   or
   a blowup of ${\Bbb P}^2$ at $k\le 8$ points. (Here, all projective spaces are over ${\Bbb C}$.)
  Let
   $Y:= \Spec (\oplus_{i=0}^\infty \varGamma (Z, {\cal L}^{\otimes i}))$
      be the  Gorenstein complex cone over $Z$,
	 where ${\cal L}$ is the ample line bundle ${\cal O}_Z(1)$
	                (resp.$\hspace{.7ex}{\cal O}_Z(1,1)$, ${\omega_Z^{-1}}$)
	    on $Z$  for $Z={\Bbb P^2}$
		   (resp.$\hspace{.7ex}{\Bbb P}^1\times {\Bbb P}^1$, else)
	 and
   $\widetilde{Y}:= \boldSpec (\Sym^{\tinybullet}{\cal L})$
     be the standard resolution of $Y$.
 Then there are natural  morphisms through the construction:
  $$
    \xymatrix{
	  & \widetilde{Y}\ar[rd]^-{\pi_Z} \ar[ld]_-{\pi_Y} \\
	  Y     && Z\,.
	}
	\hspace{6em}   \raisebox{-2em}{In illustration:}\hspace{1em}
	   \raisebox{-4em}{\includegraphics[width=0.30\textwidth]{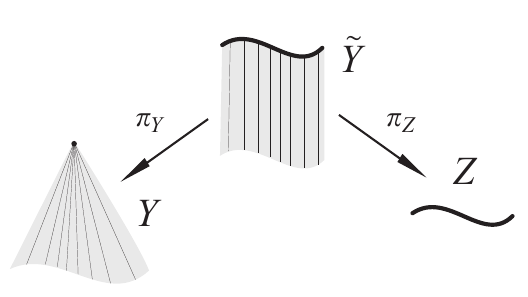}}
  $$
 For $Z={\Bbb P}^2$, $Y={\Bbb A}^3$ is regular.
 For  $Z={\Bbb P}^1\times {\Bbb P}^1$, $Y$ is the conifold in Example~2.1.4.
 Thus, we may assume here that $Z$ is neither ${\Bbb P}^2$ nor ${\Bbb P}^1\times {\Bbb P}^1$.
  In this case, {\sl Michel Van den Bergh} constructed a vector bundle $F_0$ on $Z$
  from an exceptional collection of vector bundles on $Z$ generating the derived category $D(\Qch(Z))$
  of quasi-coherent sheaves on $Z$ such that $F:= \pi_Z^\ast F_0$ is a tilting object on $\widetilde{Y}$.
 A noncommutative crepant resolution of $Y$ is then given by the endomorphism ring
   $\End_Z(F)$($=\End_{\widetilde{Y}}(F)$ in the current case) of $F$
   as an ${\cal O}_Y$-module.
 The construction of $F_0$ involves mutations on exceptional pairs on $Z$ ([K-O: Sec.$\hspace{.7ex}$1])
     --
	  beginning with something that is more accessible but does not pull back to a tilting object on Y,
	     e.g.$\hspace{.7ex}{\cal F}_0^{\mbox{\tiny\it\, initial}}
		     = \pi_{{\Bbb P}^2}^\ast({\cal O}_{{\Bbb P}^2})
			      \oplus \pi_{{\Bbb P}^2}^\ast ({\cal O}_{{\Bbb P}^2}(1))
				  \oplus \pi_{{\Bbb P}^2}^\ast ({\cal O}_{{\Bbb P}^2}(2))
				  \oplus {\cal O}_Z(C_1)\oplus \cdots \oplus {\cal O}_Z(C_k)$,
			where $\pi_{{\Bbb P}^2}: Z\rightarrow {\Bbb P}^2 $
			    is the blowing-up of ${\Bbb P}^2$ at $k$ distinct points, $k\le 8$,
			   and $C_1\,,\,\cdots\,,\, C_k\subset Z$ are the exceptional curves
	 --
   which makes the understanding of the apical algebra in this case formidable.
  {\it Assuming that in the end ${\cal F}_0$ remains to have an invertible sheaf ${\cal L}$ as a direct summand}
     and write ${\cal F}_0={\cal L}_0\oplus {\cal F}_0^\prime$.
   Then
      ${\cal F}:=  \pi_Z^\ast {\cal F}_0
	      = \pi_Z^\ast {\cal L}\oplus \pi_Z^\ast {\cal F}_0 =: {\cal L}\oplus {\cal F}^\prime$\,;\;
     $\End_Y({\cal F})
	     = \End_{\widetilde{Y}}({\cal F}\otimes_{{\cal O}_{\widetilde{Y}}}\!{\cal L}^{-1})$\,;\;
     and we may instead take
	   ${\cal F}\otimes_{{\cal O}_{\widetilde{Y}}}\!{\cal L}^{-1}
	       \simeq  {\cal O}_{\widetilde{Y}}
		                       \oplus  ({\cal F}^\prime \otimes_{{\cal O}_{\widetilde{Y}}}\!{\cal L}^{-1})$
	as the fundamental module on $\widetilde{Y}$.
 Since ${\pi_Y}_{\!\ast} {\cal O}_{\widetilde{Y}}={\cal O}_Y$ in the current case,
   the fundamental module as a ${\cal O}_Y$-module would then have ${\cal O}_Y$ as a direct summand.
   
 \noindent\hspace{40.8em}$\square$ 		
}\end{remark}

\medskip

\begin{motif-question-project} {\bf [apical algebra]}\; {\rm
 Since every commutative ${\Bbb C}$-algebra has a ${\Bbb C}$-algebra-quotient ${\Bbb C}$,
    existence of a ${\Bbb C}$-point on $Y^\nc$ over $y_\ast\in Y$  is a necessary condition for the apical algebra
	to have a nonzero commutative ${\Bbb C}$-algebra as a quotient.
 Thus:  {\it Does a ${\Bbb C}$-point always exist for the apical algebra of a noncommutative crepant resolution $Y^\nc$
   of an isolated Gorenstein singularity $Y$?}
   {\it If yes, can one count them?}
  {\it If yes, are  these ${\Bbb C}$-points related to the irreducible components
              of the exception divisor of a commutative crepant resolution of $Y$?}
  {\it If not true in general, can one give some general criterion for the existence of a ${\Bbb C}$-point for the apical algebra?}
  In particular, study the problem for ADE surface singularities, as special cases of Example~2.1.5 above,
      for the algebra associated to a dimer model in Example~2.1.6 above.  and
	  for the endomorphism algebra associated to the cone over a del Pezzo surface as reviewed in Remark~2.1.7 above.
}\end{motif-question-project}

\bigskip

\subsection{Seed systems over a Gorenstein singular Calabi-Yau space}

The local models of  noncommutative crepant resolutions can be selected and augmented to seed systems over a
  Calabi-Yau space with isolated Gorenstein singularities.

\bigskip

\begin{flushleft}
{\bf General procedure to construct a seed system}
\end{flushleft}
Let $Y$ be a a Gorenstein singular Calabi-Yau space.
A  general procedure to construct a seed system $H_{\frak U}$ for $Y$ goes as follows:
(The Calabi-Yau condition is not used in this section; but we'll stay with it for our goal beyond the current notes.)

\bigskip

\noindent
(1) {\it $C^\infty$-chart of a local model as a mildly singular manifold}\hspace{2em}
Let  $Y_0^\prime$ be a local model of an isolated Gorenstein singularity
  either as an algebraic variety over ${\Bbb C}$ or as an analytic space, or as a formal scheme over ${\Bbb C}$,
  or a formal analytic germ.
In the first two representations, there is already a built-in embedding
  $\iota^\prime : Y^\prime_0\hookrightarrow {\Bbb A}^{k^\prime}_{\Bbb C}$
    or ${\Bbb C}^{k^\prime}$ (as a complex manifold), for some $k^\prime\in{\Bbb Z}_{>0}$,
  either directly from the presentation of $Y^\prime_0$ or from the GAGA principle of Serre.
In the remaining two formal scheme/analytic-space representations, one simply truncate it to a formal neighborhood
  of high enough order from which one can also obtain an embedding
    $\iota^\prime : Y^\prime_{0}\hookrightarrow {\Bbb A}^{k^\prime}_{\Bbb C}$
	or  ${\Bbb C}^{k^\prime}$ (as a complex manifold) for some $k\in {\Bbb Z}_{>0}$.
This gives a global $C^\infty$-chart of $Y^\prime_0$
       covering in particular the singularity $\mathbf{0}^\prime\in Y^\prime_0$     and, hence,
	   realizes $Y^\prime_0$ as a mildly singular $C^\infty$-manifold with one single chart.
For convenience, we'll always assume that
  $\iota^\prime(\mathbf{0}^\prime)=$ the origin $\mathbf{0}\in {\Bbb R}^{2k^\prime}$,
   the underlying $C^\infty$-manifold/scheme of ${\Bbb A}^{k^\prime}_{\Bbb C}$ and ${\Bbb C}^{k^\prime}$.

Note that the local $C^\infty$-chart $\iota^\prime: Y^\prime_0\hookrightarrow {\Bbb R}^{2k^\prime}$
  depends on the presentation of $Y^\prime_0$. We simply fixed one.
Once $\iota^\prime$ is fixed, a given noncommutative crepant resolution ${\cal A}$ of $Y^\prime_0$,
    as an ${\cal O}_{Y^\prime_0}$-algebra,
 is then realized as the torsion ${\cal O}_{{\Bbb R}^{2k^\prime}}^{\,{\Bbb C}}$-algebra
 $\iota^\prime_\ast {\cal A}$,  supported on $\iota^\prime(Y^\prime_0)$.

\bigskip

\noindent
(2) {\it Regular covering ${\frak U}$ of $Y$}\hspace{2em}
A simplest regular covering of $Y$ can be constructed as follows:
  \begin{itemize}
    \item[(2.{\it a})]   [{\it grafting local model of singularity to $Y$}\,]\hspace{2em}
	Let $Y_\scriptsizesing= \{y_{\ast1}\,,\, \cdots\,,\, y_{\ast l}\}$  and,
	for each $y_{\ast i}\in Y_\scriptsizesing$,
	  let $\alpha_i:U_i\hookrightarrow  {\Bbb R}^{2k_i}$ be a local $C^\infty$-chart at $y_{\ast i}\in Y$ and
	       $\iota^\prime_i: U^\prime_i \hookrightarrow {\Bbb R}^{2k^\prime_i }$ be a local $C^\infty$-chart at
		   $\mathbf{0}^\prime_i$ in a local model in Item (1) for the singularity $y_{\ast i}\in Y$.		   
	By shrinking if necessary,
     one may assume that the closures $\overline{U_1}\,,\, \cdots\,,\, \overline{U_l}$ are mutually disjoint;
       that each $\partial \overline{U_i}$ is an embedded $C^\infty$-submanifold of codimension $1$; and
       that there is a $C^\infty$-map
	     $\beta_i: B^\prime_i \subset {\Bbb R}^{2k^\prime_i}\rightarrow {\Bbb R}^{2k_i}$,
		  with the origin $\mathbf{0}\in {\Bbb R}^{2k_i^\prime}$
		     sent to the origin ${\mathbf{0}}\in {\Bbb R}^{2k_i}$, such that
		 \begin{itemize}
		  \item[$\cdot$]
		   the restriction of $\beta_i$ to $B^\prime_i\cap \iota^\prime_i(U^\prime_i)$ is an embedding
		   into ${\Bbb R}^{2k_i}$,
		
		  \item[$\cdot$]
		   $\alpha_i(\overline{U_i}) \subset \beta_i(B^\prime_i\cap \iota^\prime_i(U^\prime_i))$,
		
		  \item[$\cdot$]
		    with an abuse of notation, $\beta_i$ admits a partially defined $C^\infty$ inverse
		     $$
			    \beta_i^{\;-1} \; :\;  \alpha_i(\overline{U_i} -\{y_{\ast i}\})\,\subset\, {\Bbb R}^{2k_i}\;
			    \longrightarrow\;
				B^\prime_i\cap \iota^\prime_i(U^\prime_i-\{\mathbf{0}^\prime_i\})\,
				     \subset\, {\Bbb R}^{2k^\prime_i}
			 $$
		     (which is then a diffeomorphism to its image).
		 \end{itemize}
	   Here, $B^\prime_i$ is a small open ball at the origin $\mathbf{0}\in {\Bbb R}^{2k^\prime_i}$.
	
	\item[(2.{\it b})]  [{\it adding one more open set to cover $Y$}\,]\hspace{2em}
	 For all $i\in \{1\,,\,\cdots\,,\, l\}$,
	 let $U_i^\flat $ be another open neighborhood of $y_{\ast i}$ in $Y$ such that
	   $\overline{U_i^\flat}  \subset  U_i$    and that
	   $\overline{U_i}-\overline{U_i^\flat}$ is a collar of $\partial U_i$ in $\overline{U_i}$.
	 Let $U_0 := Y-\amalg_{i=1}^l \overline{U_i^\flat}$ and
	   ${\frak U}:= \{U_0\,,\, U_1\,,\, \cdots\,,\, U_l\}$.
	 Then ${\frak U}$ is a regular covering of $Y$ in the sense of Definition~1.4.7.
																												 
    \item[]\hspace{1.2em}Note
	that the regular partition ${\cal P}_{\frak U}$ of  $Y$ associated to ${\frak U}$ thus constructed
	 has only codimension-$1$ $C^\infty$ walls and no higher codimensional strata.
	It is, hence, a simplest kind of poartition of $Y$.  	
  \end{itemize}

\bigskip

\noindent
(3) {\it Construction of a seed system $H_{\frak U}$ over ${\frak U}$ }\hspace{2em}
 The mildly singular $C^\infty$-bundle block $H_i$ over $U_i$, $i\in \{1, \cdots, l\}$, arises from
    the push-pull  $\alpha_i^\ast   {\beta_i}_\ast {\iota^\prime_i}_\ast({\cal F}_i)$
	of the fundamental module of the selected local model for $y_{\ast i}\in Y$
	that admits a nowhere-vanishing global section.
 Together with a selection of a $C^\infty$ ${\Bbb C}$-vector bundle $F_0$ on $U_0$
      that admits a nowhere-vanishing global section,
   this defines a seed system $H_{\frak U}$ on $Y$.

{\sc Figure}~2-2-1.
    %
    %
\begin{figure}[htbp]
  \bigskip
   \centering
   \includegraphics[width=0.80\textwidth]{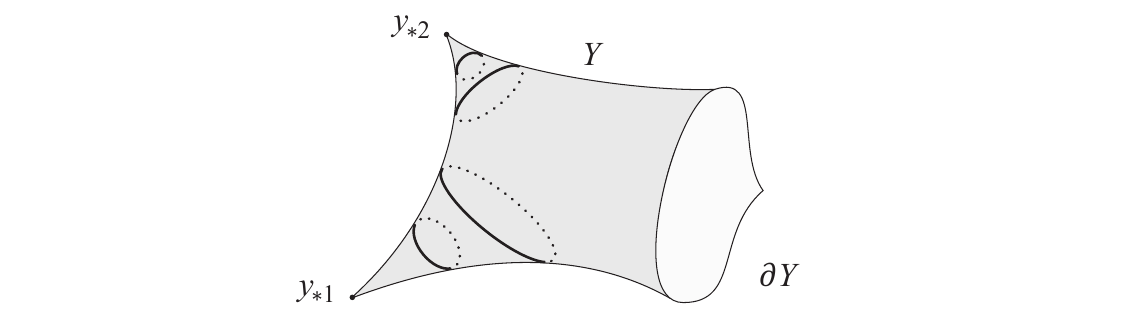}
 
  \bigskip
 \centerline{\parbox{13cm}{\small\baselineskip 12pt
  {\sc Figure}~2-2-1.
   Illustrated is the regular partition ${\cal P}_{\frak U}$ of a manifold $Y$ with corners
        associated to a simplest regular covering ${\frak U}$ of $Y$
      that confines each isolated Gorenstein singularity $y_{\ast i}$, $i=1\,,\,\cdots\,,\, l$.
    Assuming the existence of local noncommutative crepant resolutions of these singularities,
	  a seed system $H_{\frak U}$ over $Y$ can be easily constructed.
  }}
\end{figure}

\bigskip

\begin{example} {\bf [K3 surface with ADE singularities]}\; {\rm
 K3 surfaces are of interest both mathematically and string-theoretically (e.g.\ [Bej], [Hu], and [As1]).
 Singular K3 surfaces (or more generally, singular Calabi-Yau spaces) play even a greater role than smooth ones
  from the perspective of geometric engineering of quantum field theories, (e.g.\ [K-K-V]).
 Given a singular K3 surface $Y$ with ADE singularities, Example~2.1.5 of Sec.~2.1 gives a description of a local model
   of noncommutative crepant resolutions around each singularity.
 After grafting to $Y$,
    this gives rise to a mildly singular $C^\infty$ ${\Bbb C}$-vector-bundle block over a neighborhood
	 of each singularity of $Y$, with a nowhere-vanishing section.
 The above procedure completes the collection of mildly singular bundle blocks to a seed system $H_{\frak U}$ over $Y$,
  which then, following Sec.~1.3, defines a bundle settlement ${\frak S}_{H_{\frak U},\otimes}$ and
  hence a $C^\infty$ noncommutative ringed space $Y^\nc$.
 A systematic study of bundle settlements on singular K3 surfaces remains pending.
}\end{example}
 	
\bigskip

\begin{flushleft}
{\bf Special examples: Simultaneous local noncommutative crepant resolutions}
\end{flushleft}
Once knowing a general procedure exists, a few nonstandard examples from a special one-element seed are given below.
These examples fits closer to the notion of noncommutative resolution of singularities in [B-O: Sec.\ 5] of
  {\sl Alexei Bondal} and {\sl Dmitri Orlov} than those from a general multi-element seed on $Y$..

\bigskip

\begin{example} {\bf [Calabi-Yau conifold]}\; {\rm
 Let
  $Y$ be a Calabi-Yau $3$-fold (either compact or non-compact)
     with ordinary double point singularities $q_1,\,\cdots\,, q_{l_0}$
  (i.e.\ $Y$ is a Calabi-Yau {\it conifold}).
 Note that at this point, $Y$ is treated as a scheme over ${\Bbb C}$
    in the sense of (commutative) Algebraic Geometry in, e.g., [Hart] 	
   (with a K\"{a}hler structure on the underlying holomorphic $3$-manifold with isolated hypersurface singularities).
 For $\Lambda\subset \{1,\cdots, l_0\}$,
  let ${\cal I}_\Lambda \subset {\cal O}_Y$ be a sheaf of ideals on $Y$
      such that
  	    the blowing-up of $Y$ with respect to ${\cal I}_\Lambda$ gives rise to a partial small resolution
               $\pi_\Lambda^+: Y^+_\Lambda\rightarrow Y$ of $Y$ that resolves the singularities $q_i$, $i\in \Lambda$,
  		   while leaving other singularities $q_j$, $j\notin \Lambda$, intact.
 The exceptional locus of $\pi^+_\Lambda$ consists of a ${\Bbb P}^1$ over each $q_i$, $i \in\Lambda$,
  with normal bundle $\simeq {\cal O}_{{\Bbb P}^1}(-1)\oplus {\cal O}_{{\Bbb P}^1}(-1)$ ; and
 it is known
  that ${\cal I}_\Lambda$ is invertible on $Y-\{q_i | i\in \Lambda\}$,
  that $\pi^+_\Lambda$ is crepant, and
  that there is a diagram of morphisms
   $$
    \xymatrix{
       & \widetilde{Y}_\Lambda  \ar[ld] \ar[rd] \\
    Y^+_\Lambda \ar[rd]_-{\pi^+_\Lambda} & & Y^-_\Lambda  \ar[ld]^-{\pi^-_\Lambda}\\
      & Y
    }
   $$
  over $Y$,
   where
     $\widetilde{Y}_\Lambda$ is the blowing-up of $Y$ at $q_i$, $i\in \Lambda$,    and
     $Y^+_\Lambda$ and $Y^-_\Lambda$ are related by Atiyah flops
	  (i.e.\ blowing up along the exceptional ${\Bbb P}^1$
	       and then blowing down the curve class associated to a lifting of this exceptional ${\Bbb P}^1$)
	  over $q_i$, $i\in \Lambda$.

 Let
  $\{1,\cdots, l_0\}=\Lambda_1 \amalg \Lambda_2$ be a decomposition of $\{1,\cdots, l_0\}$
     into a disjoint union   and
  ${\cal I}_{\Lambda_1}$, ${\cal I}_{\Lambda_2}$ be ideal sheaves on $Y$ as above.
 Here, if $\Lambda_i$ is the empty set $\emptyset$,
   then ${\cal I}_{\Lambda_i}= {\cal I}_\emptyset := {\cal O}_Y$ by convention.
 Then
   $$
     {\cal A}_{Y;\,(\Lambda_1,\Lambda_2)}   \;
	   :=\; \Endsheaf_{{\cal O}_Y}({\cal I}_{\Lambda_1}\oplus {\cal I}_{\Lambda_2})
   $$
  is a sheaf of noncommutative ${\cal O}_Y$-algebras on $Y$ with center ${\cal O}_Y$.
 By construction,
  ${\cal I}_{\Lambda_1}\oplus {\cal I}_{\Lambda_2}$ is locally free of rank $2$
     on $Y_\ast := Y_{\smoothscriptsize} := Y-\{q_1,\cdots, q_{l_0}\}$ and, thus,
  the restriction ${\cal A}_{Y;\,(\Lambda_1,\Lambda_2)}|_{Y_\ast}$
     is a sheaf of Azumaya/matrix algebras of rank $2$ on $Y_\ast$.
  
   \smallskip

   \begin{itemize}
    \item[]\parbox[t]{38em}{{\bf Definition 2.2.1.1. [simultaneous local noncommutative crepant resolution of $Y$]}\;
	The noncommutative ringed space
	$Y^\nc_{(\Lambda_1, \Lambda_2)}:=(Y, {\cal A}_{Y;\,(\Lambda_1,\Lambda_2)})$
     is called a {\it simultaneous local noncommutative crepant resolution} of the conifold $Y$. 	
                       }  
   \end{itemize}

   \smallskip
					
 \noindent
 We shall think of $Y^\nc_{(\Lambda_1, \Lambda_2)}$
  as a {\it smooth noncommutative Calabi-Yau $3$-fold} over the singular commutative Calabi-Yau $3$-fold $Y$
  with the same underlying topology as that on $Y$, except over the ordinary double points $q_1, \cdots, q_{l_0}$,
     should the global issue of Calabi-Yau Condition is resolved (cf.\ Motif/Question/Project~2.2.6) and
 $Y^\nc_{(\Lambda_1, \Lambda_2)}$ does satisfy the condition.
 It follows from Example~2.1.4 that
   there are ${\Bbb C}$-points in $Y^\nc_{(\Lambda_1, \Lambda_2)}$
       over the conifold points $q_1\,,\,\cdots\,,\, q_{l_0}$ of $Y$.

\noindent\hspace{40.8em}$\square$
}\end{example}

\medskip

\begin{example} {\bf [singular Kummer surface]}\; {\rm  (E.g.\ [B-P-VdV: Sec.\ 5.16].)
 This is a special class of singular K3 surfaces, obtained from taking the ${\Bbb Z}/2$ quotient of
    a complex $2$-torus $T$ (as a complex abelian Lie group) under the involution $g \mapsto -g$.
 The quotient $Y:= T/({\Bbb Z}/2)$  has $16$ $A_1$ singularities (and smooth elsewhere).
 Denote the quotient map of (analytic spaces) by $\pi: T\rightarrow Y$.
 Then $Y^\nc  :=  (Y, \Endsheaf_{{\cal O}_Y}(\pi_\ast {\cal O}_T) )$
  gives a simultaneous local noncommutative crepant resolution of $Y$.
 Since $\pi_\ast {\cal O}_T$ is a torsion-free ${\cal O}_Y$-module of rank $2$,
   it cannot have ${\cal O}_Y^{\oplus 2}$ as a direct summand.
 It follows from Lemma~2.1.5.1 and Lemma~2.1.5.2 in Example~2.1.5
     that $Y^\nc$ has a ${\Bbb C}$-point over each singularity of $Y$.                				
}\end{example}

\bigskip

\begin{example} {\bf [abelian quotient of special complex $3$-torus]}\;   {\rm
  (Cf.\ [G-S-W: Sec.\ 9.5.2]; [Pr], [Y-Y].)
  For $\tau\in {\Bbb H}:= \{z\in {\Bbb C}\,|\, \Imaginary z >0 \}$,
     let $T_\tau:= {\Bbb C}/({\Bbb Z}+{\Bbb Z}\tau)$ be a complex $1$-torus of modulus $\tau$.
  Let $\tau_0= e^{2\pi\sqrt{-1}/3}=\omega$,
     then ${\Bbb Z}/3$ acts on $T_{\tau_0}$	by $z\mapsto \omega z$, with $3$ fixed points.
  Consider the quotient of the product $\widehat{Y}:= T_{\tau_0}\times T_{\tau_0}\times T_{\tau_0}$
     by the diagonal ${\Bbb Z}/3$-action $(z^1, z^2, z^3)\mapsto (\omega z^1, \omega z^2, \omega z^3)$,
	  $\pi: \widehat{Y}\longrightarrow Y:= \widehat{Y}/({\Bbb Z}/3)$.
  The diagonal action now has $27$ fixed points,
      which gives rise to  $27$ abelian isolated Gorenstein quotient singularities on $Y$.
   Let $Y^\nc := (Y, \Endsheaf_{{\cal O}_Y}({\cal O}_{\widehat{Y}}))$.
   Then the built-in $\pi: Y^\nc\rightarrow Y$ realizes $Y^\nc$
       as a simultaneous local noncommutative crepant resolution of $Y$.
}\end{example}

\medskip

\begin{motif-question-project} {\bf [existence of simultaneous local noncommutative cre-pant resolution]}\; {\rm
 Motivated by Example~2.2.2, Example~2.2.3. and Example~2.2.4:
 {\it Does a normal variety with only isolated Gorenstein singularities of the same local model always
         admit a simultaneous local noncommutative crepant resolution?}
  Of particular interest from heterotic string theory is the question for a K3 surface with exactly two $E_8$-singularities.
  {\it If such a simultaneous local noncommutative crepant resolution cannot be realized
             in the realm of Algebraic or Analytic Geometry, can it be realized in the much flexible $C^\infty$ Algebraic Geometry?}
   (Cf.\ [Bo-K], [FS-J], [Joy2], [Ka], )			
}\end{motif-question-project}

\medskip

\begin{motif-question-project} {\bf [(global) Calabi-Yau condition on $(Y, {\cal O}_Y^\nc)$]}\; {\rm
  To be more specific, let $Y$ be a compact Calabi-Yau manifold with isolated Gorenstein singularities
       whose local noncommutative crepant resolutions in the sense of  Van den Bergh [VdB2] exist.
  Purely algebraically/algebrogemetrically, there is no hope that these local noncommutative resolutions
     can glue to a global object over $Y$ unless each local model is enhanced via Morita-type equivalence so that at least
   all the fundamental module involved are of the same rank as ${\cal O}_Y$-modules.
  Part of the motivation of the bundle settlements via a seed system over $Y$ is exactly to deal with this issue of
    {\it gluing local ${\cal O}_Y$-modules of different ranks up to Morita equivalence}.
 However, one cannot deny that something distorting, not exactly Morita-like, happens along the wall of the bundle settlement
    in the  construction of  the noncommutative structure sheaf ${\cal O}_Y^\nc$ over $Y$ from a bundle settlement.
 Thus, one definitely needs a global condition to justify that $Y^\nc ;= (Y, {\cal O}_Y^\nc)$
    is (resp.\  is not good enough to be) a {\it `Calabi-Yau noncommutative ringed space'}
	(resp.\ and a re-choice of the the seed system, while keeping the local noncommutative crepant resolutions intact, is necessary).
 With the existing notion of Calabi-Yau category (cf. [Ko], [Ke2]) and
  the development of Microlocal Geometry (cf.\ [Na1], [Na2], [N-Z]), one has a guide:
 Construct or give a condition for the existence of a triangulated category of (a special class of ) modules over
   the structure sheaf ${\cal O}_Y^\nc$ associated to a bundle settlement
  (so that $\Hom$ in this category are all finite-dimensional) that allows a Serre functor.
 The condition that $Y^\nc$ is Calabi-Yau is then that the construction works and produces a Calabi-Yau category.
  Well, easier said than done, details remain pending.
}\end{motif-question-project}

\medskip

\begin{motif-question-project} {\bf [compactification on $Y^\nc$]}\;  {\rm
 Back to the stringy origin of the problem   and
 recall the lower-dimensional effective field theory aspect of Mirror Symmetry:
 \begin{itemize}
   \item[{\bf Q.}]
   {\it Can any of the Type I, Type IIA, Type IIB, $E_8\times E_8$ heterotic, $\SO(32)$ heterotic string theories, M-theory, F-theory
            be compactified on some $Y^\nc$?\\[.8ex]			
	$\mbox{\hspace{-1.4em}}\tinybullet\hspace{.6em}$	
          How may this be related to the global Calabi-Yau condition on $Y^\nc$?\\[.8ex]
	$\mbox{\hspace{-1.4em}}\tinybullet\hspace{.6em}$     		
	      Does such a compactification, if allowed, lead to the lower-dimensional effect field theory aspect of Noncommutative Mirror Symmetry?}
  \end{itemize}		
}\end{motif-question-project}

\bigskip
\bigskip
  
\section{Dynamical D-branes on a $C^\infty$ Azumaya-type noncommutative ringed space}

We now explore dynamical D-branes moving on a $C^\infty$ Azumaya-type noncommutative ringed space
  as constructed in Sec.~1 and Sec.~2.
This is an uncharted territory but one can at least try to lay down a few basic ingredients
  -- as guided by Classical Mechanics, Polyakov strings, and the study of harmonic maps -- to begin with.
In Sec.\ 3.1 we explain how dynamical D-branes in string theory are described in the current setting.
In Sec.\ 3.2 some basics of noncommutative Riemannian geometry are reviewed,
 based on which an action functional for dynamical D-branes in the current setting is then constructed, in parallel to the Polyakov string.

\bigskip
  
\subsection{($C^\infty$) maps from a D-brane world-volume to a $C^\infty$ Azumaya-type noncommutative ringed space}

Let
 \begin{itemize}
  \item[$\cdot$]
   $X$ be a $C^\infty$-manifold,
      ${\cal O}_X$ its structure sheaf of $C^\infty$-functions;
   $E$ a ${\Bbb C}$-vector bundle of rank $r$ over $X$,
     ${\cal E}$ the ${\cal O}_X^{\,\Bbb C}$-module of $C^\infty$-sections of $E$;
   $X^{\!\Azscriptsize} := (X, {\cal O}_X^\Azscriptsize := \Endsheaf_{{}_X^{\,\Bbb C}}({\cal E}))$
    the $C^\infty$ Azumaya manifold associated to the pair $(X, E)$ (or equivalently $(X,{\cal E})$);
     
  \item[$\cdot$]
   $Y$ be a mildly singular $C^\infty$-manifold, ${\cal O}_Y$  its structure sheaf of $C^\infty$-functions;
   ${\frak S}:={\frak S}_{\cal P}$ a mildly singular $C^\infty$ ${\Bbb C}$-vector-bundle settlement over $Y$
              with a regular partition ${\cal P}$,
       ${\cal S}_{\frak S}$ the ${\cal O}_Y^{\,\Bbb C}$-module of $C^\infty$-sections of ${\frak S}$;
  $Y^\nc:=(Y, {\cal O}_Y^\nc :=\Endsheaf_{{\cal O}_Y^{\,\Bbb C}}({\cal S}_{\frak S}))$
    the $C^\infty$ Azumaya-type noncommutative ringed space associated to the pair $(Y, {\frak S})$.
 \end{itemize}
For terminology, ${\cal E}$ (resp.$\hspace{.7ex}{\cal S}_{\frak S}$) is called the {\it fundamental module}
 of $X^{\!\Azscriptsize}$ (resp.$\hspace{.7ex}Y^\nc$).

\bigskip

\begin{flushleft}
{\bf Dynamical D-branes on $Y^\nc$ presented as as maps $\varphi: X^{\!A\!z}\rightarrow Y^\nc$}
\end{flushleft}
Translating
   \begin{itemize}
    \item[$\cdot$]
     the {\it matrix-enhancement behavior}, when stacked, of the scalar fields on a D-brane world-volume
          that describe its space-time coordinates   and
	
    \item[$\cdot$]	
     the {\it Higgsing-unHiggsing behavior} of D-branes under deformations
    \end{itemize}
  from String Theory to Algebraic Geometry
 leads one to take an Azumaya manifold with a fundamental module $(X^{\!A\!z}, {\cal E})$
       as the world-volume for stacked D-branes    and
   describe a dynamical D-brane in a target-space(-time) as a map from $X^{\!A\!z}$, defined as follows:
(See [Liu: Sec.\ 1], [L-Y1: Sec.\ 1 \& Sec.\ 2] (D(1)), [L-Y4: Sec.\ 1] (D(11.1)) for detailed explanations,
   beginning with [Pol4: vol.\ I, Sec. 8.7, p.272] of Polchinski;
   see also [G-Sh], [H-W] and [Wi2].)
   
\bigskip

\begin{definition}  {\bf [map/morphism from $X^{\!\Azscriptsize}$ to $Y^\nc$]}\; {\rm
 A  {\it map} (synonymously {\it morphism})
     $$
	   \varphi\::\: X^{\!\Azscriptsize}\;  \longrightarrow\; Y^\nc
	 $$
	is defined contravariantly as an equivalence class of gluing systems of  $C^\infty$ ${\Bbb C}$-algebra-homomorphi-sms,
	in notation
	 $$
	  \varphi^\sharp\::\: {\cal O}_Y^\nc \; \longrightarrow\; {\cal O}_X^\Azscriptsize\,,
	 $$
 ({\it without} a specification of a continuous map $X\rightarrow Y$ as topological spaces).	
{\sc Figure} 3-1-1.
      %
	  %
  \begin{figure}[htbp]
    \bigskip
     \centering
     \includegraphics[width=0.80\textwidth]{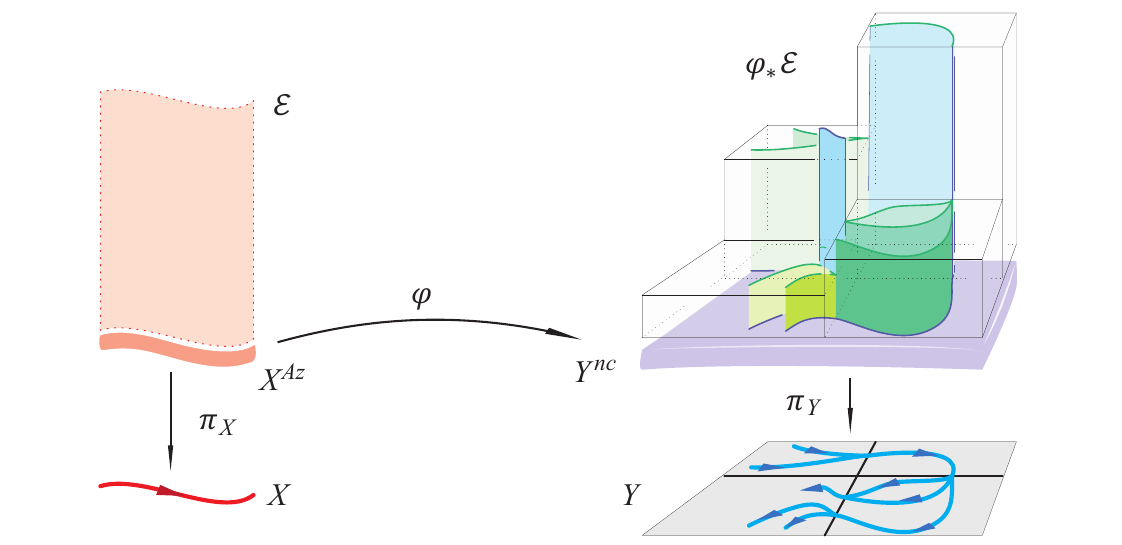}
 
    \bigskip
   \centerline{\parbox{13cm}{\small\baselineskip 12pt
    {\sc Figure}~3-1-1.
      A map/morphism $\varphi$ from $X^{\!A\!z}$ to $Y^{nc}$ is indicated.
	  Here, both the rank of the Chan-Paton sheaf ${\cal E}$ on $X$ and the rank of bundle blocks
	    on the regularly partitioned target-space $Y$
	     are presented via the height in the illustration.
	  The push-forward $\varphi_\ast{\cal E}$ of ${\cal E}$ to $Y^\nc$ and its scheme-theoretic support on $Y$
	    are indicated by {\sc Blue-Green}.
	  To make clearer the image $C^\infty$-scheme of  $X^{\!A\!z}$ in $Y$ under $ \pi_Y \circ \varphi$,
	     an orientation on $X$ is added and the induced orientation on the image is indicated. \\
      $\mbox{\hspace{1.2em}}$Note
	  that $\varphi$ is defined contravariantly via an equivalence class
	    $\varphi^\sharp: {\cal O}_Y^{nc}\rightarrow {\cal O}_X^{Az}$ of gluing systems of
         ring-homomorphisms; thus in general $\varphi$ does not induce any map $X\rightarrow Y$.
	 It is this feature that allows $\varphi$ to correctly depict the Higgsing-unHiggsing behavior of D-branes under deformations
	    as described in String Theory.
   }}
  \end{figure}	
}\end{definition}

\medskip

\begin{remark} {$[$morphism between gluing systems of rings$]$}\; {\rm
 In Definition~3.1.1,
 one should think of ${\cal O}_X^{A\!z}$	(resp.\ ${\cal O}_Y^\nc$)
   as defining an atlas of $C^\infty$ ${\Bbb C}$-algebras
      (i.e~${\Bbb C}$-algebras with center the complexification of a $C^\infty$-ring) via central localizations on $X$ (resp.~$Y$).
 A morphism $X^{\!\Azscriptsize}\rightarrow Y^\nc$ is defined contravariantly as an equivalence class of  gluing systems
  of ${\Bbb C}$-algebra-homomorphisms whose restriction to the $C^\infty$-rings in the centers  are  $C^\infty$-ring-homomorphism
      to its image.
 See [L-Y1: Sec.~1.2] (D(1)) for details and
         [L-L-S-Y] (D(2)), [L-Y4] D(11.1) for further explanations.
 Cf.~{\sc Figure}~3-1-1.		
}\end{remark}
 
\bigskip

Since $X$ and $Y$ are determined by $C^\infty(X)$ and $C^\infty(Y)$ respectively, as $C^\infty$-schemes,   and
  ${\cal O}_X^\Azscriptsize$ and ${\cal O}_Y^\nc$ are locally finite extension of
  ${\cal O}_X^{\,\Bbb C}$ and ${\cal O}_Y^{\,\Bbb C}$ respectively, as sheaves of ${\Bbb C}$-algebras,
 the data $\varphi^\sharp:{\cal O}_Y^\nc\rightarrow {\cal O}_X^\Azscriptsize$ in Definition~3.1.1,
  on one hand, induces and, on the other hand, is determined by a $C^\infty$ ${\Bbb C}$-algebra-homomorphism
   $$
    \varphi^\sharp\::\:
	  C^\infty(Y^\nc):=\varGamma({\cal O}_Y^\nc)\; \longrightarrow\;
	  C^\infty(X^{\!\Azscriptsize}):= \varGamma({\cal O}_X^\Azscriptsize)\,,
   $$
   where
     $\varGamma(\tinybullet)$ is the $C^\infty$ global-section functor and
     the induced ${\Bbb R}$-algebra-homomorphism
	 $$
	   \underline{\varphi}^\sharp:
	     C^\infty(Y)\; \longrightarrow\;  \varphi^\sharp (C^\infty(Y))\subset C^\infty(X^{\!\Azscriptsize})
	 $$
	 is a $C^\infty$-ring-homomorphism.	
$\underline{\varphi}^\sharp$ defines the composition $\underline{\varphi}$ of
  $X^{\!A\!z}\stackrel{\varphi}{\longrightarrow} Y^\nc
     \stackrel{\pi_Y}{\longrightaarrow}Y$.
	
\bigskip

\begin{definition} {\bf [surrogate $X_\varphi$ of $X^{\!A\!a}$ associated to $\varphi$\,; morphism $f_\varphi$]}\;
{\rm
 (Cf.~[L-Y1: Definition~1.1.1, Definition/Example 1.1.2] (D(1)),
          [L-L-S-Y: Definition 2.1.1] (D(2)),
		  [L-Y4: {\sc Figure} 1-2, Definition~5.1.4, Lemma/Definition 5.3.1.7] (D(11.1)).)\hspace{1em}
  (1) Let
	\begin{itemize}
	 \item[$\cdot$]
      $A_{\underline{\varphi}}\subset C^\infty(X^{\!A\!z})$ be the $C^\infty$-ring generated by
	     the subrings $C^\infty(X)$ and $\underline{\varphi}^\sharp(C^\infty(Y))$;
	
     \item[$\cdot$]	
	  $A_\varphi\subset C^\infty(X^{\,A\!z})$ be the subring generated by
	    $A_{\underline{\varphi}}^{\Bbb C}$ and $\varphi^\sharp(C^\infty(Y^\nc))$;
	
	 \item[$\cdot$]
	  $X_{\underline{\varphi}}$ be the $C^\infty$-scheme associated to the $C^\infty$-ring
	    $A_{\underline{\varphi}}$.	
   \end{itemize}
  (Note that $X_{\underline{\varphi}}$ is the surrogate of $X^{\!A\!z}$ associated to
       $\underline{\varphi}: X^{\!A\!z}\rightarrow Y$, as studies in [L-Y4] (D(11.1)).)
  Denote the structure sheaf of $X_{\underline{\varphi}}$ by ${\cal O}_{X_{\underline{\varphi}}}$
     (or simply ${\cal A}_{\underline{\varphi}}$).
  Then the $A_{\underline{\varphi}}$-algebra $A_\varphi$  determines
    an ${\cal O}_{X_{\underline{\varphi}}}$-algebra ${\cal A}_\varphi$.
 The {\it surrogate} of $X^{\!A\!z}$ {\it associated to} $\varphi$ is defined to be
    the ringed space
	  $$
	     X_\varphi\; :=\;  (X_{\underline{\varphi}}, {\cal A}_{\varphi})\,.
	  $$
  
  (2)
    By construction there are a sequence of dominant morphisms
    $X^{\!\Azscriptsize}\rightarrow X_{\varphi}\rightarrow  X_{\underline{\varphi}} \rightarrow X$
	associated to the built-in inclusions
	${\cal O}_X\hookrightarrow  {\cal A}_{\underline{\varphi}}
	      \hookrightarrow {\cal A}_\varphi  \hookrightarrow {\cal O}_X^{A\!z}$
	 of ${\cal O}_X$-algebras.
  Through these inclusions,
    the fundamental (left) ${\cal O}_X^{A\!z}$-module ${\cal E}$ is naturally
       a (left) ${\cal A}_\varphi$-module and an ${\cal A}_{\underline{\varphi}}$-module,
    both still denoted by ${\cal E}$
	(or  $_{{\cal A}_{\varphi}}{\cal E}$ and $_{{\cal A}_{\underline{\varphi}}}{\cal E}$
	          respectively when need to distinguish them).
			
  (3) The built-in ring-homomorphism
                 $f_{\varphi}^\sharp := \varphi^\sharp: C^\infty(Y^\nc)\rightarrow A_{\varphi}$ and
			  $C^\infty$-ring-homomorphism
              	 $f_{\underline{\varphi}}:= \varphi^\sharp|_{C^\infty(Y)}: C^\infty(Y)
				                                                         \rightarrow A_{\underline{\varphi}}$
			  define natural morphisms
			  $$
			     f_\varphi\;:\; X_\varphi\; \longrightarrow\; Y^\nc
				   \hspace{2em}\mbox{and}\hspace{2em}
				 f_{\underline{\varphi}}\;:\; X_{\underline{\varphi}}\; \longrightarrow\; Y
			  $$
			  respectively.
  They fit into the the following commutative diagram:			
     $$
	   \xymatrix{
	    {\cal E}\ar@{.}[dr] \ar@{.}[ddr]  \ar@{.}@/_.1ex/[dddr]    \ar@{.}@/_/[ddddr]  \\
	           & X^{\!A\!z}\ar@{->>}[d] \ar[drrrr]^-\varphi \ar'[drr][ddrrrr]_(.3){\underline{\varphi}}\\
		       & X_{\varphi} \ar@{->>}[d] \ar[rrrr]_(.3){f_\varphi}         &&&& Y^\nc\ar@{->>}[d]^-{\pi_Y} \\
		       & X_{\underline{\varphi}}\ar@{->>}[d] \ar[rrrr]^(.3){f_{\underline{\varphi}}}      &&&& Y  \\
		       & X      &&&& \hspace{2em}\;\;.
	   }
	 $$
}\end{definition}

\bigskip

\begin{flushleft}
{\bf The graph of $\varphi:X^{\!A\!z}\rightarrow Y^\nc$}
\end{flushleft}
(Cf.~[L-L-S-Y: Sec.~2.2] (D(2)),
           [L-Y4: Sec.~1, Sec.~5] (D(11.1)).)		
Let $X\stackrel{pr_1}{\longlongleftarrow}X\times Y \stackrel{pr_2}{\longlongrightarrow} Y$
 be the projection maps.
Then
 the partition ${\cal P}$  on $Y$ pulls back via $\pr_2$ to a product partition
   $\pr_2^\ast{\cal P}= \{X \times Y_i\,|\, Y_i\in {\cal P}\}$ on $X\times Y$    and
 the $C^\infty$ ${\Bbb C}$-vector-bundle settlement ${\frak S}$ on $_{\cal P}Y$ pulls back to
   a $C^\infty$ ${\Bbb C}$-vector-bundle settlement $\pr_2^\ast{\frak S}$ on $X\times Y$,
 whose associated ${\cal O}_{X\times Y}^{\,\Bbb C}$-module is canonically identical to
   $\pr_2^\ast {\cal S}_{\frak S}$.
The associated noncommutative structure sheaf ${\cal O}_{X\times Y}^\nc$ associated to
  $\pr_2^\ast{\frak S}$  is thus canonically identical to the ${\cal O}_{X\times Y}^{\,\Bbb C}$-algebra
  $\pr_2^\ast {\cal O}_Y^\nc$.
This defines the produce $C^\infty$ Azumaya-type noncommutative ringed space $X\times Y^\nc$
 ($:= (X\times Y, {\cal O}_{X\times Y}^\nc)$).

\bigskip

\begin{lemma} {\bf [unique lifting of $\varphi$ to $X\times Y^\nc$]}\;
 The map $\varphi: X^{\!A\!z}\rightarrow Y^\nc$
    lifts uniquely to a map $\widetilde{\varphi}: X^{\!A\!z}\rightarrow X\times Y^\nc$
  that makes the following diagram commute:
   $$
     \xymatrix{
	   X^{\!A\!z}\ar[dr]^(.6){\widetilde{\varphi}} \ar@/^/[drrr]^-\varphi \ar@/_/ @{->>}[ddr]_-{\pi_X} \\	
	     &  X\times Y^\nc \ar@{->>}[rr]_-{pr_2} \ar@{->>}[d]^-{pr_1}   && Y^\nc \\
		 & X                            && \hspace{2em}\,.
	 }
   $$
\end{lemma}
  
\medskip

\noindent{\it Proof}.
  From [L-Y4: Sec.$\hspace{.7ex}$5] (D(11.1)),
  the ring-homomorphism $\underline{\varphi}^\sharp: C^\infty(Y)\rightarrow C^\infty(X^{\!A\!z})$
    pushes out uniquely to a ring-homomorphism
	$$
	   \widetilde{\underline{\varphi}}\; :\;  C^\infty(X\times Y)\;\longrightarrow\; C^\infty(X^{\!A\!z})\,
	$$
   that is a $C^\infty$-ring-homomorphism to its image  and makes the following diagram commute:
  $$
    \xymatrix{
	  C^\infty(X^{\!A\!z}) \\
	     & C^{\infty}(X\times Y) \ar[ul]_(.3){\widetilde{\underline{\varphi}}^\sharp}
		       && \:\: C^\infty(Y)\ar@/_/[ulll]_-{\underline{\varphi}^\sharp} \ar@{^{(}->}[ll]^-{pr_2^\sharp} \\
		 & \rule{0ex}{1.2em}C^{\infty}(X)\ar@/^/@{_{(}->}[uul]^-{\iota=\pi_X^\sharp}
		                                   \ar@{_{(}->}[u]_-{pr_1^\sharp}          && \hspace{3em}\:\;\;.
	}
  $$
 Since ${\cal O}_Y^\nc$ is a locally finite-dimensional algebraic extension of ${\cal O}_Y$ and
	  ${\cal O}_{X\times Y}^\nc= \pr_2^\ast {\cal O}_Y^\nc$, 	
   the above diagram pushes out by $\pi_Y^\sharp: C^\infty(Y)\rightarrow C^\infty(Y^\nc)$ to a unique commutative diagram
   $$
    \xymatrix{
	  C^\infty(X^{\!A\!z}) \\
	     & C^{\infty}(X\times Y^\nc) \ar[ul]_(.3){\widetilde{\varphi}^\sharp}
		       && \:\: C^\infty(Y^\nc)\ar@/_/[ulll]_-{\varphi^\sharp} \ar@{^{(}->}[ll]^-{pr_2^\sharp} \\
		 & \rule{0ex}{1.2em}C^{\infty}(X)\ar@/^/@{_{(}->}[uul]^-{\iota=\pi_X^\sharp}
		                                   \ar@{_{(}->}[u]_-{pr_1^\sharp}          && \hspace{3em}\:\;\;\;\;\;,
	}
  $$
  with $\widetilde{\varphi}^\sharp \circ \pi_Y^\sharp = \widetilde{\underline{\varphi}}^\sharp$,
  that completes the given diagram
   $$
    \xymatrix{
	  C^\infty(X^{\!A\!z}) \\
	     & C^{\infty}(X\times Y^\nc)
		       && \:\: C^\infty(Y^\nc)\ar@/_/[ulll]_-{\varphi^\sharp} \ar@{^{(}->}[ll]^-{pr_2^\sharp} \\
		 & \rule{0ex}{1.2em}C^{\infty}(X)\ar@/^/@{_{(}->}[uul]^-{\iota=\pi_X^\sharp}
		                                   \ar@{_{(}->}[u]_-{pr_1^\sharp}          && \hspace{3em}\:\;\;\;\;\;.
	}
  $$
  The lemma now follows.
  
\noindent\hspace{40.8em}$\square$
  
\bigskip
	
It follows that
 $$
   \begin{array}{rclclcl}
     C^\infty(X^{\!A\!z})    & \supset
	       & \widetilde{\varphi}^\sharp(C^\infty(X\times Y^\nc))
		        & =   & A_{\varphi}   &  =   & C^\infty(X_{\varphi})\\[1.2ex]
     & \supset
        &	C^\infty(X_{\underline{\varphi}}) &  =   & A_{\underline{\varphi}}
	          & =  &  \widetilde{\underline{\varphi}}^\sharp(C^\infty(X\times Y)) \\[1.2ex]
	& \supset   &  C^\infty(X) \,.
   \end{array}			
 $$
Thus,
  the surrogate $X_{\widetilde{\varphi}}$ of $X^{\!A\!z}$ associated to $\widetilde{\varphi}$
     is identical to $X_{\varphi}$ and
  the surrogate $X_{\widetilde{\underline{\varphi}}}$ of $X^{\!A\!z}$ associated to
     $\widetilde{\underline{\varphi}}$ is identical to $X_{\underline{\varphi}}$\,;    and
 one has the commutative diagram
   $$
	   \xymatrix{
	    {\cal E}\ar@{.}[dr] \ar@{.}[ddr]  \ar@{.}@/_.1ex/[dddr]    \ar@{.}@/_/[ddddr]  \\
	           & X^{\!A\!z}\ar@{->>}[d] \ar[drrrr]^-{\widetilde{\varphi}}
			                            \ar'[drr][ddrrrr]_(.3){\widetilde{\underline{\varphi}}}\\
		       & X_{\varphi} \ar@{->>}[d] \ar@{^{(}->}[rrrr]_(.3){\widetilde{f_\varphi}}
			          &&&& X\times Y^\nc\ar@{->>}[d]^-{\pi_{X\times Y}} \\
		       & X_{\underline{\varphi}}\ar@{->>}[d]
			                                                          \ar@{^{(}->}[rrrr]^(.3){\widetilde{f_{\underline{\varphi}}}}
                 			   &&&& X\times Y \ar[dllll]^-{pr_1}  \\
		       & X      &&&& \hspace{4em}\;\;\;\;.
	   }
	 $$

\bigskip

\begin{definition} {\bf [graph of map $\varphi$\,]}\; {\rm
  Recall the defining fundamental (left) module ${\cal E}$ on $X^{\!A\!z}$.
  The {\it graph} of $\varphi: X^{\!A\!z}\rightarrow Y^\nc$ is by definition
    the pushforward (left) ${\cal O}_{X\times Y}^\nc$-module
    $$
	   \widetilde{\cal E}_\varphi\;  :=\;  \widetilde{\varphi}_{\ast\,} {\cal E}
	      =  {\widetilde{f_\varphi}}_\ast  ({_{{\cal A}_\varphi}}{\cal E})
	$$
	on $X\times Y^\nc$.
 From the above discussion,
  its {\it scheme-theoretic support}  $\Supp(\widetilde{\cal E}_\varphi)$,
     a ringed space with the structure sheaf
	  ${\cal O}_{X\times Y}^\nc/
	      \Ker({\cal O}_{X\times Y}^{nc}
		                \rightarrow \Endsheaf_{{\cal O}_{X\times Y}^{nc}}(\widetilde{\cal E}_\varphi))$,
  is isomorphic to $X_{\varphi}$, as ringed spaces, via the embedding $\widetilde{f_\varphi}$.
 By construction, $\pr_{1\ast}\widetilde{\cal E}_\varphi = {\cal E}$,
    where $\pr_1: X\times Y^\nc \rightarrow X$ is the projection map.
}\end{definition}	

\medskip

\begin{remark} $[$characteristic properties of graph $\widetilde{\varphi}_{\ast\,}{\cal E}$ of $\varphi\,]$\;
{\rm
 Call a ringed space $(Z, {\cal O}_Z)$ {\it of dimension zero}
   if its underlying topological space $Z$ is of dimension zero and ${\cal O}_Z$ is simply a finite-dimensional algebra.
 An ${\cal O}_{X\times Y}^\nc$-module  is said to be {\it of dimension zero}
    if its scheme-theoretic support is of dimension zero.
 Then the graph $\widetilde{\cal E}_\varphi$ of $\varphi$ is characterized by two properties:
   \begin{itemize}
    \item[(1)]
	 The relative dimension of $\widetilde{\cal E}_\varphi$ on $(X\times Y^\nc)/X$ is zero.
	
	\item[(2)]
	$\pr_{1\ast}\widetilde{\cal E}_\varphi = {\cal E}$,
	   where $\pr_1: X\times Y^\nc \rightarrow X$ is the projection map.
   \end{itemize}
  Any (left) ${\cal O}_{X\times Y}^\nc$-module $\widetilde{\cal E}^\prime$
    of relative dimension zero on $(X\times Y^\nc)/X$
	that satisfies the condition $\pr_{1\ast}\widetilde{\cal E}^\prime \simeq {\cal E}$
	defines a map $\varphi^\prime: X^{\!A\!z}\rightarrow Y^\nc$ with graph $\widetilde{\cal E}^\prime$
	once an ${\cal O}_X^{\,\Bbb C}$-module-isomorphism
	  $\pr_{1\ast}\widetilde{\cal E}^\prime \stackrel{\sim}{\longrightarrow} {\cal E}$ is specified.	
 Cf.~[L-L-S-Y: Sec.~2.2] (D(2)), [L-Y4: Sec.~5.3] (D(11.1)).
}\end{remark}

%
%
%

\bigskip

We end this subsection with two remarks that are consequences of Morita equivalence.

\bigskip

\begin{remark} $[$D-brane probe, i.e.$\hspace{.7ex}$existence of map from $X^{\!A\!z}\,]$\; {\rm
 Due to that an Azumaya algebra over ${\Bbb C}$ is a simple algebra,
   in general there may be no maps from $X^{\!A\!z}$ to $Y^\nc.$.
 A simplest example to illustrate this is when $Y^\nc$ is an Azumaya manifold, say,
   with the fundamental ${\cal O}_Y^{\,\Bbb C}$-module of rank $l$.
 It follows from the Morita equivalence
       between ${\cal O}_{X\times Y}^{\,\Bbb C}$-{\it Mod} and ${\cal O}_{X\times Y}^\nc$-{\it Mod}
   that, in this case,
    unless the rank $r$ of  the fundamental module ${\cal E}$ of $X^{\!A\!z}$ is an integer multiple of $l$,
	there is no map from $X^{\!A\!z}$ to $Y^\nc$.
 From this perspective, the $Y^\nc$ constructed from a seed system over a singular Calabi-Yau space $Y$ in
    Sec.$\hspace{.7ex}$2.2 could be more ``{\it D-brane accommodating}''.
 When the apical algebra of a local model of  noncommutative crepant resolutions has the field ${\Bbb C}$ as a quotient,
   the $Y^\nc$ associated to a seed system from such local resolutions has a ${\Bbb C}$-point
   over $Y_\scriptsizesing$.
 In such a case, there is always a map from the D-brane world-volume $X^{\!A\!z}$ to $Y^\nc$.
 In other words, such $Y^\nc$ can be probed by D-branes with the Chan-Paton bundle of arbitrary rank.
}\end{remark}

\medskip

\begin{remark} $[$scattering behavior of $\varphi$ when crossing the wall of $_{\cal P}Y]$\; {\rm
 Recall (e.g.~[P-S: Chapter~4]) that a scattering of particles in Particle Physics (for the $4$-dimensional Minkowski space-time)
     phrased in terms of Quantum Field Theory
	is governed by the action functional for the quantum fields associated to the particles and,
   at the event (i.e.$\hspace{.7ex}$space-time point) of scattering,
	the total momentum of the particle system has to be conserved.
 Recall (e.g.~[L-Y4: Sec.~5] (D(11.1)) )
  that even when $X$ is connected,
    the ${\cal O}_X^{\,\Bbb C}$-algebra ${\cal O}_X^{A\!z}$ can have
    an ${\cal O}_X^{\,\Bbb C}$-subalgebra ${\cal A}$ whose center corresponds to a scheme $Z$
	  that is finite over $X$
	 (i.e.$\hspace{.7ex}(Z,{\cal A})$ is a ringed space that is finite over $X$).
  Thus, even though $X^{\!A\!z}$ is connected,
     the surrogate $X_\varphi$ of  $X^{\!A\!z}$ associated to a map $\varphi: X^{\!A\!z}\rightarrow \tinybullet$
	   can be locally (in the $C^\infty$ or analytic topology) disconnected over $X$   and
     $\varphi(X^{\!A\!z})$  in this case will represent a scattering process of
	    a collection of D-branes with Chan-Paton sheaves of lower ranks on the target-space $\tinybullet$.
  Though no ``Brane Field Theory" is in existence yet, some basic conservation laws should still be obeyed at the event
   of brane fusion or splitting on $\tinybullet$.
  
  When the D-brane target-space $\tinybullet$ is taken to be $Y^\nc$ in the current work,
    a new scattering phenomenon occurs when D-branes cross the wall of $_{\cal P}Y$
	  for the bundle settlement ${\frak S}$ that defines $Y^\nc$.
 The additional ``{\it conservation law}" in this case is governed by
	 \begin{itemize}
       \item[(1)]	
       {\it the Chan-Paton ${\cal O}_X^{\,\Bbb C}$-module ${\cal E}$ on $X$ is of constant rank}, and

	   \item[(2)]
        the {\it Morita equivalence between
          ${\cal O}_{X\times (Y-Y_{sing}- {\,_{\cal P}^{\;\ge 1}Y})}^\nc$-{\it Mod}  and
		  ${\cal O}_{X\times (Y-Y_{sing}- {\,_{\cal P}^{\;\ge 1}Y})}^{\,\Bbb C}$-{\it Mod}\\
	    has to be obeyed,
              in particular when applied to the graph $\widetilde{\cal E}_\phi$ of $\varphi$ on $X\times Y^\nc$}.
     \end{itemize}
 Some {\it scattering diagrams} of D$0$-branes at the wall of $_{\cal P}Y$ are illustrated in {\sc Figure}~3-1-2.
      %
      %
 There, the additional conserved quantity is the {\it total length} of the local D$0$-brane system involved in the scattering.
 Similar statement holds for higher dimensional  D-branes scattering at the wall. 	
 
\begin{figure}[htbp]
 \bigskip
  \centering
  \includegraphics[width=0.80\textwidth]{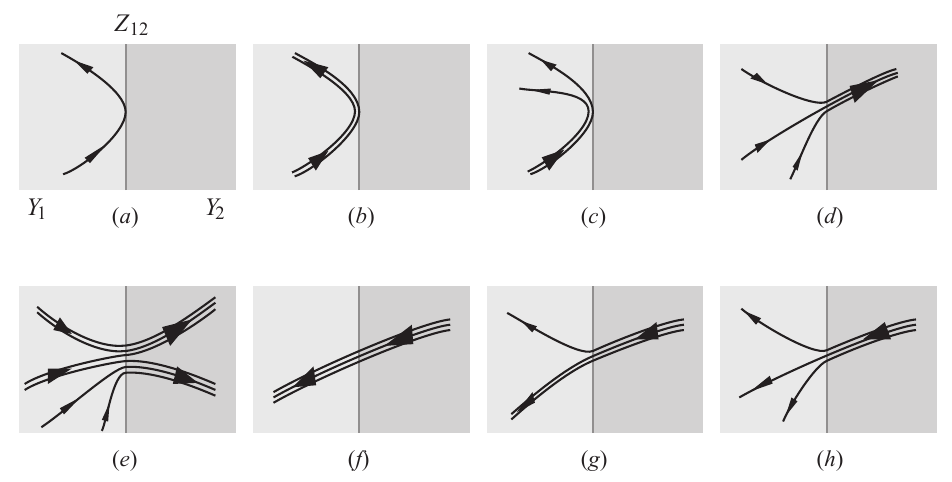}
  
  \bigskip
  \bigskip
 \centerline{\parbox{13cm}{\small\baselineskip 12pt
  {\sc Figure}~3-1-2. {\sc [scattering diagram of D$0$-branes at wall of $_{\cal P}Y$]}\;
   Examples of D-particle scattering at the wall of $_{\cal P}Y$ are illustrated.
   Recall that the {\it length} of a module $M$ over a finite-dimensional ${\Bbb C}$-algebra is simply
      the dimension of $M$ as a ${\Bbb C}$-vector space.
   Such a module, when represented as a $0$-dimensional torsion ${\cal O}_Y^{nc}$-module,
      describes a D-particle on $Y^{nc}$.
   The related world-line is an Azumaya $1$-manifold $X^{\!A\!z}$   and
      the dynamical D-particle in $Y^{nc}$ is then described by a map $\varphi: X^{\!A\!z}\rightarrow Y^{nc}$.
   Here, we focus on a neighborhood of the site at the wall where the scattering occurs.
   The two blocks in ${\cal P}$ involved are denoted by $Y_1$, $Y_2$;
       the wall between them by $Z_{12}$;
	   the $C^\infty$ ${\Bbb C}$-vector-bundle block  over $Y_i$ by $F_i$, $i=1,2$.
   We {\it assume that the rank of $F_1$ is $1$ and the rank of $F_2$ is $3$}.
   The length of each D-particle involved in the scattering is labelled by a numeral and indicated by the number of threads
       of the image world-line in $Y$ under $\pi_Y: Y^{nc}\rightarrow Y$.
   The time-direction of the scattering is indicated by an arrow on the threads.	
   In Examples ({\it a}), ({\it b}), ({\it c}), the incident D-particle system is bounced back by the wall
     since, by Morita equivalence, any $0$-dimensional ${\cal O}_Y^{nc}$-module with support inside $Y_2$
	  must have length an integer multiple of $3$ ($=$ the rank of $F_2$).
	Thus, in such cases, the wall of the bundle settlement  behaves like an infinite potential well,
	   trapping/confining D-particles that do not stack properly at its bottom. \\
    $\mbox{\hspace{1.2em}}$For
	physics-minded readers, note that a D0-brane carries a Ramond-Ramond $1$-form charge.
	Their trajectories may bend when there is a background Ramond-Ramond $1$-form on $Y$,
	which is implicit in all these figures.
  }}
\end{figure}
}\end{remark}

\bigskip

\subsection{Basic ingredients for noncommutative Riemannian geometry on $Y^\nc$}

The main purpose of this subsection is to introduce of the notion of {\it metric tensor}
 on a $C^\infty$ Azumaya-type noncommutative ringed space $Y^\nc$
    as defined in Sec.$\hspace{.7ex}$1   and
	     examplified in Sec.$\hspace{.7ex}$2 via the seed from local noncommutative crepant resolutions.

\bigskip

\begin{flushleft}
{\bf Special case: $Y^\nc$ an Azumaya manifold}
\end{flushleft}
Let $Y^\nc:= (Y, {\cal O}_Y^\nc:=\Endsheaf_{{\cal O}_Y^{\,\Bbb C}}(\cal F))$,
    where $Y$ is a $C^\infty$-manifold    and
             ${\cal F}$ is the sheaf of $C^\infty$-sections of a $C^\infty$ ${\Bbb C}$-vector bundle $F$ on $Y$,
 be a $C^\infty$ Azumaya manifold.
The differential calculus for $Y^\nc$ that fits our need when $Y^\nc$ is taken as a D-brane world-volume
  was given in [L-Y4: Sec.~4]  (D(11.1)), based on earlier works	
  [DV-K-M], [DV-M], [Mas]
    of {\it Michel Dubois-Violette}, {\it Richard~Kerner}, {\it John Madore}, {\it Thierry Masson}
 and  the works 	
  [S\'{e}1], [S\'{e2}]
    of {\it Emmanuel S\'{e}ri\'{e}}.
See also the related works, e.g.,~[B-M2],  [Con], [GB-V-F], [Mad], and [M-M-M].		
The same differential calculus will be used here.
Some necessary facts and immediate consequence for the purpose of Sec.~3.3 are reviewed/given below.
Readers are referred to [L-Y4: Sec.~4] (D(11.1))and related references for details.
 
\bigskip

\begin{definition} {\bf [tangent sheaf \& cotangent sheaf]}\; {\rm
 ([L-Y4: Definition 4.2.2] (D(11.1)).)
  The sheaf of derivations of ${\cal O}_Y^\nc$ is naturally an ${\cal O}_Y^{\,\Bbb C}$-module,
    called the {\it tangent sheaf} of $Y^\nc$, denoted ${\cal T}_\ast Y^\nc$.
  The sheaf of differentials of ${\cal O}_Y^\nc$ is naturally a bi-${\cal O}_Y^\nc$-module,
    called the {\it cotangent sheaf} of $Y^\nc$, denoted ${\cal T}^\ast Y^\nc$.
}\end{definition}
    
\medskip

\begin{lemma} {\bf [splitting of ${\cal T}_{\ast}Y^\nc$]}\;
{\rm ([DV-M: Proposition 3]; [L-Y4: Lemma 4.2.5] (D(11.1)).)}
 Let $\Innsheaf({\cal O}_Y^\nc)$ be the sheaf of inner derivations of ${\cal O}_Y^\nc$.
 Then there is a natural exact sequence of ${\cal O}_Y^{\,\Bbb C}$-modules
  $$
    0\; \longrightarrow\; \Innsheaf({\cal O}_Y^\nc)\;
	     \longrightarrow \; {\cal T}_{\ast}Y^\nc\;
         \longrightarrow\; {\cal T}_{\ast}Y^{\Bbb C}\;  \longrightarrow\; 0\,.
  $$
 Furthermore, any connection
  $\nabla^\prime:{\cal F}\rightarrow {\cal T}_{\ast}Y\otimes_{{\cal O}_Y^{\,\Bbb C}}{\cal F}$
  on ${\cal F}$ induces an embedding (i.e.$\hspace{.7ex}${\sl horizontal lift})
  $$
    \xymatrix{
	 \iota^{\nabla^\prime}\;:\; {\cal T}_{\ast}Y^{\Bbb C}\; \ar@{^{(}->}[r]
	       & \;{\cal T}_{\ast}Y^\nc
	 }
  $$
  as ${\cal O}_Y^{\,\Bbb C}$-modules,
   with  $\xi\mapsto  \nabla^{\prime\prime}_{\xi}$,
  that splits the above short exact sequence.
 Here,
  $\nabla^{\prime\prime}:{\cal O}_Y^\nc \rightarrow
      {\cal T}_{\ast}Y\otimes_{{\cal O}_Y^{\,\Bbb C}}{\cal O}_Y^\nc$
  is the induced connection of $\nabla^\prime$ on ${\cal O}_Y^\nc := \Endsheaf_{{\cal O}_Y^{\,\Bbb C}}{\cal F}$.
\end{lemma}

\bigskip
  
Recall that in our setting for the notion of  {\it differentials of a ring},
  two additional relations are added to the usual linearity over the base ring and the Leibniz rule;
  namely, {\it center commutativity} and  {\it chain rule}:

\bigskip

\begin{definition} {\bf [module of differentials of ring]}\;  {\rm
 ([L-Y4: Definition 4.1.4] (D(11.1)).)
 Let $R$ be a(n associative, unital) $S$-algebra,
  with the center $Z_R$ of $R$ a $C^\infty$-ring that contains $S$ as a $C^\infty$-subring
 Then,
  the {\it module of differentials}, denoted by $\Omega_{R/S}$, is the bi-$R$-module
   generated by the set
   $$
     \{ d(r)\,|\,  r\in R   \}
   $$
   subject to the relations
   $$
    \begin{array}{llll}	
     \mbox{($S$-linearity)}
	  &&& d(a_1r_1+a_2r_2)\;
	            =\; a_1\,d(r_1)\,+\,a_2\,d(r_2)\;
				=\; d(r_1)a_1\,+\, d(r_2)\,a_2 				\,,   \\[.6ex]
	 \mbox{(Leibniz rule)}
	  &&&  d(r_1r_2)\;=\;  d(r_1)\, r_2\,+\, r_1\, d(r_2) \,, \\[.6ex]
	 \mbox{($Z_R$-commutativity)}
	  &&&  d(r_1)\,r_3\;=\; r_3\,d(r_1)
    \end{array}
   $$	
   for all $a_1,\,a_2\in S$, $r_1,\,r_2\in R$, and $r_3\in Z_R$,  and
  $$
   \begin{array}{llll}
    \mbox{(chain rule)}\hspace{4em}
      &&&  d(h(r_1,\,\cdots\,,\,r_s))   \\
	  &&&  \hspace{1em}
                 =\;  \partial_1 h(r_1,\,\cdots\,,\,r_s)\, d(r_1)\;
	                     +\; \cdots\; +\;
			           \partial_s h(r_1,\,\cdots\,,\,r_s)\, d(r_s)	 \hspace{1em}
   \end{array}			
  $$
  for all $h\in C^\infty({\Bbb R}^s)$, $s\in {\Bbb Z}_{\ge 1}$, and $r_1,\,\cdots\,,\,r_s\in Z_R$.
 Denote the image of $d(r)$ under the quotient by $dr$.
 Then, by definition, the built-in map
  $$
    \begin{array}{ccccc}
	 d & :&  R & \longrightarrow  & \Omega_{R/S}\\[.6ex]
	    &&    r  & \longmapsto        & dr
	\end{array}
  $$
  is an $S$-derivation from $R$ to $\Omega_{R/S}$.
}\end{definition}

\bigskip

\noindent
Note that the Leibniz rule and the $Z_R$-commutativity imply that
 $$
   dr_3\,r_1\;=\; r_1\,dr_3  \hspace{4em}
	 \mbox{for all $r_1\in R$ and $r_3\in Z_R$.}
 $$
Note also that, using $dr_1\,r_2= d(r_1r_2)-r_1 dr_2$ for all $r_1,\,r_2\in R$,
 the bi-$R$-module $\Omega_{R/S}$ can be regarded as a left $R$-module generated by
 $\{dr\,|\,r\in R\}$.

\bigskip

It is these additional relations and their consequences as noted above
 that allow us to have a short exact sequence of sheaves of differentials that is dual to the one in Lemma~3.2.2:
  
\bigskip

\begin{lemma} {\bf [splitting of ${\cal T}^{\ast}Y^\nc$]}\;
 Let ${\cal T}^\ast_{Y^{nc}/Y}$
     be the sheaf of relative differentials of ${\cal O}_Y^\nc$ over ${\cal O}_Y^{\,\Bbb C}$.
 Then there is a natural exact sequence of bi-${\cal O}_Y^\nc$-modules
  $$
    0\; \longrightarrow\;
    {\cal O}_Y^\nc \otimes_{{\cal O}_Y}\!{\cal T}^\ast Y   \;	                                          	     
	     \longrightarrow \; {\cal T}^{\ast}Y^\nc\;
         \longrightarrow\; {\cal T}^{\ast}_{Y^{nc}/Y}\;  \longrightarrow\; 0\,.
  $$
 A connection $\nabla^\prime$ on ${\cal F}$ induces a connection on ${\cal O}_Y^\nc$
  that gives a bi-${\cal O}_Y^\nc$-module epimorphism
   $\pi^{\nabla^\prime}: {\cal T}^{\ast}Y^\nc
                                                  \rightarrow {\cal O}_Y^\nc \otimes_{{\cal O}_Y}\!{\cal T}^\ast Y$,
    $\omega \mapsto \omega \circ \iota^{\nabla^\prime}$,
   which splits the short exact sequence.
 Here, $\iota^{\nabla^\prime}:{\cal T}_\ast Y^{\Bbb C}\hookrightarrow {\cal T}_\ast Y^\nc$
   is the horizontal-lift map in Lemma~3.2.2.
\end{lemma}

\medskip

\begin{remark} $[$relative exterior differential\,$]$\; {\rm
  As an immediate corollary of
       [B-M2: Example$\hspace{.7ex}$1.8 - end of Sec.$\hspace{.7ex}$1.1, pp.$\hspace{.7ex}$7-9],
  {\it ${\cal T}^\ast_{Y^{nc}/Y}$ is locally inner}
  in the sense that,
	over a small enough open set of $Y$, there exists a (non-unique) local section
	$\theta\in {\cal T}^\ast_{Y^{nc}/Y}$ such that the relative exterior differential map
	$d_{Y^{nc}/Y}: {\cal O}_Y^\nc\rightarrow {\cal T}^\ast_{Y^{nc}/Y}$
	is locally expressible as $f\mapsto [\theta, f]:= \theta f-f \theta$.
}\end{remark}

\bigskip

The module $\Omega_{R/S}$ of differentials as constructed in Definition~3.2.3
 comes with a universal property as well, which determines $\Omega_{R/S}$ uniquely up to a unique isomorphism:
 
\bigskip

\begin{lemma} {\bf [universal property of $\Omega_{R/S}$]}\;
 {\rm ([L-Y4: Lemma 4.1.5] D(11.1).)}
 Let $\Theta: R\rightarrow M$ be an $S$-derivation from $R$ to $M$.
 Then there exists a unique bi-$R$-module homomorphism $h:\Omega_{R/S}\rightarrow M$
  such that the following diagram commutes
  $$
   \xymatrix{
    & &&& \Omega_{R/S}\ar[dd]^-{h} \\
    & R \ar[urrr]^-{d}   \ar[drrr]_-{\Theta}                \\
    & &&& M  &.
   }
  $$
\end{lemma}

\bigskip

Translated to the $C^\infty$ Azumaya manifold $Y^\nc$, this implies the natural {\it exterior differential map}
  $d: {\cal O}_Y^\nc \rightarrow {\cal T}^\ast Y^\nc$ that fits into the following commutative diagram:
  $$
   \xymatrix{
    {\cal O}_Y^\nc \otimes_{{\cal O}_Y}\!{\cal T}^\ast_Y\hspace{1ex}  \ar  @{^{(}->} [rr]
	       &&{\cal T}^{\ast}Y^\nc  \\
	 \hspace{1em}{\cal O}_Y\hspace{1.1em}\ar  @{^{(}->} [rr]^-\iota  \ar[u]^-{\iota\,\otimes\,d}
	       &&   \;{\cal O}_Y^\nc\;;  \ar[u]_-d
    }
  $$
  and a canonical ${\cal O}_Y^{\,\Bbb C}$-module isomorphism:
  
\bigskip

\begin{lemma} {\bf [derivation vs.$\hspace{.7ex}$differential; ${\cal O}_Y^{\,\Bbb C}$-bilinear pairing]}\;
 {\rm ([L-Y4: Lemma 4.1.6 \& Lemma 4.1.7] (D(11.1)).)}
 There is a canonical isomorphism
  $$
    \imath\::\;
   {\cal T}_\ast Y^\nc\; \stackrel{\sim}{\longrightarrow}\;
      \Homsheaf_{\,\mbox{\scriptsize bi-}{\cal O}_Y^{nc}}({\cal T}^\ast Y^\nc, {\cal O}_Y^\nc)
  $$
  as ${\cal O}_Y^{\,\Bbb C}$-modules.
 This defines a pairing
     \begin{eqnarray*}
    {\cal T}_\ast Y^\nc \otimes _{{\cal O}_Y^{\,\Bbb C}} {\cal T}^\ast Y^\nc
	                                         &  \longrightarrow        & {\cal O}_Y^\nc \\
        \xi\otimes \omega\hspace{3.4em}          & \longmapsto              & \imath(\xi)(\omega)\,,	 	
     \end{eqnarray*}
     that is nondegenerate in the ${\cal T}_\ast Y^\nc$-factor
	  (namely, if $\imath(\xi)(\omega)= \imath(\xi^\prime)(\omega)$ for all $\omega$, then $\xi=\xi^\prime$).
 This pairing coincides with the fundamental ${\cal O}_Y^{\,\Bbb C}$-bilinear pairing
     \begin{eqnarray*}
    {\cal T}_\ast Y^\nc \otimes _{{\cal O}_Y^{\,\Bbb C}} {\cal T}^\ast Y^\nc
	                                         &  \longrightarrow        & {\cal O}_Y^\nc \\
        \xi\otimes f_1 (df_2) f_3\hspace{.4ex}           & \longmapsto              & f_1 (\xi f_2) f_3\;.
     \end{eqnarray*}
\end{lemma}

\medskip

\begin{remark} $[$not dual  as ${\cal O}_Y^{\,\Bbb C}$-modules$]$\; {\rm
 Caution that, as an ${\cal O}_Y^{\,\Bbb C}$-module,
   the rank of ${\cal T}^\ast Y^\nc$ is greater than the rank of ${\cal T}_\ast Y^\nc$ by $l^2-1$,
   where $l=\rank {\cal F}$ as an ${\cal O}_Y^{\,\Bbb C}$-module.
 Hence, the two are not dual  as ${\cal O}_Y^{\,\Bbb C}$-modules.
 Cf.$\hspace{.7ex}$[L-Y4: Lemma 4.2.3] (D(11.1)).
   }\end{remark}

\medskip

\begin{remark} $[$differential graded algebra associated to ${\cal O}_Y^\nc]$\; {\rm
 For conceptual completion, it should be noted that the exterior differential map
  $d:{\cal O}_Y^\nc\rightarrow {\cal T}^\ast Y^\nc$ extends to a {\it differential graded bi-${\cal O}_Y^\nc$-algebra}
     (in short, {\it DG-algebra})
	$(\Omega^\tinybullet_{Y^{nc}}, d)$ with $d: \Omega^0_{Y^{nc}}\rightarrow \Omega^1_{Y^{nc}}$
    coinciding with the exterior differential map.
 Again, the center-commutativity relation induces a DG-algebra monomorphism
  $\Omega^\tinybullet_Y\otimes_{{\cal O}_Y}\!{\cal O}_Y^\nc
       \hookrightarrow \Omega^\tinybullet_{Y^{nc}}$
  that extends the inclusion
     ${\cal T}^\ast Y\otimes_{{\cal O}_Y}\!{\cal O}_Y^\nc
	      \hookrightarrow {\cal T}^\ast Y^\nc$ in Lemma~3.2.4.
 Cf.~L-Y4: Sec.~4, theme: {\sl Differential graded algebra associated to ring with center a $C^k$-ring}] (D(11.1)).
}\end{remark}

\bigskip

The notion of {\it metric tensor} on a $C^\infty$ Azumaya manifold, adapted from [S\'{e}1] and [S\'{e}2], is given in
 [L-Y4: Sec.~4.3];
see also [Mad:  Sec.$\hspace{.7ex}$3.4] and the more recent [B-M1], [B-M2] for related discussions.
We recall the definition with an updating remark:

\bigskip

\begin{definition} {\bf [metric tensor].} {\rm
 (Cf.$\hspace{.7ex}$[S\'{e}1: Definition 4.2.1], [S\'{e}2: Sec.~3.1.1: Definition$\hspace{.7ex}$3.1];
           [L-Y4: Definition$\hspace{.7ex}$4.3.1] (D(11.1)).)
 A {\it metric tensor} on the Azumaya $C^\infty$-manifold $Y^\nc$
  is a nondegenerate symmetric ${\cal O}_Y^{\,\Bbb C}$-bilinear map
  $$
    g\;:\; {\cal T}_{\ast}Y^\nc
	            \otimes_{{\cal O}_Y^{\,\Bbb C}}{\cal T}_{\ast}Y^\nc\;
				\longrightarrow\; {\cal O}_Y^{\,\Bbb C}\,.
  $$
 Here, we say that $g$ is {\it nondegenerate}  if
  \begin{itemize}
   \item[(1)]
    the induced ${\cal O}_Y^{\,\Bbb C}$-module-homomorphism
	  $\hat{g}:{\cal T}_{\ast}Y^\nc \rightarrow  {\cal T}^{\ast}Y^\nc$
	  over the built-in ${\cal O}_Y^{\,\Bbb C}\hookrightarrow {\cal O}_Y^\nc$,
	 which sends a local section $\xi$ to the functional $\,\cdot\,\mapsto\,g(\xi,\,\cdot\,)$,
	is injective and
    	
   \item[(2)]
    ${\cal T}^{\ast}Y^\nc$ is generated by $\Image\hat{g}$
	both as a left-${\cal O}_Y^\nc$-module and as a right-${\cal O}_Y^\nc$-module.
 \end{itemize}
}\end{definition}

\medskip

\begin{remark} $[$work of Beggs \& Majid\,$]$\; {\rm
 The more recent work [B-M1], [B-M2] of
   {\sl Edwin$\hspace{.7ex}$Beggs} and {\sl Shahn Majid}
  introduce also the notion of metric
  (or `{\it quantum metric}' in their term) for a noncommutative algebra in a very thorough and systematic manner.
 Their notion of metric can also be adapted to a noncommutative algebra with a (possibly complexifed) $C^\infty$-ring center
   to give another, inequivalent, definition of metric tensors for a $C^\infty$ Azumaya manifold $Y^\nc$.
 Both notions and their respective associated Noncommutative Riemannian Geometry on Azumaya manifolds
   should be studied further.
}\end{remark}

\bigskip

\begin{flushleft}
{\bf General case: $Y^\nc$ a $C^\infty$ Azumaya-type noncommutative ringed space}
\end{flushleft}
Let $Y^\nc$ be a $C^\infty$ Azumaya-type noncommutative ringed space as given in the the beginning of
 Sec.$\hspace{.7ex}$3.1.
[L-Y4: Sec.$\hspace{.7ex}$4.1: {\sl Differential calculus on noncommutative ring with center a $C^\infty$-ring}]
 (D(11.1))
 applies to ${\cal O}_Y^\nc$  to give
   both the ${\cal O}_Y^{\,\Bbb C}$-module of {\it derivations} of ${\cal O}_Y^\nc$
    (i.e.$\hspace{.7ex}$ the {\it tangent sheaf} of $Y^\nc$, denoted ${\cal T}_\ast Y^\nc$)
  and the bi-${\cal O}_Y^\nc$-module of {\it differentials} of ${\cal O}_Y^\nc$
   (i.e.$\hspace{.7ex}$ the {\it cotangent sheaf} of $Y^\nc$, denoted ${\cal T}^\ast Y^\nc$).
The universal property of ${\cal T}^\ast Y^\nc$ implies that:

\bigskip

\begin{lemma} {\bf [tangent sheaf, cotangent sheaf, pairing]}\;
 Lemma~3.2.7 remains to hold for a $C^\infty$ Azumaya-type noncommutative ringed space $Y^\nc$.
\end{lemma}
 
\bigskip
 
However,  ${\cal T}_\ast Y^\nc$ in general is no longer a locally free ${\cal O}_Y^{\,\Bbb C}$-module and
  thus the notion of metric tensor takes some care as to what the nondegeneracy condition should be imposed
  over the locus $Y_\scriptsizesing \cup {\,_{\cal P}^{\:\ge 1}}Y$,
  where the rank of the ${\cal O}_Y^{\,\Bbb C}$-module ${\cal T}_\ast Y^\nc$ may jump.
We'll leave this technical issue for the future and be contented with the following weakest form for the moment:

\bigskip

\begin{definition} {\bf [(weakest) metric tensor].} {\rm
 A {\it (weakest) metric tensor} on the $C^\infty$ Azumaya-type noncommutative ringed space $Y^\nc$
  is a symmetric ${\cal O}_Y^{\,\Bbb C}$-bilinear map
  $$
    g\;:\; {\cal T}_{\ast}Y^\nc
	            \otimes_{{\cal O}_Y^{\,\Bbb C}}{\cal T}_{\ast}Y^\nc\;
				\longrightarrow\; {\cal O}_Y^{\,\Bbb C}
  $$
  such that
  its restriction to over $Y-(Y_\scriptsizesing \cup {\,_{\cal P}^{\:\ge 1}}Y)$ is a metric in the sense of
  Definition~3.2.10.
}\end{definition}

\medskip

\begin{motif-question-project} {\bf [noncommutative geometry on $Y^\nc$]}\; {\rm
  On $Y$, when smooth, there are Riemannian geometry, complex geometry, symplectic geometry,   K\"{a}hler geometry,
   calibrated geometry for studying.
  When $Y$ has Gorenstein isolated singularities,
    a noncommutative crepant resolution $Y^\nc$ of $Y$ is meant to be a noncommutative smooth resolution of the singular $Y$.
 Thus,  (cf.\ [B-M1], [DV-K-M], [DV-M], [Giu], [Mas], [M-M-M], [S\'{e}1],~[S\'{e}2])
   \begin{itemize}
    \item[{\bf Q.}] {\it
	 Is there a good noncommutative version of
	    Riemannian geometry, complex geometry, symplectic geometry,   K\"{a}hler geometry, and calibrated geometry on $Y^\nc$?}
   \end{itemize}
   In particular, one would need a version of complex geometry, symplectic geometry, and calibrated geometry to bring in the notion of
    {\it holomorphic cycles, Lagrangian or special Lagrangian submanifolds, cycles associated to a calibration}  to $Y^\nc$
	for studying dynamical open D-branes on $Y^\nc$, in which case, the underlying topology $X$ of the D-brane-world-volume
	  will be a manifold with corners and the map $\varphi$ will take $\partial X$ to these special cycles/submanifolds in $Y^\nc$.	
}\end{motif-question-project}

\bigskip

\subsection{An action  functional for maps $\varphi$, given a connection $\nabla$ and metrics}

With the preparations made in Sec.$\hspace{.7ex}$3.1 and Sec.$\hspace{.7ex}3.2$,
  given a connection $\nabla$ on the fundamental module ${\cal E}$ over the D-brane world-volume
    with underlying topology a $C^\infty$ manifold $X$,
   a metric tensor
      $h: {\cal T}_\ast X\otimes_{{\cal O}_X}\!{\cal T}_\ast X \rightarrow {\cal O}_X$ on $X$,     and
   a metric tensor
      $g: {\cal T}_\ast Y^\nc\otimes_{{\cal O}_Y^{\,\Bbb C}}{\cal T}_\ast Y^\nc
	     \rightarrow {\cal O}_Y^{\,\Bbb C}$
     on $Y^\nc$,
 we are now ready to construct an action functional for maps
   $\varphi: X^{\!A\!z}=(X;{\cal O }_X^{A\!z}:= \Endsheaf_{{\cal O}_X^{\,\Bbb C}}({\cal E}))
      \rightarrow Y^\nc:=(Y, {\cal O}_Y^\nc)$
	 -- defined contravariantly by a ${\Bbb C}$-algebra-homomorphism
     $$
       \varphi^\sharp\;:\; C^\infty(Y^\nc)\;\longrightarrow\;  C^\infty(X^{\!A\!z})
     $$	
	  that restricts to a $C^\infty$-ring-homomorphism $C^\infty(Y)\rightarrow \varphi^\sharp(C^\infty(Y))$ --
  that mimics the Theory of Bosonic Strings (e.g.~[B-B-S], [G-S-W], [Pol4])   and
        the Theory of Harmonic Maps (e.g.~[E-L], [E-S], [L-W])
 and generalizes part of [L-Y6] (D(13.3)).

The map $\varphi: X^{\!A\!z}\rightarrow Y^\nc$
 induces a (left) ${\cal O}_X^{A\!z}$-module-homomorphism
  $$
     d\varphi\; :\;  {\cal O}_X^{A\!z}\otimes_{\varphi^\sharp,\, {\cal O}_Y^{nc}}\!{\cal T}^\ast Y^\nc\;
                                 \longrightarrow\;  {\cal T}^\ast X^{\!A\!z} \,,
 	\hspace{2em}			
    1\otimes (f_1df_2)= \varphi^\sharp(f_1)\otimes df_2
		 \longmapsto \varphi^\sharp(f_1) d\varphi^\sharp(f_2)\,.
 $$
I.e.$\hspace{.7ex}\varphi$ induces a
   $d\varphi\in \Hom_{{\cal O}_X^{A\!z}}(
     {\cal O}_X^{A\!z}\otimes_{\varphi^\sharp,\, {\cal O}_Y^{nc}}\!{\cal T}^\ast Y^\nc ,
    {\cal T}^\ast X^{\!A\!z})$.
Since $X^{\!A\!z}$ is an Azumaya manifold,
 the connection $\nabla$ on ${\cal E}$ induces a bi-${\cal O}_X^{A\!z}$-module epimorphism
  $\pi^\nabla: {\cal T}^\ast X^{\!A\!z}\rightarrow {\cal O}_X^{A\!z}\otimes_{{\cal O}_X}\!{\cal T}^\ast X$,
  as in Lemma~3.2.4
After the composition of $d\varphi$ and $\pi^{\nabla}$, one obtains
 $$
   \nabla\varphi\; :=\;   \pi^\nabla\!\circ  d\varphi \;
   \in\; 	  	
   \Hom_{{\cal O}_X^{A\!z}}(
      {\cal O}_X^{A\!z}\otimes_{\varphi^\sharp,\, {\cal O}_Y^{nc}}\!{\cal T}^\ast Y^\nc ,
     {\cal O}_X^{A\!z}\otimes_{{\cal O}_X}\!{\cal T}^\ast X	  )\,.
 $$

The canonical filling of the diagram,
 i.e.$\hspace{.7ex}$the {\it pullback}  by $\varphi$ or equivalently the {\it pushout} by $\varphi^\sharp$
 $$
  \xymatrix @R=3ex{
    {\cal T}^\ast Y^\nc \ar[d]_-\alpha    &    \\
	    {\cal O}_Y^\nc  \ar[r]^-{\varphi^\sharp}     & {\cal O}_X^{A\!z}
   }
      \hspace{3em}\raisebox{-1.6em}{$\Longrightarrow$}\hspace{3em}
   \xymatrix @R=3ex{
    {\cal T}^\ast Y^\nc \ar[d]_-{\alpha}
	       & {\cal O}_X^{A\!z}  \otimes_{\varphi^\sharp,\,{\cal O}_Y^{nc}}\!{\cal T}^\ast Y^\nc
		          \ar[d]^-{\varphi^\ast \alpha}\\
	{\cal O}_Y^\nc \ar[r]^-{\varphi^\sharp}      & \:\:\:\:{\cal O}_X^{A\!z}
   }
 $$
 induces a natural (left) ${\cal O}_X^{A\!z}$-module-homomorphism
     %
     %
 $$
   \xymatrix @R=3ex{
  {\cal O}_X^{A\!z}\otimes_{\varphi^\sharp,\,{\cal O}_Y^{\,\Bbb C}}
         \Homsheaf_{bi\mbox{-}{\cal O}_Y^{nc}}({\cal T}^\ast Y^\nc, {\cal O}_Y^\nc)
	  \ar[rr]
	 && \Homsheaf_{{\cal O}_X^{A\!z}}(
	        {\cal O}_X^{A\!z}\otimes_{\varphi^\sharp,\,{\cal O}_Y^{nc}}\!{\cal T}^\ast Y^\nc\,,\,
	         {\cal O}_X^{A\!z}) \\
    {\cal O}_X^{A\!z}\otimes_{\varphi^\sharp,\,{\cal O}_Y^{\,\Bbb C}}\!{\cal T}_\ast Y^\nc
	   \ar[u]^-\wr
      && \hspace{3.2em}
	         \raisebox{2.4ex}{\hspace{6em}..}\hspace{-4em}
	       ({\cal O}_X^{A\!z}
		         \otimes_{\varphi^\sharp,\, {\cal O}_Y^{nc}}\!{\cal T}^\ast Y^\nc)^\vee
		      \ar @{=}[u]
	        \hspace{2em}\,.\;\;\hspace{.6em}
    }		
 $$
Note that when ${\cal T}^\ast Y^\nc$ is restricted to over $Y-(Y_\scriptsizesing \cup  {_{\cal P}^{\,\ge 1}}Y)$,
 the horizontal arrow is an isomorphism.
One has also a natural (left) ${\cal O}_X^{A\!z}$-module-homomorphism
  $${\small
  \begin{array}{rcl}
   ({\cal O}_X^{A\!z}\otimes_{\varphi^\sharp,\,{\cal O}_Y^{nc}}\!{\cal T}^\ast Y^\nc)^{\vee_L}
      \otimes_{{\cal O}_X^{\,\Bbb C}} ({\cal O}_X^{A\!z}\otimes_{{\cal O}_X}\!{\cal T}^\ast X)
	        & \longrightarrow  &
	    \Homsheaf_{{\cal O}_X^{A\!z}}(
      {\cal O}_X^{A\!z}\otimes_{\varphi^\sharp,\, {\cal O}_Y^{nc}}\!{\cal T}^\ast Y^\nc\, , \,
     {\cal O}_X^{A\!z}\otimes_{{\cal O}_X}\!{\cal T}^\ast X	  )\,, \\[1.2ex]
	    \alpha\otimes \beta    \hspace{7.2em}
	     & \longmapsto  &
	   \hspace{11.2em}   (\: \tinybullet\;  \mapsto\;  \alpha(\tinybullet)\cdot \beta\, ) \hspace{2.3em}\,,
  \end{array}
  }$$
 which is also an isomorphism
  when ${\cal T}^\ast Y^\nc$ is restricted to over $Y-(Y_\scriptsizesing \cup  {_{\cal P}^{\,\ge 1}}Y)$.
Combining these homomorphisms, one then obtains a natural (left) ${\cal O}_X^{A\!z}$-module-homomorphism
 $$
   ({\cal O}_X^{A\!z}\otimes_{\varphi^\sharp,\,{\cal O}_Y^{nc}}\!{\cal T}_\ast Y^\nc)
      \otimes_{{\cal O}_X^{\,\Bbb C}} ({\cal O}_X^{A\!z}\otimes_{{\cal O}_X}\!{\cal T}^\ast X)\;
	     \longrightarrow\; 	
	    \Homsheaf_{{\cal O}_X^{A\!z}}(
      {\cal O}_X^{A\!z}\otimes_{\varphi^\sharp,\, {\cal O}_Y^{nc}}\!{\cal T}^\ast Y^\nc ,
     {\cal O}_X^{A\!z}\otimes_{{\cal O}_X}\!{\cal T}^\ast X	  )\,,
 $$
 which is an isomorphism
  when ${\cal T}^\ast Y^\nc$ is restricted to over $Y-(Y_\scriptsizesing \cup  {_{\cal P}^{\,\ge 1}}Y)$.

The metric $h$ on ${\cal T}_\ast X$ induces a metric, also denoted by $h$, on ${\cal T}^\ast X$;
 together with the metric $g$ on ${\cal T}_\ast Y^\nc$,
 they define a pairing $\langle\,\tinybullet\,,\,\tinybullet\,\rangle_{(h,g)}$ on
 $({\cal O}_X^{A\!z}\otimes_{\varphi^\sharp,\,{\cal O}_Y^{nc}}\!{\cal T}_\ast Y^\nc)
      \otimes_{{\cal O}_X^{\,\Bbb C}} ({\cal O}_X^{A\!z}\otimes_{{\cal O}_X}\!{\cal T}^\ast X)$
 and hence a {\it rational pairing}, still denoted by $\langle\,\tinybullet\,,\,\tinybullet\,\rangle_{(h,g)}$,  on
 $\Hom_{{\cal O}_X^{A\!z}}(
      {\cal O}_X^{A\!z}\otimes_{\varphi^\sharp,\, {\cal O}_Y^{nc}}\!{\cal T}^\ast Y^\nc ,
     {\cal O}_X^{A\!z}\otimes_{{\cal O}_X}\!{\cal T}^\ast X	  )$
	that takes value in ${\cal O}_{U_\varphi}^{A\!z}$
	for some open set $U_\varphi\subset X$, depending on $\varphi$.
  
\bigskip

\begin{definition} {\bf [energy density; standard action functional for $\varphi$]}\;
{\rm
 Define the {\it $(\nabla$-adjusted, kinetic$)$ energy density} for $\varphi: X^{\!A\!z}\rightarrow Y^\nc$ to be
   $$
      {\frak e}^{(\nabla, h; g)} (\varphi)\, \vol_h\;
	   :=\;  \mbox{$\frac{1}{2}$}\,
	         T\cdot \Real  \left(\Tr \langle \nabla\varphi, \nabla\varphi \rangle_{(h, g)}\right) \vol_h\,,
   $$
  where
    $T$ is called the {\it tension} of the $D$-brane, a constant that depends only on the dimension of $X$;
    $\Real$ indicates the real part;  and $\vol_h$ is the volume form on $X$ associated to the metric $h$.
The {\it standard action functional} $S^{(\nabla, h; g)}$  for $\varphi$, given $(\nabla, h; g)$, is then
 defined to be the {\it $(\nabla$-adjusted, kinetic$)$ energy} $E^{(\nabla, h; g)}(\varphi)$ of $\varphi$:
 $$
    S^{(\nabla, h; g)}(\varphi)\;
	 :=\; E^{(\nabla, h; g)}(\varphi)\;
	 :=\; \int_X {\frak e}^{(\nabla, h; g)}(\varphi)\,\vol_h\;
	  =\; \frac{1}{2}\,T\, \int_X \Real \left(\Tr
	                      \langle \nabla\varphi, \nabla\varphi \rangle_{(h, g)}\right) \vol_h\,.
 $$
 (Cf.~[L-Y6: Definition 4.2] D(13.3).)
}\end{definition}
    
\medskip

\begin{remark} $[$Dirac-Born-Infeld action\,$]$\; {\rm
  For the completeness of discussion, one should mention that there is also an analogue of  the Dirac-Born-Infeld action
     for D-branes in the current setting, cf.\ [L-Y5] (D(13.1)); see also [S\'{e}1].
  However, from the lesson of Nambu-Goto string vs.\ Polyakov string,  it is likely the energy functional
    that will be more suitable for quantization.
}\end{remark}

\medskip
 
\begin{motif-question-project} {\bf [new analytic tool tailored to noncommutative geometry]}\; {\rm
 Formally and before quantization, once an action functional is written down, one may follow Classical Mechanics to address
  its equations of motion and their solutions.
 As a lesson from [L-Y5] (D(13.1)), this could be very messy and non-illuminating.
 \begin{itemize}
  \item[{\bf Q.}] {\it
   Can one develop a new analytic tool really tailored to the type of noncommutative geometry under study
     so that the computation is tidier and more insightful?}
 \end{itemize}
}\end{motif-question-project}

\medskip
 
\begin{motif-question-project} {\bf [beyond prelude]}\; {\rm
  Before anything, one should bring in world-volume supersymmetry into the construction.
  From this point on, the nature of the supersymmetric D-brane world-volume theory will be quite world-volume-dimension dependent.
  Dynamical super D-instantons, dynamical super D-strings, dynamical super D$2$-branes, and dynamical super D$3$-branes
    are likely to convey different features of noncommutative mirror symmetry.
 One should point out that the current notes focus only on the map $\varphi: X^{\!A\!z}\rightarrow Y^\nc$
     and leave the connection $\nabla$ on ${\cal E}$ fixed.
 In a complete theory, $\nabla$ is also dynamical and the full action functional for the pair $(\varphi, \nabla)$ contains
    other terms, including the Yang-Mills and the Chern-Simons terms.
 The D-brane world-volume aspect of Noncommutative Mirror Symmetry, as explained in Sec.~0 Introduction,  is likely
   to mix the duality on the gauge field sector (cf.~Higgs phase $\leftrightarrow$ Coulomb phase; symplectic duality) and
              the duality on the nonlinear sigma-model sector (cf.~target-space duality, noncommutative mirror symmetry).
 For an application to Algebraic Geometry and Symplectic Geometry, this should lead in the end to
    a noncommutative version of  enumerative geometry and
    a Gromov-Witten-type theory for a noncommutative ringed space, or a soft noncommutative scheme,
	or any class of noncommutative spaces that can serve as targets of a map from $X^{\!A\!z}$,
    cf.~[L-Y2]	(D(10.1)) and [L-Y3] (D(10.2)). \dotfill \\
 All these remain very challenging and are open for destined readers to explore.
}\end{motif-question-project}

\newpage
\baselineskip 13pt
{\footnotesize

\vspace{1em}

\noindent
chienhao.liu@gmail.com; 
  \\
yau@math.harvard.edu, styau@tsinghua.edu.cn 

\end{document}